\documentclass[10pt,leqno]{article}

\usepackage[francais]{babel}
\usepackage{amsmath,amssymb,amscd}
\usepackage{a4wide}
 
\usepackage[latin1]{inputenc}
\usepackage{entete-lfp}
\usepackage[all]{xy} %%%pour une sortie pdf, il faut enlever dvips
\title{Sur le comptage des fibrés de Hitchin nilpotents}
\author{Pierre-Henri Chaudouard et Gérard Laumon}
\date{}
\begin{document}

\maketitle

\begin{abstract}
Cet article est une contribution à la fois au calcul du nombre de fibrés de Hitchin sur une courbe projective et à l'explicitation de  la partie nilpotente de la formule des traces d'Arthur-Selberg pour une fonction test très simple. Le lien entre les deux questions a été établi dans \cite{scfh}. On décompose cette partie nilpotente en une somme d'intégrales adéliques indexées par les orbites nilpotentes. Pour les orbites de type \og régulières par blocs\fg{}, on explicite complètement ces intégrales en terme de la fonction zêta de la courbe.
\end{abstract}

\renewcommand{\abstractname}{Abstract}
\begin{abstract}
This paper is concerned with two problems. One is to count  Hitchin bundles on a projective curve and the other is to get an explicit formula for the nilpotent part of the Arthur-Selberg  trace formula for a simple test function. The fact that the two problems are in fact related has been noticed in a previous paper  cf. \cite{scfh}. We expand the nilpotent part of the Arthur-Selberg trace formula in a sum of adelic integrals indexed by nilpotent orbits. For \og regular by blocks\fg{} orbits, we get an explicit formula for these integrals in terms of the zeta function of the curve.
\end{abstract}

\tableofcontents

\section{Introduction}

\begin{paragr} Soit $C$ une courbe projective, lisse, géométriquement connexe sur un corps $k$.  Soit $D$ un diviseur sur $C$. Soit $n\in \NN^*$ et $e\in \ZZ$. Par fibré de Hitchin de rang $n$ et degré $e$, on entend la donnée d'un couple $(\ec,\theta)$ où
  \begin{itemize}
  \item $\ec$ est un fibré vectoriel de degré $e$ et rang $n$ ;
  \item $\theta: \ec\to\ec\otimes_{\oc_C}\oc_C(D)$ est un morphisme de $\oc_C$-module.
  \end{itemize}

Lorsque le rang $n$ est premier au degré $e$ et que le diviseur $D$ est canonique, l'espace de  modules des fibrés de Hitchin \emph{stables}, est une variété algébrique quasi-projective et lisse sur $k$ (cf. \cite{Nit}). Un problème fondamental (et complètement ouvert en rang $>4$), est le calcul de ses nombres de Betti. Hausel et Rodriguez-Villegas (cf. \cite{HRV}) ont proposé une conjecture, plutôt sophistiquée, sur ce que devraient être ces nombres de Betti. Précisons que ces auteurs ne travaillent pas directement sur l'espace de modules en question mais sur une certaine variété de caractères, qui est une variété algébrique affine complexe qui lui est difféomorphe lorsque le corps de base $k$ est le corps des nombres complexes. Ils ont su calculer le $E$-polynôme de la  variété de caractères qu'ils considèrent ; sur la base de ce calcul, ils formulent une expression conjecturale pour son polynôme de Hodge mixte et donc, \emph{ipso facto}, pour son polynôme de Poincaré. Par ailleurs, Garcia-Prada, Heinloth et Schmitt (cf. \cite{GPHS}) ont donné un algorithme pour calculer le motif de l'espace de modules des fibrés de Hitchin. Cela leur a permis de vérifier la conjecture de Hausel et Rodriguez-Villegas en rang $4$ et petit genre. Leur approche exploite la décomposition de Byalinicki-Birula sous l'action du groupe multiplicatif et la géométrie des espaces de chaînes. Toutefois, il ne semble pas évident d'obtenir la  conjecture de Hausel et Rodriguez-Villegas à partir de leur méthode. Notons que Mozgovoy, sur la base du travail de Hausel et Rodriguez-Villegas,  a donné une formule conjecturale pour le motif (cf. \cite{Moz}).

Dans cet article, nous suivons une autre stratégie, expliquée dans \cite{scfh},  pour (tenter de) démontrer les conjectures de Hausel et Rodriguez-Villegas.
\end{paragr}

\begin{paragr}[Comptage de points.] --- Supposons toujours que $n$ est premier au degré $e$ et que le diviseur est canonique ou de degré strictement supérieur à $2g_C-2$, où $g_C$ est le genre de la courbe $C$. Dans cette situation, Nitsure a construit l'espace de modules des fibrés de Hitchin stables, qui est encore quasi-projective et lisse, ainsi qu'un morphisme propre, dit de Hitchin, de cette variété  vers un espace affine. Supposons de plus que le corps de base $k$ est un corps fini. Un argument d'homotopie qui utilise une action du groupe multiplicatif montre que la cohomologie de cet espace de modules est pure. Dès lors,  il suffit de calculer le nombre de points de l'espace de modules sur les corps finis pour obtenir ses nombres de Betti.  Dans \cite{scfh}, on explique comment ce nombre de points s'exprime à l'aide d'une intégrale nilpotente, réminiscente de la partie unipotente de la formule des traces d'Arthur. Cette intégrale est une intégrale adélique de la forme suivante
  \begin{equation}
    \label{eq:int-arthurienne}
    \int_{G(F)\back G(\AAA)^e} K_D(g)\, dg
  \end{equation}
  où $F$ est le corps des fonctions de $C$, $\AAA$ l'anneau des adèles de $F$ et $G$ le groupe $GL(n)$. L'ensemble $G(\AAA)^e$ est l'ensemble des éléments adéliques de $G$ de degré $-e$. Le quotient $G(F)\back G(\AAA)^e$ est muni d'une mesure invariante convenablement normalisée pour laquelle le volume est fini.  La fonction $K_D$ est la valeur en $T=0$  de la fonction  suivante 
$$K_{D,T}(g)= \sum_P (-1)^{\dim(a_P^G)} \sum_{\delta\in P(F)\back G(F)}\hat{\tau}_P(H_P(\delta g)-T) K^P_D(\delta g)$$
où $P$ parcourt les sous-groupes paraboliques standard de $G$ et 
$$K^P_D(g)=   \sum_{X\in \nc^{M_P}} \int_{\ngo_P(\AAA)} \mathbf{1}_D(g^{-1}(X+U)g)\,dU.$$
On a noté $P=M_PN_P$ la décomposition de Levi standard de $P$ où $N_P$ est le radical unipotent de $P$, $\nc^{M_P}\subset \mgo_P(F)$ est l'ensemble des éléments nilpotents $F$-rationnels de l'algèbre de Lie $\mgo_P$ de $M_P$ et $dU$ est une mesure de Haar sur l'algèbre de Lie $\ngo_P(\AAA)$ de $N_P(\AAA)$. La fonction $\mathbf{1}_D$ est une fonction caractéristique liée au diviseur $D$. Pour le lecteur peu familié des travaux d'Arthur, les notations $a_P^G$, $\hat{\tau}_P$ et $H_P$ se trouvent aux §§\ref{S:notations2}, \ref{S:notations3} et \ref{S:reseaux}. Soit $(\nc^{G})$ l'ensemble fini des $G(F)$-orbites nilpotentes  dans  l'algèbre de Lie $\ggo(F)$ de $G(F)$. Pour tout $\oc\in (\nc^G)$, on introduit la fonction suivante de la variable $g\in G(\AAA)$
$$ K^{P}_{D,\oc}(g)=\sum_{X\in \mathcal{N}^{M_P}, I_P^G(X)=\oc} \int_{\ngo_P(\AAA)}  \mathbf{1}_D(g^{-1}(X+U)g)\,dU,$$
où la somme porte sur les $X\in \mathcal{N}^{M_P}$ dont l'induite de Lusztig-Spaltenstein selon $P$ (cf. section \ref{sec:induction}) est l'orbite $\oc$. On a évidemment
$$K^{P}_{D}(g)=\sum_{\oc\in (\nc^G)}K^{P}_{D,\oc}(g).$$
On pose aussi 
$$K_{D,T,\oc}(g)= \sum_{P} (-1)^{\dim(a_P^G)} \sum_{\delta\in P(F)\back G(F)} \hat{\tau}_P(H_B(\delta g)-T) K^{P}_{D,\oc}(\delta g),$$
de sorte qu'on a 
$$K_{D,T}(g)=\sum_{\oc\in (\nc^G)}K_{D,T,\oc}(g).$$
Voici le premier résultat de notre article.

\begin{theoreme} (cf. corollaire \ref{cor:approx})
Pour tout paramètre $T$, l'intégrale   $$\int_{G(F)\back G(\AAA)^e}   K_{D,T,\oc}(g)\, dg$$
converge absolument. De plus, la dépendance en le paramètre $T$ de l'intégrale est quasi-polynomiale.
\end{theoreme}

En particulier, on a le développement suivant
$$\int_{G(F)\back G(\AAA)^e}   K_{D,T}(g)\, dg=\sum_{\oc\in (\nc^G)}\int_{G(F)\back G(\AAA)^e}   K_{D,T,\oc}(g)\, dg.$$
Le mérite du théorème précédent est de ramener le problème du comptage des fibrés de Hitchin au calcul orbite par orbite $\oc$ de l'intégrale 
$$\int_{G(F)\back G(\AAA)^e}   K_{D,T,\oc}(g)\, dg$$
en $T=0$. On a une interprétation en terme de comptage de  cette intégrale. C'est l'objet du théorème suivant. On y utilise les notions de cardinal d'un groupoïde et de $T$-semi-stabilité d'un fibré vectoriel pour lesquelles on renvoie simplement à \cite{scfh}.

\begin{theoreme} (cf.  théorème \ref{thm:approx}) \label{thm-intro:approx}
  La partie quasi-polynomiale en $T$ du cardinal du groupoïde des fibrés de Hitchin $(\ec,\theta)$ de rang $n$ et degré $e$, dont le fibré vectoriel sous-jacent $\ec$ est $T$-semi-stable et l'endomorphisme $\theta$ est génériquement dans l'orbite $\oc$ est  donnée par l'intégrale 
  \begin{equation}\label{eq:KOduthm}
    \int_{G(F)\back G(\AAA)^e}   K_{D,T,\oc}(g)\, dg.
  \end{equation}
\end{theoreme}

Pour faire le lien avec le théorème \ref{thm:approx} du corps de l'article, le cardinal du-dit groupoïde s'exprime à l'aide du dictionnaire adèles-fibrés de Weil par l'intégrale suivante :

$$\int_{G(F)\back G(\AAA)^e}   F^G(g,T) \sum_{X\in \oc} \mathbf{1}_D(g^{-1}Xg) \, dg$$
où $F^G(g,T)$ est la fonction caractéristique du compact de $G(F)\back G(\AAA)^e$ formé des $g$ tels que le fibré correspondant est $T$-semi-stable (pour plus de détails sur ces aller-retours entre adèles et fibrés, comptage et intégrales, on renvoie encore une fois le lecteur à \cite{scfh}). En particulier, l'intégrale ci-dessus est évidemment convergente.

Pour certaines orbites $\oc$, il est possible de donner la valeur de l'intégrale \eqref{eq:KOduthm}. C'est l'autre résultat majeur de cet article. Ces  orbites sont de la forme suivante : on fixe un diviseur $d\geq 1$ de $n$ et soit $r=n/d$ ; ce sont les orbites des éléments nilpotents de $\ggo(F)$ dont la décomposition de Jordan possède $d$ blocs de taille $r$. Parmi celles-ci, on trouve l'orbite nulle qui correspond à $d=n$ et l'orbite régulière qui correspond à $d=1$. La réponse obtenue s'exprime en terme du groupe $\Gc'=SL(n,\CC)$ dont l'apparition n'est pas surprenante : c'est le dual complexe de Langlands du groupe $PGL(n)$ et ce dernier intervient naturellement car la notion de stabilité est invariante par tensorisation par un fibré en droites. Soit $\Mc\subset \Gc'$ le sous-groupe de Levi \og standard\fg{} de $\Gc'$, qui stabilise les $r$ sous-espaces de $\CC^n$ de dimension $d$ en somme directe engendrés par les $d$-uplets de vecteurs \og consécutifs\fg{} de la base canonique de $\CC^n$. Pour $\Pc\subset \Gc'$ sous-groupe parabolique standard de $\Gc'$ de facteur de Levi $\Mc$, on définit une fonction rationnelle sur le tore $Z_{\Mc}^0$ (composante neutre du centre de $\Mc$) par
$$\Phi_{C,D}^d(t)=\prod_{\varpi } Z_{C,D}^d(t^{\varpi}),$$
où le produit est pris sur l'ensemble des poids  fondamentaux $\Pc$-dominants et où l'on introduit la fraction rationnelle de la variable formelle $X$
$$Z_{C,D}^d(X)=X^{-d\deg(D)}Z_C(q^{-1}X)Z_C(q^{-2}X)\ldots Z_C(q^{-d}X),$$
où $Z_C$ est la fonction zêta de la courbe $C$ et $q$ est le cardinal du corps de base $k$. On pose aussi pour tout $e\in \ZZ$ et tout $t\in Z_{\Mc}^0$
$$\Psi_{C,D}^{d,e}(t)=\frac{1}{n\cdot r!} \sum_{z\in \mathbf{\mu}_n} \sum_{w\in \mathfrak{S}_r}z^{-de}  \Phi_{C,D}^d( w\cdot( t z)).$$
Ici on identifie $Z_{\Mc}^0$ à $(\CC^\times)^r$ sur lequel agit naturellement le groupe de permutation $\mathfrak{S}_r$ d'ordre $r!$. Le groupe  $\mathbf{\mu}_n$ des racines $n$-ièmes de l'unité s'identifie naturellement au centre de $\Gc'$. \emph{A priori} la fonction rationnelle $\Psi_{C,D}^{d,e}(t)$ a un pôle en $t=1$ sauf si $d=n$ auquel $\Pc=\Gc$ et $\Psi_{C,D}^{d,e}(t)=1$.

\begin{theoreme}(cf. théorème \ref{cor:calcul})\label{thm:intro-calcul}
  \begin{enumerate}
  \item La fonction rationnelle $\Psi_{C,D}^{d,e}(t)$ n'a pas de pôle en $t=1$.
  \item Soit $\oc\subset \ggo(F)$ l'orbite nilpotente des éléments dont la décomposition de Jordan est formée de $d$ blocs de taille $r$. Lorsque \emph{$e$ est premier à $r$}, l'intégrale \eqref{eq:KOduthm} ci-dessus associée à $\oc$ et à $T=0$, est égale à 
$$q^{n(n-d)\deg(D)/2} q^{nd(g_C-1)} Z_C^*(q^{-1})Z_C(q^{-2})\ldots Z_C(q^{-d}) \Psi_{C,D}^{d,e}(1).$$
où $Z_C^*(X)=(1-qX) Z_C(X)$.
  \end{enumerate}
\end{theoreme}
Lorsque $d=n$, l'intégrale  \eqref{eq:KOduthm} est égale au volume de $G(F)\back G(\AAA)^e$, on a $\Psi_{C,D}^{d,e}(1)=1$ et la formule
$$ q^{n^2(g_C-1)} Z_C^*(q^{-1})Z_C(q^{-2})\ldots Z_C(q^{-n})$$
n'est autre que la formule de Siegel pour ce volume.
 \end{paragr}

\begin{paragr} Cet article résout donc partiellement le problème du comptage des fibrés de Hitchin : il manque encore le  calcul des intégrales \eqref{eq:KOduthm} pour les orbites restantes (c'est-à-dire non régulières par blocs). Mentionnons que dans \cite{scfh}, un raffinement des conjectures de Hausel-Rodriguez-Villégas et Mozgovoy est énoncé qui donne conjecturalement la valeur des intégrales \eqref{eq:KOduthm} en $T=0$. Dans \cite{scfh}, il est également vérifié que le théorème \ref{thm:intro-calcul} est compatible en rang $\leq 3$ avec cette conjecture raffinée. 

Comme on l'a compris, notre méthode consiste à expliciter l'intégrale arthurienne \eqref{eq:int-arthurienne} qui est essentiellement la partie nilpotente de la formule des traces pour la fonction $\mathbf{1}_D$. De fait, les méthodes employées ici doivent certainement s'adapter à des fonctions test autres que les fonctions $\mathbf{1}_D$ et également au cas où $F$ est un corps de nombres. Signalons que dans le cas d'un corps de nombres et de l'orbite régulière, une formule dans le même esprit que celle du théorème \ref{thm:intro-calcul} était connue indépendamment de E. Lapid et T. Finis qui en a parlé lors d'un exposé à Orsay.
\end{paragr}

\begin{paragr}[Organisation de l'article.] --- Le c\oe ur de l'article se trouve à la section \ref{sec:asymp} où est prouvé le théorème \ref{thm-intro:approx}. Une fois ce résultat en poche, on peut dans le cas d'une orbite régulière par blocs \og renormaliser\fg{} l'intégrale \eqref{eq:KOduthm}. C'est l'objet de la section \ref{sec:reg-blocs} dans laquelle on trouve la démonstration du théorème \ref{thm:intro-calcul}. L'intégrale apparaît alors comme une valeur spéciale d'une série qu'on a calculée auparavant dans la section \ref{sec:calculI}. Cependant, pour obtenir une expression raisonnablement simple pour cette valeur spéciale, il nous a fallu développer à la section \ref{sec:combi} une variante des résultats combinatoires qu'Arthur utilise dans ses développements sur la formule des traces. Cette section a donc aussi un intérêt propre. Les sections restantes \ref{sec:induction} et \ref{sec:aux} sont consacrées à quelques rappels utiles.
\end{paragr}

\begin{paragr}[Remerciements.] --- Le premier auteur nommé souhaite remercier T. Finis, W. Hoffmann, E. Lapid et J.-L. Waldspurger pour des discussions enrichissantes.   
\end{paragr}

\section{Combinatoire des $(G,M)$-cofamilles}\label{sec:combi}

\begin{paragr}[Introduction.] ---\label{S:notations} Dans cette section, on développe  la notion de  $(G,M)$-cofamille, très proche de celle  de  $(G,M)$-famille d'Arthur (cf. définitions \ref{def:cofam} et \ref{def:cofamper}). Les résultats de cette section seront utilisés ensuite pour obtenir des formules explicites pour certaines intégrales adéliques. Bien que notre article se focalise sur le groupe $GL(n)$, les méthodes employées doivent pouvoir s'adapter à d'autres groupes réductifs. À des fins de référence future, nous avons rédigé cette section  dans le cadre général d'un groupe réductif sur un corps quelconque. Les principaux résultats de cette section sont formulés dans les théorèmes \ref{thm:ppal} et \ref{thm:ppalper}.

  Les notations sont, à quelques variantes près, celles qui se sont imposées depuis les travaux d'Arthur sur la formule des traces (cf. \cite{ar1}, \cite{trace_inv} ou encore \cite{ar-intro}). Nous les rappelons brièvement dans les prochains paragraphes.
\end{paragr}

  \begin{paragr}[Notations.] ---\label{S:notations2}
   Soit $G$ un groupe réductif sur un corps $F$ quelconque. On fixe un sous-tore déployé maximal de $G$ dont on note $T$ le centralisateur. Sauf indication contraire, par \og sous-groupe parabolique de $G$ \fg{}, on entend un sous-groupe parabolique de $G$ au sens usuel qui, de plus, est défini sur $F$ et contient $T$. Un tel sous-groupe parabolique possède alors un unique facteur de Levi défini sur $F$ et contenant $T$. On appelle simplement sous-groupes de Levi de tels sous-groupes de $G$. Pour un sous-groupe $H$ de $G$, on note $\lc^G(H)$, resp. $\fc^G(H)$, l'ensemble des sous-groupes de Levi, resp. sous-groupes paraboliques, de $G$ qui contiennent $H$. Lorsque $H$ est un sous-groupe de Levi, on note $\pc^G(H)$ l'ensemble des sous-groupes paraboliques de $G$ de facteur de Levi $H$. Ces notations valent encore si l'on remplace $G$ par un sous-groupe de Levi de $G$ ou par un sous-groupe parabolique $Q$ de $G$ (dans ce dernier cas, l'exposant $Q$ signifie qu'on se limite à des sous-groupes inclus dans $Q$). Lorsque le contexte est clair, l'omission de l'exposant (resp. du groupe $H$) signifie qu'on considère une situation relative à $G$ (resp. qu'on prend $H=T$). Par exemple, l'ensemble des sous-groupes de Levi de $G$ est noté simplement $\lc$ au lieu de $\lc^G(T)$.
\medskip

Pour tout $P\in \fc$, on a la décomposition de Levi $P=M_PN_P$ où $M_P$ est le facteur de Levi de $P$ qui contient $T$ et $N_P$ est le radical unipotent de $P$. Soit $A_P$ le tore central déployé maximal de $M_P$ et $X^*(A_P)$ le groupe des caractères rationnels de ce tore. Ce dernier est un réseau dans l'espace vectoriel réel
$$a_{P,\RR}^*=X^*(A_P)\otimes_\ZZ \RR.$$
Soit $X^*(M_P)$ le groupe des caractères rationnels de $M_P$ définis sur $F$. Le morphisme de restriction $X^*(M_P)\to X^*(A_P)$ induit un isomorphisme $X^*(M_P)\otimes_\ZZ \RR\to a_{P,\RR}^*$. Comme cet espace ne dépend  que de $M_P$, il sera aussi noté $a_{M_P,\RR}^*$. Pour tous sous-groupes paraboliques  $P\subset Q$, on a un morphisme de restriction $X^*(A_P)\to X^*(A_Q)$ d'où une projection $a_{P,\RR}^*\to a_{Q,\RR}^*$. En utilisant l'inclusion $M_P\subset M_Q$ et dualement la restriction des groupes de caractères correspondant, on a également  une inclusion $a_{Q,\RR}^*\to a_{P,\RR}^*$. On note sans exposant $*$ les espaces vectoriels réels duaux. L'orthogonal de $a_{Q,\RR}$ dans $a_{P,\RR}^*$ est noté indifféremment $a_{P,\RR}^{Q,*}$ ou $a_{M_P,\RR}^{M_Q,*}$. On a des décompositions en somme directe
$$ a_{P,\RR}^{*}= a_{P,\RR}^{Q,*}\oplus a_{Q,\RR}^{*},$$
et de même pour les espaces duaux. Finalement, on réserve la notation sans indice $\RR$ pour les $\CC$-espaces obtenus par extension des scalaires, par exemple
$$ a_{P}^{*}=a_{P,\RR}^{*}\otimes_{\RR}\CC.$$
On a donc une décomposition en $\RR$-espaces $ a_{P}^{*}=a_{P,\RR}^{*}\oplus ia_{P,\RR}^{*}$ où $i\in \CC$ vérifie $i^2=-1$, et on note $\Re(\la)$ et $\Im(\la)$ les parties réelle et imaginaire de $\la\in  a_{P}^{*}$.

Le groupe de Weyl de $(G,T)$ agit sur $a_{P,\RR}^{*}$ et on fixe sur cet espace un produit euclidien invariant par $W$. La décomposition ci-dessus est alors orthogonale. On met sur $ a_{P,\RR}^{*}$ la mesure euclidienne.

Pour tous $P\in \pc(T)$ et $Q\in \fc(P)$, soit  $\Delta_P^Q$, resp. $\Delta_P^{Q,\vee}$, l'ensemble des racines, resp. des coracines, simples de $T$ dans $N_P\cap M_Q$. Ces ensembles forment des bases respectivement de $a_P^{Q,*}$ et $a_P^Q$. Soit $\hat{\Delta}_P^Q$ et $\hat{\Delta}_P^{Q,\vee}$ les bases duales, respectivement bases  de $a_P^Q$  et $a_P^{Q,*}$. Lorsque l'on omet l'exposant $Q$, cela sous-entend qu'on prend $Q=G$.

On obtient des bases  $\Delta_Q$ et $\Delta_Q^\vee$  de $a_Q^{G,*}$ et $a_Q^G$ par projection des ensembles $\Delta_P-\Delta_P^Q$ et $\Delta_P^\vee-\Delta_P^{Q,\vee}$. Comme la notation le suggère, ces bases ne dépendent du choix de $P\in \pc^Q(T)$. Par dualité, on obtient des bases $\hat{\Delta}_Q$ et $\hat{\Delta}_Q^\vee$ de $a_Q^G$  et $a_Q^{G,*}$. Plus généralement, par le même procédé, pour tout $R\in \fc(Q)$, on définit des ensembles $\Delta_Q^R$, $\hat{\Delta}_Q^R$ etc., le cas $R=G$ redonnant la construction précédente. On observera qu'on a par construction $\Delta_Q^R=\Delta_{Q\cap M_R}^{M_R}$ où dans le membre de droite la construction est relative au groupe réductif $M_R$.

Soit $\ZZ(\Delta_P^{Q,\vee})$ le sous-$\ZZ$-module de $a_P^{Q}$ engendré par $\Delta_P^{Q,\vee}$.
\end{paragr}

\begin{paragr}[Fonctions $\tau$ et $\hat{\tau}$.] --- \label{S:notations3}Pour tous sous-groupes paraboliques $P\subset Q$, soit $\tau_P^Q$ et $\hat{\tau}_P^Q$ les fonctions caractéristiques respectives des ensembles 
$$\{H\in a_{T,\RR} \mid \bg\al,H\bd>0 \,\,\forall \, \al\in \Delta_P^Q\}$$
et
$$\{H\in a_{T,\RR} \mid \bg\varpi,H\bd >0 \,\,\forall \, \varpi \in \hat{\Delta}_P^Q\}.$$
Soit  $a_P^{Q,+}\subset a_{P,\RR}^Q$ la chambre de Weyl ouverte définie par la condition $\tau_P^Q=1$. Soit 
$$(a_P^{Q,*})^+ \subset a_{P,\RR}^{Q,*}$$
 la chambre de Weyl ouverte définie par la condition $\bg \la , \al^\vee\bd >0$ pour tout $\al^\vee \in   \Delta_P^{Q,\vee}$. On note  $\overline{a_P^{Q,+}}$ la chambre de Weyl fermée, donc définie par des inégalités larges.
 
Les fonctions $\tau$ et $\hat{\tau}$ vérifient \og les relations d'orthogonalité\fg{} suivantes où $H\in a_T$ et où l'on somme sur les sous-groupes paraboliques $R$ tels que $P\subset R\subset Q$ (cf. \cite{ar-intro} formules (8.10) et (8.11))

\begin{equation}
  \label{eq:orthog1}
\sum_{\{R \mid P\subset R\subset Q\}  } (-1)^{\dim(a_P^R)}\tau_P^R(H)\hat{\tau}_R^Q(H)=\left\lbrace\begin{array}{l} 0 \text{ si } P \subsetneq Q \\ 1 \text{ si } P= Q
  \end{array}\right.
\end{equation}
et

\begin{equation}
  \label{eq:orthog2}
\sum_{\{R \mid P\subset R\subset Q\}  } (-1)^{\dim(a_P^R)}\hat{\tau}_P^R(H)\tau_R^Q(H)=\left\lbrace\begin{array}{l} 0 \text{ si } P \subsetneq Q \\ 1 \text{ si } P= Q
  \end{array}\right. .
\end{equation}
\end{paragr}

\begin{paragr}[Projection.] --- Pour tous sous-groupes paraboliques $P\subset Q$ et tout  $\la\in a_T^*$, soit  $\la_P^Q$ la  projection de $\la$ sur $(a_P^Q)^*$ selon la décomposition 
$a_T^*=a_T^{P,*}\oplus (a_P^Q)^*\oplus a_Q^*$. Cette décomposition ne dépend que des facteurs de Levi de $P$ et $Q$ ; il en est de même pour la projection. L'omission de l'exposant (resp. de l'indice) signifie qu'on prend $Q=G$, resp. $P\in \pc(T)$.
\end{paragr}

\begin{paragr}[Fractions rationnelles $\theta_P$ et $\hat{\theta}_P$.]\label{S:theta} --- Pour tous sous-groupes paraboliques $P\subset Q\subset G$, on introduit les fractions rationnelles de la variable $\la\in a_T^*$
$$\hat{\theta}_P^Q(\la)=  \frac{\hat{v}_P^Q}{ \prod_{\varpi^\vee\in \hat{\Delta}_P^Q} \bg\la,\varpi^\vee \bd } $$
où $\hat{v}_P^Q$ est le covolume dans $a_P^Q$ du réseau engendré par  $\hat{\Delta}_P^Q$, et 
$$\theta_P^Q(\la)=  \frac{v_P^Q}{ \prod_{\al\in  \Delta_P^Q} \bg\la,\al^\vee \bd } $$
où $v_P^Q$ est le covolume dans $a_P^Q$ du réseau engendré par  $\Delta_P^Q$.

\begin{remarque}
  Ces fonctions $\hat{\theta}_P^Q(\la)$ et $\theta_P^Q(\la)$ sont en fait les inverses de celles introduites et notées de la même façon par Arthur (cf. \cite{trace_inv} p. 15). Par ailleurs, elles  ne dépendent que de $\la_P^Q$. On a
  \begin{equation}
    \label{eq:desc-theta}
    \hat{\theta}_P^Q(\la)=\hat{\theta}_{P\cap M_Q}^{M_Q}(\la_P^Q) \  \text{   et     }     \theta_P^Q(\la)=\theta_{P\cap M_Q}^{M_Q}(\la_P^Q).
  \end{equation}
\end{remarque}

À l'instar des fonctions $\tau$ et $\hat{\tau}$, les fonctions $\hat{\theta}_P^Q(\la)$ et $\theta_P^Q(\la)$ vérifient le lemme suivant.

\begin{lemme}\label{lem:tautau}(Arthur)

\begin{equation}
  \label{eq:orthog1-theta}
\sum_{\{R \mid P\subset R\subset Q\}  } (-1)^{\dim(a_P^R)}\hat{\theta}_P^R(\la)  \theta_R^Q(\la)  =\left\lbrace\begin{array}{l} 0 \text{ si } P \subsetneq Q \\ 1 \text{ si } P= Q
  \end{array}\right.
\end{equation}
et

\begin{equation}
  \label{eq:orthog2-theta}
\sum_{\{R \mid P\subset R\subset Q\}  } (-1)^{\dim(a_P^R)}\theta_P^R(\la)\hat{\theta}_R^Q(\la)=\left\lbrace\begin{array}{l} 0 \text{ si } P \subsetneq Q \\ 1 \text{ si } P= Q
  \end{array}\right. .
\end{equation}
\end{lemme}

\begin{preuve}
  La relation \eqref{eq:orthog1-theta} se déduit de \eqref{eq:orthog1} par une transformation de Fourier (cf. la formule (6.1) p.33 de \cite{trace_inv}). De même \eqref{eq:orthog2-theta} se déduit de \eqref{eq:orthog2}.
\end{preuve}

\end{paragr}

\begin{paragr}[Fonctions $c_P^Q(\varphi,\la)$.] ---\label{S:cphi}
  Soit $P_0$ un sous-groupe parabolique et $\varphi$ une fonction sur $a_{P_0}^*$.  Soit $\varphi^P$ la fonction sur $a_{P_0}^*$ définie par
$$\varphi^P(\la)=\varphi(\la^P).$$
Cette fonction ne dépend que de la projection $\la^P$ de $\la$.
Pour tous sous-groupes paraboliques $P_0\subset P\subset Q\subset G$, on pose
\begin{equation}
  \label{eq:cPQ}
  c_P^Q(\varphi,\la)=\sum_{\{R\mid P\subset R\subset Q  \}} (-1)^{\dim(a_R^Q)} \hat{\theta}_P^R(\la) \varphi^R(\la) \theta_R^Q(\la)
\end{equation}

\begin{remarque}
  Cette définition est formellement très proche de la définition (6.3) p. 33 de \cite{trace_inv}. Hormis un signe, la seule différence notable est qu'on considère la fonction $ \varphi(\la^R)$ et non $\varphi(\la_R)$ comme chez Arthur. On prendra garde cependant que nos fonctions $\theta$ sont les inverses de celles d'Arthur. Lorsqu'on prend la fonction $1$ (constante égale à $1$), on a $c_P^Q(1,\la)=0$ sauf si $P=Q$ auquel cas on a $c_Q^Q(1,\la)=1$ (c'est une conséquence immédiate du lemme \ref{lem:tautau}.
\end{remarque}

Dans la suite, en s'inspirant de  \cite{trace_inv} section 6, on développe les propriétés de ces fonctions.

\begin{lemme} \label{lem:Levi}On a 
$$c_P^Q(\varphi,\la)=c_{M_Q\cap P}^{M_Q}(\varphi^Q,\la^Q).$$
\end{lemme}

\begin{preuve} Elle découle immédiatement  de \eqref{eq:desc-theta}.

\end{preuve}

\begin{lemme}Pour tous sous-groupes paraboliques $P_0\subset P\subset Q\subset G$, on a la formule d'inversion suivante 
$$ \varphi^Q(\la) \hat{\theta}_P^Q(\la)=\sum_{\{R\mid P\subset R\subset Q  \}}c_P^R(\varphi,\la)\hat{\theta}_R^Q(\la) $$
  \end{lemme}

  \begin{preuve}
En utilisant la formule \eqref{eq:cPQ}, on obtient que le second membre vaut
$$\sum_{P\subset S\subset R\subset Q} (-1)^{\dim(a_S^R)} \hat{\theta}_P^S(\la) \varphi^S(\la) \theta_S^R(\la) \hat{\theta}_R^Q(\la)$$
$$=\sum_{P\subset S\subset Q} (-1)^{\dim(a_S^Q)} \hat{\theta}_P^S(\la) \varphi^S(\la) \sum_{S\subset R\subset Q}(-1)^{\dim(a_R^Q)} \theta_S^R(\la) \hat{\theta}_R^Q(\la)$$
 Par \eqref{eq:orthog2-theta}, la somme intérieure vaut $0$ si $S\subsetneq Q$ et $1$ si $S=Q$. On trouve donc
$$ \hat{\theta}_P^Q(\la) \varphi^Q(\la)$$
qui est bien le premier membre.

  \end{preuve}

  \begin{theoreme}
    \label{thm:lissite}
Pour tous sous-groupes paraboliques $P_0\subset P\subset Q\subset G$ et toute fonction $\varphi$ holomorphe sur un voisinage de $0$ dans $a_{P_0}^*$, la fonction $c_P^Q(\varphi,\la)$ est également holomorphe  sur un voisinage de $0$ dans $a_{P_0}^*$.
  \end{theoreme}

La démonstration du théorème \ref{thm:lissite} sera donnée au § \ref{S:demolissite}. Auparavant, nous aurons besoin des constructions auxiliaires du paragraphe suivant.
\end{paragr}

\begin{paragr}[Quelques constructions auxiliaires.] --- Dans ce paragraphe, on fixe $P_0\in \fc(T)$ un sous-groupe parabolique et  $\varphi$ une fonction holomorphe sur un voisinage de $0$ dans $a_{P_0}^*$.

Soit $P$ un sous-groupe parabolique contenant $P_0$.  Soit $(\varpi_\al^{\vee})_{\al\in \Delta^G_P}$ la base de $a_P^G$ duale de  $\Delta^G_P$.
  Pour tout $\la\in a^{*}_T$, on a 
$$\la-\la^P-\sum_{\al \in \Delta^G_P} \bg \la, \varpi_\al^{\vee}   \bd \al\in a_G^*.$$
Pour tout sous-groupe parabolique $Q$ de $G$ contenant $P$, on pose
$$\tilde{\la}^{Q/P}=\la^P+\sum_{\al \in \Delta^Q_P} \bg\la, \varpi_\al^{\vee}   \bd \al.$$
Si $\la\in (a^Q_T)^*$, on a $\tilde{\la}^{Q/P}=\la$. L'application $\la \mapsto \tilde{\la}^{Q/P}$ est donc un projecteur de $a_T^*$ sur $(a^Q_T)^*$ qu'on ne confondra pas avec le projecteur $\la \mapsto \la^Q$. On a cependant
$$\tilde{\la}^{P/P}=\la^P.$$
Dualement, on définit un opérateur analogue sur les fonctions par
\begin{equation}
  \label{eq:phitilde}
  \tilde{\varphi}^{Q/P}(\la)=\varphi(\tilde{\la}^{Q/P}).
\end{equation}
On introduit également la fonction
\begin{equation}
  \label{eq:ctilde}
  \tilde{c}^{Q/P}(\varphi,\la)= \hat{\theta}_P^Q(\tilde{\la}^{Q/P}) \sum_{\{R\mid P\subset R\subset Q\}} (-1)^{\dim(a_R^Q)} \tilde{\varphi}^{R/P}(\la).
\end{equation}

\begin{lemme}\label{lem:lissite}
La fonction $\tilde{c}^{Q/P}(\varphi,\la)$ est holomorphe sur un voisinage de $0$ dans $a_{P_0}^*$.
\end{lemme}

\begin{preuve}
Il s'agit de prouver que la somme alternée  
$$\sum_{\{R\mid P\subset R\subset Q\}} (-1)^{\dim(a_R^Q)} \tilde{\varphi}^{R/P}(\la)$$
s'écrit comme le produit d'une fonction holomorphe sur un voisinage de $0$ dans $a_{P_0}^*$ par 
$$\prod_{\al\in \Delta_P^Q} \bg \tilde{\la}^{Q/P}, \varpi_\al^\vee \bd.$$
Par la formule de Taylor, il faut et il suffit de montrer que cette somme alternée s'annule pour les $\la$ tels qu'il existe $\al\in \Delta_P^Q$
\begin{equation}
  \label{eq:alla}
  \bg \tilde{\la}^{Q/P}, \varpi_\al^\vee \bd=0,
\end{equation}
où $(\varpi_\al^\vee)_{\al\in \Delta_P^Q}$ est la base duale de $a_P^Q$ duale de $\Delta_P^Q$. On a 
$$ \tilde{\la}^{Q/P}=\la^P+\sum_{\beta\in \Delta_P^Q} \bg \la , \varpi_\beta^\vee \bd \beta$$
et la condition \eqref{eq:alla} implique 
\begin{equation}
  \label{eq:alla2}
  \bg \la, \varpi_\al^\vee \bd=0,
\end{equation}
Soit $\al\in \Delta_P^Q$ et $\la\in a_T$ tel que  $\bg \tilde{\la}^{Q/P}, \varpi_\al^\vee \bd=0$. Les sous-groupes paraboliques $R$ tels que $P\subset R\subset Q$ vont par paires $(R,R')$ caractérisées par $\Delta_P^{R'}= \Delta_P^{R}\cup \{\al\}$. On a 
$$\tilde{\la}^{R'/P}=\la^P+\sum_{\beta\in \Delta_P^R} \bg \la , \varpi_\beta^\vee \bd \beta+ \bg \la , \varpi_\al^\vee \bd \al.$$
Donc \eqref{eq:alla2} implique qu'on a 
$$\tilde{\la}^{R'/P}=\tilde{\la}^{R/P}$$
et donc 
$$ \tilde{\varphi}^{R'/P}(\la)= \tilde{\varphi}^{R/P}(\la)$$
d'où l'annulation cherchée.
\end{preuve}

\begin{lemme}\label{lem:egalite}
 Pour tous sous-groupes paraboliques $P\subset R\subset Q$, on a 
$$\hat{\theta}_P^R(\tilde{\la}^{R/P}) \hat{\theta}_R^Q(\tilde{\la}^{Q/P})=\hat{\theta}_P^Q(\tilde{\la}^{Q/P}).$$
\end{lemme}

\begin{preuve}
Il résulte des définitions du §\ref{S:theta} et de l'inclusion $\hat{\Delta}_R^Q\subset \hat{\Delta}_P^Q$ que pour prouver l'égalité cherchée il suffit de vérifier l'égalité
\begin{equation}
  \label{eq:egalite}
  \prod_{\varpi^\vee\in \hat{\Delta}_P^R} \bg \tilde{\la}^{R/P},\varpi^\vee \bd=\prod_{\varpi^\vee\in \hat{\Delta}_P^Q-\hat{\Delta}_R^Q} \bg \tilde{\la}^{Q/P},\varpi^\vee \bd.
\end{equation}
L'ensemble  $\hat{\Delta}_P^Q-\hat{\Delta}_R^Q$ se décrit encore comme $\{\varpi^\vee_\al \mid \al \in \Delta_P^R\}$. Comme on a 
$$ \tilde{\la}^{Q/P}=\la^P+\sum_{\al\in \Delta_P^Q} \bg \la , \varpi_\al^\vee \bd \al $$
et
$$ \tilde{\la}^{R/P}=\la^P+\sum_{\al\in \Delta_P^R} \bg \la , \varpi_\al^\vee \bd \al $$
on a 
$$\prod_{\varpi^\vee\in \hat{\Delta}_P^R} \bg \tilde{\la}^{R/P},\varpi^\vee \bd= \prod_{\beta\in \Delta_P^R} \bg \la , \varpi_\al^\vee \bd$$
ce qui donne bien l'égalité \eqref{eq:egalite}.
\end{preuve}

\begin{lemme} \label{lem:invtilde} On a la formule d'inversion
   $$\hat{\theta}_P^Q(\tilde{\la}^{Q/P})\tilde{\varphi}^{Q/P}(\la)=\sum_{\{R\mid P\subset R\subset Q\}} \hat{\theta}_R^Q(\tilde{\la}^{Q/P})
  \tilde{c}^{R/P}(\varphi,\la)$$
où la somme est prise sur les sous-groupes paraboliques $R$.

\end{lemme}

\begin{preuve}Le second membre s'écrit
  \begin{equation}
    \label{eq:form}
    \sum_{\{R\mid P\subset S\subset  R\subset Q\}} (-1)^{\dim(a_S^R)} \tilde{\varphi}^{S/P}(\la) \hat{\theta}_P^R(\tilde{\la}^{R/P}) \hat{\theta}_R^Q(\tilde{\la}^{Q/P})
  \end{equation}
 En utilisant le lemme \ref{lem:egalite} on obtient
$$ \sum_{\{R\mid P\subset S\subset  R\subset Q\}} (-1)^{\dim(a_S^R)} \tilde{\varphi}^{S/P}(\la) \hat{\theta}_P^Q(\tilde{\la}^{Q/P})$$
$$= \hat{\theta}_P^Q(\tilde{\la}^{Q/P}) \sum_{\{R\mid P\subset  S\subset Q\}}  \tilde{\varphi}^{S/P}(\la)\sum_{\{S \mid S\subset  R\subset Q\}}(-1)^{\dim(a_S^R)}.$$
Comme la somme alternée sur $S$ est nulle sauf si $S=Q$, on obtient bien le second membre cherché.

\end{preuve}
\end{paragr}

\begin{paragr}[Démonstration  du théorème \ref{thm:lissite}.] --- \label{S:demolissite} En raisonnant par récurrence sur la dimension de $G$, on peut supposer que l'énoncé est vrai pour les sous-groupes de Levi propres de $G$. En utilisant le lemme \ref{lem:Levi}, on en déduit que la fonction $c_P^Q(\varphi,\la)$ est holomorphe dans un voisinage de $0$ dans $a_{P_0}^*$ pour $Q\subsetneq G$. Il reste donc à prouver que la fonction $c_P^G(\varphi,\la)$ holomorphe au voisinage de $0$. On a $c_G^G(\varphi,\la)=\varphi(\la)$ et, là encore par récurrence, on peut supposer que pour $P\subsetneq Q$  la fonction $c_Q^G(\varphi,\la)$ est holomorphe au voisinage de $0$ pour toute fonction $\varphi$ holomorphe au voisinage de $0$.

Soit $Q$ un  sous-groupe parabolique contenant $P$. Soit $\mu=\la^Q$. Comme on a $\tilde{\mu}^{Q/P}=\mu$, on obtient
$$\hat{\theta}_P^Q(\la) \varphi^Q(\la)= \hat{\theta}_P^Q(\mu) \varphi(\mu)=\hat{\theta}_P^Q(\tilde{\mu}^{Q/P}) \tilde{\varphi}^{Q/P}(\mu).$$
Le lemme \ref{lem:invtilde} implique qu'on a
\begin{eqnarray*}
 \hat{\theta}_P^Q(\la) \varphi^Q(\la)&=& \hat{\theta}_P^Q(\tilde{\mu}^{Q/P})\tilde{\varphi}^{Q/P}(\mu)\\
&=&\sum_{\{R\mid P\subset R\subset Q\}} \hat{\theta}_R^Q(\tilde{\mu}^{Q/P})
  \tilde{c}^{R/P}(\varphi,\mu).\\
&=&\sum_{\{R\mid P\subset R\subset Q\}} \hat{\theta}_R^Q(\la)  \tilde{c}^{R/P}(\varphi,\la^Q).\\
\end{eqnarray*}
Si l'on insère cette dernière expression  dans l'égalité
$$ c_P^G(\varphi,\la)=\sum_{\{Q \mid P\subset Q\subset G \}} (-1)^{\dim(a_Q^G)} \hat{\theta}_P^Q(\la) \varphi^Q(\la) \theta_Q^G(\la),$$
on obtient
\begin{equation}
  \label{eq:devissagec}
   c_P^G(\varphi,\la)=\sum_{\{R,Q \mid P\subset R\subset Q\subset G \}} (-1)^{\dim(a_Q^G)}  \hat{\theta}_R^Q(\la)  \tilde{c}^{R/P}(\varphi,\la^Q)   \theta_Q^G(\la).
 \end{equation}
 Le terme $ \tilde{c}^{P/P}(\varphi,\la^Q)= \tilde{\varphi}^{P/P}(\la^Q) =\varphi^P(\la^Q)=\varphi(\la^P)$ ne dépend pas de $Q$. Par conséquent la contribution de $R=P$ vaut
$$\sum_{\{Q \mid P\subset Q\subset G \}} (-1)^{\dim(a_Q^G)}  \hat{\theta}_P^Q(\la^Q)  \tilde{c}^{P/P}(\varphi,\la^Q)   \theta_Q^G(\la)=\varphi(\la^P) c_P^G(1,\la)=0$$
par le lemme \ref{lem:tautau}.
Pour tout $R$, on peut poser $\psi_R(\la)= \tilde{c}^{R/P}(\varphi,\la)$. Par le lemme \ref{lem:lissite}, c'est une fonction holomorphe au voisinage de $0$. On a alors l'égalité

$$ c_P^G(\varphi,\la)=\sum_{\{R \mid P\subsetneq R\subset  G \}} c_R^G(\psi_R,\la)$$
qui est bien holomorphe au voisinage de $0$ puisque par hypothèse de récurrence chaque $c_R^G(\psi_R,\la)$ l'est.

\end{paragr}

\begin{paragr}[Sous-groupes paraboliques co-adjacents.] --- Soit $M\in \lc$ un sous-groupe de Levi. Introduisons la définition suivante.

  \begin{definition}
    On dit que deux éléments $P_1$ et $P_2$ de $\pc(M)$ sont co-adjacents s'il existe $\al\in \Delta_{P_1}$ tel que 
$$\Delta_{P_1}\cap \Delta_{P_2}=\Delta_{P_1}-\{\al\}.$$
  \end{definition}

  \begin{lemme}\label{lem:lesequi}
Soit $P_1$ et $P_2$ deux éléments distincts  de $\pc(M)$.
Les conditions suivantes sont équivalentes :
\begin{enumerate}
\item Les sous-groupes paraboliques $P_1$ et $P_2$ sont co-adjacents.
\item Il existe  un unique sous-groupe de Levi maximal $L\in \lc(M)$ tel que
  \begin{enumerate}
  \item $P_1\cap L=P_2\cap L$
  \item $\pc(L)=\{LP_1,LP_2\}$.
  \end{enumerate}
\item Il existe un unique copoids $\varpi^\vee\in \hat{\Delta}^\vee_{P_1}\cap(-\hat{\Delta}^\vee_{P_2})$ tel que les ensembles  respectifs des restrictions des éléments de $\hat{\Delta}_{P_1}^\vee-\{\varpi^\vee\}$ et $\hat{\Delta}_{P_2}^\vee-\{-\varpi^\vee\}$ à l'hyperplan de $a_M$ défini par $\varpi^\vee$ coïncident.
\end{enumerate}
\end{lemme}

\begin{preuve}
  Montrons que 1 implique 2. Soit $Q_1$ et $Q_2$ les sous-groupes paraboliques maximaux $Q_1$ et $Q_2$ qui contiennent respectivement $P_1$ et $P_2$ et qui sont déterminés par la condition 
  
$$\Delta_{P_1}^{Q_1}= \Delta_{P_1}\cap \Delta_{P_2}= \Delta_{P_2}^{Q_2}.$$

  Les groupes  $Q_1$ et $Q_2$ ont le même facteur de  Levi standard qu'on note  $L$. Ce Levi $L$ est maximal. Par le lemme \ref{lem:fibre} assertion 1 ci-dessous, les groupes $Q_1$ et $Q_2$ sont distincts donc nécessairement opposés au sens où $Q_1\cap Q_2=L$. On a par construction $P_1\cap L=P_2\cap L$. Le lemme \ref{lem:fibre}, assertion 2, montre qu'on a $P_1=(P_1\cap L) N_{Q_1}$ d'où $Q_1=LN_{Q_1}=LP_1$ de même pour $P_2$. L'unicité de $L$ résulte du fait qu'on a 
$$\Delta_{P_1}^{LP_1}=\Delta^L_{L\cap P_1}=\Delta^L_{L\cap P_2}\subset \Delta_{P_1}\cap \Delta_{P_2}$$
et pour des raisons de cardinalité cette inclusion est une égalité.

Montrons que 2 implique 3. Soit $Q_i=LP_i$ pour $i=1,2$. Les sous-groupes paraboliques maximaux et distincts $Q_1$ et $Q_2$ ont même facteur de Levi $L$ : ils sont donc opposés. On a donc l'égalité de singleton $\Delta_{Q_1}=-\Delta_{Q_2}$ et dualement $\hat{\Delta}_{Q_1}^\vee=-\hat{\Delta}_{Q_2}^\vee$. Soit  $\varpi^\vee_0 \in \hat{\Delta}_{P_1}^\vee$ tel que  $\hat{\Delta}_{Q_1}^\vee=\{\varpi^\vee_0  \}$.

Pour $i=1,2$, on a
\begin{equation}
    \label{eq:deltaPi}\hat{\Delta}_{Q_i}^\vee= \{  \varpi^\vee \in \hat{\Delta}_{P_1}^\vee \mid   \bg \al,\varpi^\vee\bd =0 \text{ } \forall \al \in \Delta_{P_i}^{Q_i}\}
  \end{equation}
de sorte que $\varpi^\vee_0\in \hat{\Delta}^\vee_{P_1}\cap(-\hat{\Delta}^\vee_{P_2})$. L'hyperplan de $a_M$ défini par  $\varpi^\vee_0$ est le sous-espace $a_M^L$. La projection de $\hat{\Delta}_{P_i}^\vee-\hat{\Delta}_{Q_i}^\vee$ sur $a_M^{L,*}$ n'est autre que la base duale de la base $\Delta_{P_i}^{Q_i}$. Comme cette dernière ne dépend pas de $i$ (c'est simplement $\Delta_{P_i\cap L}^L$), la projection en question ne dépend pas non plus de $i$. L'unicité sera démontrée un peu plus bas.

Enfin 3 implique 1. Soit
$$\Delta_{P_i}'=\{\al \in \Delta_{P_i} \mid  \bg \al, \varpi^\vee\bd=0 \}.$$
C'est un sous-ensemble de $\Delta_{P_i}$ de complémentaire un singleton.
Cet ensemble de racines  est une base de l'hyperplan $\bg \la, \varpi^\vee\bd=0$, duale des restrictions des éléments de $\hat{\Delta}_{P_i}-\{\pm\varpi^\vee\}$. Comme ces restrictions ne dépendent pas de $i$, il en est de même de $\Delta_{P_i}'$. L'assertion 1 s'en déduit.

Revenons à la question de l'unicité dans l'assertion 3. Supposons qu'on ait un autre copoids $\pi^\vee\not=\varpi^\vee$ qui satisfait à la partie existence de l'assertion 3. Dans ce cas, on peut lui associer comme ci-dessus un sous-ensemble  $\Delta_{P_i}''\subset \Delta_{P_i}$ distinct de $\Delta_{P_i}'$ et qui ne dépend de $i$. On a donc $\Delta_{P_i}=\Delta_{P_i}'\cup \Delta_{P_i}''$ ne dépend pas de $i$ d'où $\hat{\Delta}_{P_1}^\vee=\hat{\Delta}_{P_2}^\vee$ ce qui n'est pas possible sous l'hypothèse $\varpi^\vee \in \hat{\Delta}^\vee_{P_1}\cap(-\hat{\Delta}^\vee_{P_2})$.

\end{preuve}

Dans la démonstration précédente, on a utilisé le lemme suivant dont on laisse la preuve au lecteur.

\begin{lemme}\label{lem:fibre}
  Soit $L\in \lc(M)$. Alors
  \begin{enumerate}
  \item Les ensembles $\pc^Q(M)$, pour $Q\in \pc(L)$ sont deux-à-deux disjoints.
  \item Pour tout $Q\in \pc(L)$, l'application
$$P\mapsto P\cap L$$
induit une bijection de $\pc^Q(M)$ sur $\pc^L(M)$ d'inverse $P'\mapsto P'N_Q$.
  \end{enumerate}
\end{lemme}
\end{paragr}

\begin{paragr}[Les $(G,M)$-cofamilles.] --- Arthur a introduit la notion de $(G,M)$-famille : ce sont des familles de fonctions qui sont indexées par les sous-groupes paraboliques de même facteur de Levi $M$ et qui vérifient certaines conditions de recollement pour des sous-groupes paraboliques adjacents (cf. \cite{trace_inv} p.36). Dans la suite, nous aurons besoin de la  notion suivante, qui s'inspire de la définition d'Arthur et qui est relative à la co-adjacence. 

\begin{definition}\label{def:cofam}
    Soit $M$ un sous-groupe de Levi de $G$. Une $(G,M)$-cofamille sur un ouvert $\Omega$ de $a_M^*$ est une famille $(\varphi_P)_{P\in \pc(M)}$ de fonctions holomorphes sur $\Omega$, qui est indexée par les sous-groupes paraboliques $P$ de Levi $M$ et  qui vérifie la condition de recollement suivante :

Pour tout couple $(P_1,P_2)\in \pc(M)^2$ formé d'éléments co-adjacents,  on a 
$$\varphi_{P_1}(\mu)=\varphi_{P_2}(\mu)$$
pour tout $\mu\in \ago^{L,*}_M\cap \Omega$ où $L\in \lc(M)$ est le sous-groupe de Levi $L\in \lc^G(M)$ maximal déterminé par $P_1$ et $P_2$ (cf. assertion 2 du lemme \ref{lem:lesequi}).
\end{definition}

\begin{remarque}
On peut reformuler la condition de recollement ainsi : pour tout couple  $(P_1,P_2)\in \pc(M)^2$ de sous-groupes co-adjacents, la fonction  $\varphi_{P_1}-\varphi_{P_2}$ est divisible par $\bg\la, \varpi^\vee \bd$ dans l'anneau des fonctions holomorphes sur $\Omega$ où $\varpi^\vee$ est le copoids donné par l'assertion 3 du lemme \ref{lem:lesequi}.
\end{remarque}

\begin{lemme}\label{lem:ex} Soit $\varphi$ une fonction holomorphe sur un ouvert $\Omega$ de $\CC$. Pour tout $P\in \pc(M)$, soit $\varphi_P$ la fonction sur $a_M^*$ définie par
$$\varphi_P(\la)=\prod_{\al\in \Delta_P} \varphi(\bg\la, \varpi_\al^\vee \bd).$$
La famille $(\varphi_P)_{P\in \pc(M)}$ est une $(G,M)$-cofamille l'ouvert des $\la\in a_M^*$ tels que $\bg\la, \varpi_\al^\vee \bd\in \Omega$.
\end{lemme}

\begin{preuve}
   Soit $P$ et $P'$ adjacents.  Soit $\varpi^\vee \in \hat{\Delta}_P^\vee \cap (-\hat{\Delta}_{P'}^\vee)$ qui satisfait l'assertion 3 du lemme \ref{lem:lesequi}. Il s'agit de vérifier que les fonctions $\varphi_P$ et $\varphi_{P'}$ coïncident sur l'hyperplan défini par $\varpi^\vee$. C'est évident pour les facteurs associés à  $\pm \varpi^\vee$ qui valent tous deux $\varphi(0)$. Les autres facteurs sont indexés par des coracines qui ont même restriction à l'hyperplan défini par $\varpi^\vee$ et sont donc égaux deux à deux. 
\end{preuve}

\begin{lemme}\label{lem:indpceQ}
Soit $(\varphi_P)_{P\in \pc(M)}$ une $(G,M)$-cofamille. Soit $L\in \lc(M)$ et  $Q\in \pc(L)$.  Pour tout $P\in \pc^L(M)$ et $\la\in (a_M^L)^*$, on pose
$$\varphi^Q_P(\la)=\varphi_{PN_Q}(\la).$$
La famille $\varphi^Q=(\varphi^Q_P)_{P\in \pc^L(M)}$ est une $(L,M)$-cofamille. Elle ne dépend pas du choix de  $Q\in \pc(L)$. On la note $\varphi^L$ dans la suite.
\end{lemme}

\begin{preuve}
Soit $P_1$ et $P_2$ deux éléments de $\pc^L(M)$ co-adjacents. Pour $i\in \{1,2\}$, soit $\tilde{P}_i=P_i N_Q$. On a l'égalité 
$$\Delta_{\tilde{P_1}}-\Delta_{\tilde{P_1}}^Q=\Delta_{\tilde{P_2}}-\Delta_{\tilde{P_2}}^Q.$$
Comme $\Delta_{\tilde{P_i}}^Q=\Delta_{P_i}$, les sous-groupes paraboliques  $\tilde{P}_1$ et $\tilde{P}_2$ sont co-adjacents. Soit $\tilde{S}$ le Levi maximal qu'ils définissent. Alors $S=\tilde{S} \cap L$ est le Levi maximal défini par $P_1$ et $P_2$. Pour tout $\mu\in a_M^{S,*}\subset a_M^{\tilde{S},*}$, on a 
$$\varphi^Q_{P_1}(\mu)=\varphi_{\tilde{P_1}}(\mu)=\varphi_{\tilde{P_2}}(\mu)=\varphi^Q_{P_2}(\mu),$$
les égalité extrêmes résultent des définitions alors que l'égalité du milieu est la condition de recollement pour la cofamille de départ. La famille $\varphi^Q$ est donc une $(L,M)$-cofamille.

Montrons ensuite qu'elle ne dépend pas du choix de $Q$. Pour tout couple $(Q,Q')$ d'éléments de $\pc(L)$, il existe une suite finie d'éléments $Q_1,\ldots,Q_r$ de $\pc(L)$ tels que $Q_1=Q$, $Q_2=Q'$ et $Q_i$ et $Q_{i+1}$ sont adjacents (au sens usuel: les chambres aiguës dans $a_L^G$ associées à $Q_i$ et $Q_{i+1}$ ont un mur en commun). Il suffit donc de traiter le cas où  $Q$ et $Q'$ sont adjacents. Dans ce cas, il existe $R\in \fc(L)$  tel que 
\begin{enumerate}
\item $L$ est un sous-groupe de Levi maximal de $M_R$ ;
\item les sous-groupes $S=M_R\cap Q$ et $S'=M_R\cap Q'$ sont les deux éléments distincts (et donc opposés) de $\pc^{M_R}(L)$ ;
\item $Q=SN_R$ et $Q'=S'N_R$.
\end{enumerate}
Soit $P\in \pc^L(M)$. Les sous-groupes $Q$ et $Q'$ déterminent les éléments $PN_Q$ et $PN_{Q'}$ dans $\pc^G(M)$. Il s'agit de voir l'égalité pour tout $\mu\in a_M^{L,*}$,
\begin{equation}
  \label{eq:avoir}
  \varphi_{PN_Q}(\mu)=\varphi_{PN_{Q'}}(\mu)
\end{equation}

On a aussi $N_Q=N_SN_R$ et $N_{Q'}=N_{S'}N_R$ et donc  $PN_Q=PN_S N_R$ et $PN_{Q'}=PN_{S'} N_R$. Les sous-groupes $PN_S$ et $PN_{S'}$ sont deux éléments de  $\pc^{M_R}(M)$. Par conséquent, l'égalité  \eqref{eq:avoir} est équivalente à l'égalité suivante pour la $(M_R,M)$-cofamille $\varphi^R$ et  $\mu\in a_M^{L,*}$
\begin{equation}
  \label{eq:avoir2}
  \varphi_{PN_S}^R(\mu)=\varphi_{PN_{S'}}^R(\mu).
\end{equation}
Les éléments $PN_S$ et $PN_{S'}$ sont en fait deux éléments co-adjacents de  $\pc^{M_R}(M)$. L'égalité \eqref{eq:avoir2} ne fait donc que traduire une condition de recollement pour $(M_R,M)$-cofamille $\varphi^R$.
\end{preuve}
\end{paragr}

\begin{paragr}[Fonctions associées à une $(G,M)$-cofamille.] --- Arthur a donné un procédé pour construire une fonction lisse à partir d'une $(G,M)$-famille et des fonctions $\theta_P$ (cf. \cite{trace_inv} p.37). Dualement, on a l'énoncé suivant pour une $(G,M)$-cofamille et les fonctions  $\hat{\theta}_P$.
  
\begin{lemme}\label{lem:fctlisse}
  Soit $(\varphi_P)_{P\in \pc(M)}$ une $(G,M)$-cofamille sur un ouvert $\Omega$ de $a_M^*$.  La fonction $$\varphi_M(\la)=\sum_{P\in \pc(M)} \hat{\theta}_P^G(\la) \varphi_P(\la)$$
est holomorphe sur $\Omega$.
\end{lemme}

\begin{preuve}
On utilise de manière répétée le lemme suivant, qui se démontre à l'aide de la formule de Taylor. 

\begin{lemme}\label{lem:Taylor} Soit $\Phi$ une fonction méromorphe sur un voisinage de $z_0\in \CC^n$ et $l_1,\ldots,l_r$ des fonctions affines sur $\CC^n$ telles que les hyperplans affines $H_i$ définis par  $l_i=0$ soient deux à deux distincts et  la fonction
$$\Psi : z\mapsto l_1(z)\ldots l_r(z) \Phi(z)$$
se prolonge en une fonction holomorphe sur un un voisinage de $z_0\in \CC^n$.

Supposons que,  pour tout indice $i$ tel que $H_i$ contient $z_0$,  la restriction de $\Psi$ à $H_i$ est nulle dans un voisinage de $z_0$. Alors $\Phi$ est holomorphe au voisinage de $z_0$.
\end{lemme}

Soit $\Pi$   un système de représentants du quotient de $\cup_{P\in \pc(M)} \hat{\Delta}_P^\vee$ par la  relation de colinéarité. Introduisons les fonctions polynomiales sur $a_M^*$
$$F_\Pi(\la)=\prod_{\varpi^\vee \in \Pi} \bg \la,\varpi^\vee\bd$$
et pour tout $P\in \pc(M)$
$$\hat{\Theta}_{P}(\la)=F_\Pi(\la) \cdot \hat{\theta}_{P}(\la).$$
D'après le lemme \ref{lem:Taylor}, il suffit de voir que la fonction
$$\psi_M(\la)=F_\Pi(\la) \varphi_M(\la)=\sum_{P\in \pc(M)} \hat{\Theta}_{P}(\la) \varphi_P(\la),$$
qui est clairement holomorphe sur $\Omega$,  s'annule sur les hyperplans  d'équation $\bg \la,\varpi^\vee\bd=0$ lorsque $\varpi^\vee$ décrit $\Pi$. Soit  $\varpi^\vee\in \Pi$ et  $\hgo$ l'hyperplan associé. Soit $P\in \pc(M)$ tel que $\varpi^\vee\in  \hat{\Delta}_P^\vee$. Soit $Q$ le  sous-groupe parabolique maximal $Q$ de Levi $L$ défini par $P\subset Q$ et $\hgo=\ago^{L,*}_M$.

Soit $P'\in \pc(M)$. Le polynôme $\hat{\Theta}_{P'}(\la)$ s'annule de manière évidente sur $\hgo$ sauf  si $\varpi^\vee$ est proportionnel à un élément de $\hat{\Delta}_{P'}^\vee$. Dans ce cas, on définit un sous-groupe parabolique $Q'$ contenant $P'$ par la condition
$$\Delta_{P'}^{Q'}=\{\al\in \Delta_P \mid \bg \al,\varpi^\vee \bd=0\}.$$
On a $Q'\in \pc(L)$. Comme $\pc(L)$ est formé de la paire de sous-groupes paraboliques opposés $\{Q,\bar{Q}\}$, l'ensemble des $P'$ tels que $\varpi^\vee$ est proportionnel à un élément de $\hat{\Delta}_{P'}^\vee$ se décrit à l'aide du lemme \ref{lem:fibre}. C'est la réunion disjointe sur $P_0\in \pc^L(M)$ des paires $\{P_0N_Q, P_0N_{Q'}\}$ de sous-groupes paraboliques co-adjacents. Soit  $P_0\in \pc^L(M)$ et $\{P_1,P_2\}$ la paire associée. Pour prouver la nullité de  $\psi_M$ sur $\hgo$, il suffit de prouver celle de 
$$\la \mapsto \big(\hat{\Theta}_{P_1}(\la) \varphi_{P_1}(\la)+\hat{\Theta}_{P_2}(\la) \varphi_{P_2}(\la)\big).$$ 
Par le lemme \ref{lem:lesequi} assertion 3, on a pour $\la\in \Omega\cap \hgo$
$$\hat{\Theta}_{P_1}(\la) \varphi_{P_1}(\la)=-\hat{\Theta}_{P_2}(\la) \varphi_{P_1}(\la).$$
Il suffit donc de prouver la nullité sur $\hgo\cap \Omega$ de 
$$\varphi_{P_1}(\la)-\varphi_{P_2}(\la)$$
qui est précisément la condition de recollement des $(G,M)$-cofamilles.
\end{preuve}

Soit $(\varphi_P)_{P\in \pc(M)}$ une  $(G,M)$-cofamille. Pour tout $P\in \pc(M)$ et tout sous-groupe parabolique$Q$ contenant $P$, soit $\varphi_P^Q$ la fonction sur $a_M^{G,*}$ définie par
$$\varphi_P^Q(\la)=\varphi_P(\la^Q),$$
où, rappelons-le,  $\la^Q=\la^{M_Q}$ est la projection de $\la$ sur $a_M^{Q,*}$.

  \begin{lemme}\label{lem:proc-indpQ}
    Soit $\varphi=(\varphi_P)_{P\in \pc(M)}$ une $(G,M)$-cofamille. Soit $L\in \lc(M)$ et $Q\in \pc(L)$. La fonction de la variable $\la\in a_M^{G,*}$
$$\sum_{P\in \pc^Q(M)} \hat{\theta}_P^Q(\la) \varphi_P^Q(\la)$$
ne dépend que du sous-groupe de Levi $L$. Elle est égale à 
$$\varphi_{M}^{L}(\la^{L}).$$
  \end{lemme}
  
\begin{remarque}
  Dans l'énoncé du lemme   \ref{lem:proc-indpQ} ci-dessus, on note $\varphi^{L}$ la $(L,M)$-cofamille déduite de $\varphi$ par le procédé du lemme  \ref{lem:indpceQ} et $\varphi^{L}_M$ la fonction qui s'en déduit par le lemme \ref{lem:fctlisse}.
  \end{remarque}
  
  \begin{preuve}
Pour chaque terme de la somme considéré on a 
$$ \hat{\theta}_P^Q(\la) \varphi_P^Q(\la)=\hat{\theta}_{P\cap L}^L(\la) \varphi_{P\cap L}^L(\la)$$
par le lemme \ref{lem:indpceQ}. L'application $P\mapsto P\cap L$ induit une bijection de $\pc^Q(M)$ sur $\pc^L(M)$. Le lemme s'en déduit.
    
  \end{preuve}

\end{paragr}

\begin{paragr}[Co-familles et fonctions $c_P$.] --- L'énoncé suivant est un analogue du corollaire 6.4 de \cite{trace_inv}.

\begin{theoreme}\label{thm:ppal}
    Soit $(\varphi_P)_{P\in \pc(M)}$ une $(G,M)$-cofamille. On a l'égalité suivante pour tout $\la\in a_M^*$
$$\varphi_M(\la^G)=\sum_{P\in \pc(M)} c_P^G(\varphi_P,\la).$$
  \end{theoreme}

  \begin{preuve}
On explicite le membre de droite à l'aide de \eqref{eq:cPQ} ce qui donne 

$$    \sum_{P\in \pc(M)} \sum_{P\subset Q\subset G} (-1)^{\dim(a_Q^G)} \hat{\theta}_P^Q(\la) \varphi_P^Q(\la) \theta_Q^G(\la).$$
Pour éviter toute ambiguïté, précisons que $\varphi_P^Q(\la)$ est la fonction définie par
$$\varphi_P^Q(\la)=\varphi_P(\la^Q).$$
En intervertissant les deux sommes, on obtient 
$$\sum_{Q\in \fc^G(M)} (-1)^{\dim(a_Q^G)} \theta_Q^G(\la) \sum_{P\in \pc^Q(M)} \hat{\theta}_P^Q(\la) \varphi_P^Q(\la).$$
Par le lemme \ref{lem:proc-indpQ}, on reconnaît dans la somme intérieure la fonction $\varphi_M^{M_Q}(\la^{M_Q})$. L'expression précédente devient

$$\sum_{Q\in \fc^G(M)} (-1)^{\dim(a_Q^G)} \theta_Q^G(\la)\varphi_{M}^{M_Q}(\la^{M_Q})$$
$$=\sum_{L\in \lc^G(M)} (-1)^{\dim(a_L^G)}\varphi_{M}^{L}(\la^L)\sum_{Q\in \pc(L)} \theta_Q^G(\la)$$
D'après \cite{trace_inv} pp. 36 et 39, preuve du corollaire 6.4, on a 
$$\sum_{Q\in \pc(L)} \theta_Q^G(\la)= \left\lbrace \begin{array}{l} 0 \text{ si } L\subsetneq G \\ 1 \text{ si } L=G \end{array}\right.. $$
Il s'ensuit que dans la somme précédente ne subsiste que le terme $L=G$ qui donne bien $\varphi_M(\la^G)$. 

     \end{preuve}

\end{paragr}

\begin{paragr}[Point $\xi$.] ---
 Soit $\xi \in  a_{T,\RR}$. Pour tous sous-groupes paraboliques $P\subset Q$, on définit un point
$$[\xi]_P^Q\in a_{P,\RR}^Q$$
 de la manière suivante. Soit $\xi_P^Q$ le projeté de $\xi$ sur $a_P^Q$ qu'on écrit dans la base $\Delta_P^{Q,\vee}$ :
$$\xi_P^Q=\sum_{\al^\vee\in \Delta_P^{Q,\vee}} \bg\varpi_{\al^\vee}, \xi_P^Q  \bd \al^\vee.$$
On pose alors 
$$[\xi]_P^Q=\sum_{\al^\vee\in \Delta_P^{Q,\vee}} [\bg\varpi_{\al^\vee}, \xi_P^Q  \bd] \al^\vee,$$
où pour tout $x\in \RR$, le symbole $[x]$ désigne le plus petit entier $i$ tel que $i>x$. 

\begin{lemme}\label{lem:projxi}
Pour tous sous-groupes paraboliques $P\subset R\subset Q$, la projection de  $[\xi]_P^Q$ sur $a_R^Q$ est égale à $[\xi]_R^Q$.
\end{lemme}

\begin{preuve} Pour tout $X\in a_P^Q$, on a 
  $$X=\sum_{\al^\vee\in \Delta_P^{Q,\vee}} \bg\varpi_{\al^\vee}, X \bd \al^\vee.$$
Soit $p_R^Q$ la projection sur $a_R^Q$.
Pour tout $\al^\vee \in \Delta_P^{Q,\vee}-\Delta_P^{R,\vee}$, on a $\varpi_{\al^\vee}\in \hat{\Delta}_R^Q$ et 
$$ \bg\varpi_{\al^\vee}, X  \bd= \bg\varpi_{\al^\vee}, p_R^Q(X)  \bd.$$
 On a donc 
\begin{eqnarray*}
  p_R^Q(X) &=&\sum_{\al^\vee\in \Delta_P^{Q,\vee}-\Delta_P^{R,\vee}} \bg\varpi_{\al^\vee}, X  \bd p_R^Q(\al^\vee)\\
&=&\sum_{\al^\vee\in \Delta_R^{Q,\vee}} \bg\varpi_{\al^\vee}, X  \bd \al^\vee.
\end{eqnarray*}
On en déduit
\begin{eqnarray*}
  p_R^Q([\xi]_P^Q) &=&\sum_{\al^\vee\in \Delta_P^{Q,\vee}-\Delta_P^{R,\vee}} \bg\varpi_{\al^\vee}, [\xi]_P^Q  \bd p_R^Q(\al^\vee)\\
&=& \sum_{\al^\vee\in \Delta_P^{Q,\vee}-\Delta_P^{R,\vee}} [\bg\varpi_{\al^\vee}, \xi_P^Q  \bd] p_R^Q(\al^\vee)\\
  &=& \sum_{\al^\vee\in \Delta_R^{Q,\vee}} [\bg\varpi_{\al^\vee}, \xi_R^Q  \bd] \al^\vee \\
&=& [\xi]_R^Q.
\end{eqnarray*}
\end{preuve}
\end{paragr}

\begin{paragr}[Fonctions $\hat{\vartheta}$ et $\vartheta^\xi$.] --- Dans toute la suite, on fixe $q$ un entier $>1$. Pour tous sous-groupes paraboliques $P\subsetneq Q$, on définit les fonctions méromorphes de la variable $\la\in a_T^*$
\begin{equation}
  \label{eq:vartheta}
  \vartheta_P^{Q,\xi}(\la)=q^{-\bg \la, [\xi]_P^Q\bd} \prod_{\al^\vee \in \Delta_P^{Q,\vee}} \frac{1}{1-q^{-\bg \la,\al^\vee\bd}}
\end{equation}
et 
\begin{equation}
  \label{eq:varthetac}
  \hat{\vartheta}_P^{Q}(\la)=\prod_{\varpi^\vee \in \hat{\Delta}_P^{Q,\vee}} \frac{1 }{1-q^{-\bg \la,\varpi^\vee\bd}}
\end{equation}
On complète ces définitions en posant $\vartheta_P^{P,\xi}(\la)=1$ et $\hat{\vartheta}_P^{P}(\la)=1$.

\begin{lemme}
  Soit $P\subset Q$ des sous-groupes paraboliques. La série  de Fourier de la variable $\la\in a_T^*$
$$\sum_{H\in \ZZ(\Delta_P^{Q,\vee}) }\hat{\tau}_P^Q(H-\xi) q^{-\bg\la, H  \bd },$$
converge pour $\Re(\la_P^Q)\in (a_P^{Q,*})^+$. Sur ce domaine de convergence, sa somme coïncide avec la fonction  $\vartheta_P^{Q,\xi}(\la)$.
\end{lemme}

\begin{preuve}
Comme la famille $\Delta_P^{Q,\vee}$ est duale de $\hat{\Delta}_P^{Q}$, la série est égale au produit sur $\al^\vee \in \Delta_P^{Q,\vee}$ des séries géométriques
$$\sum_{n_\al} q^{-n_\al \bg\la,\al^\vee   \bd }$$
où $n_\al\in \ZZ$ satisfait l'inégalité $n_\al > \bg \varpi_\al, \xi \bd$, ou encore  $n_\al \geq [\bg \varpi_\al, \xi \bd]$. En particulier, elle  converge absolument pour $\Re( \bg\la,\al^\vee\bd)>0$ pour tout $\al^\vee \in \Delta_P^{Q,\vee}$, c'est-à-dire pour  $\la \in  (a_P^{Q,*})^+$.

Ces séries géométriques ont pour somme, suivant $\al^\vee \in \Delta_P^{Q,\vee}$,
$$\frac{q^{- [\bg \varpi_\al, \xi \bd] \bg\la,\al^\vee   \bd   }}{1-q^{-\bg \la,\al^\vee\bd}}.$$
Le résultat s'en déduit vu qu'on a  $\bg \varpi_\al, \xi \bd=\bg \varpi_\al, \xi_P^Q \bd$ et 
$$\sum_{\al^\vee \in \Delta_P^{Q,\vee}} - [\bg \varpi_\al, \xi \bd] \bg\la,\al^\vee   \bd   =-\bg \la, [\xi]_P^Q\bd.$$
\end{preuve}

\begin{lemme} \label{lem:xigeneral}Soit $\xi\in a_{P,\RR}^*$ en position générale, $Q$ un sous-groupe parabolique et $M\in \lc^Q$. On a
$$\sum_{P\in \pc^Q(M)} \vartheta_P^{Q,\xi}(\la)=\left\lbrace \begin{array}{l} 0 \text{ si } M\subsetneq M_Q \\ 1 \text{ si } M=M_Q \end{array}\right..$$
  \end{lemme}

  \begin{preuve}
Le cas $M_Q=M$ est évident. On suppose donc $M\subsetneq M_Q$. Quitte à remplacer $G$ par $M_Q$, il suffit de traiter le cas $Q=G$. La démonstration  repose sur des lemmes combinatoires dus à Langlands et Arthur (cf. \cite{localtrace}  p.22). Soit $\Lambda\in a_{M,\RR}^*$ en position générale et pour  tout $P\in \pc(M)$ soit 
$$\Delta_P^{\Lambda}=\{\al\in \Delta_P\mid \bg\Lambda,\al^\vee \bd <0\}$$
Soit $(\varpi_\al)_{\al \in \Delta_P}$ la base de $a_P^{G,*}$ duale de $\Delta_P^\vee$ et $\chi_P^{\Lambda}$ la fonction caractéristique des $H\in a_{M,\RR}$ tels que $\bg \varpi_\al, H\bd>0$ si $\al\in \Delta_P^{\Lambda}$ et $\bg \varpi_\al, H\bd\leq 0$ si $\al\in \Delta_P-\Delta_P^{\Lambda}$.
D'après Arthur (\emph{loc. cit.}), la fonction de la variable $H\in a_{M,\RR}^G$
$$\sigma_M(H)=\sum_{P\in \pc(M)} (-1)^{|\Delta_P^{\Lambda}|}\chi_P^{\Lambda}(H)$$
ne dépend pas de $\Lambda$ et vérifie
\begin{equation}
  \label{eq:sigmaM}
  \sigma_M(H)= \left\lbrace
  \begin{array}{c}
    0 \text{ si } H\not=0\\
1 \text{ si } H=0
  \end{array}\right. .
\end{equation}

Soit $P_0\in \pc(M)$ défini par la condition $\Lambda\in (a_{P_0}^*)^{+}$.
Pour tout $\la \in a_T^*$ tel que $\Re(\la_P)\in -(a_{P_0}^*)^{+}$ et tout $P\in \pc(M)$, la série 
$$(-1)^{|\Delta_P^{\Lambda}|}\sum_{H\in \ZZ(\Delta_P^\vee)} \chi_P^{\Lambda}(H-\xi)q^{-\bg \la, H\bd}$$
converge absolument. Cette série, qui  est un produit de séries géométriques, a pour  somme 
$$(-1)^{|\Delta_P^{\Lambda}|} \prod_{\al\in \Delta_P^\Lambda} \frac{q^{-[\bg \varpi_\al, \xi\bd ]\bg \la, \al^\vee\bd}}{1-q^{-\bg \la, \al^\vee\bd}}\prod_{ \al\in \Delta_P-\Delta_P^{\Lambda} } \frac{q^{(-[\bg \varpi_\al, \xi\bd ]+1) \bg \la, \al^\vee\bd}}{1-q^{\bg \la, \al^\vee\bd}}.$$

$$=(-1)^{|\Delta_P|}\prod_{\al\in \Delta_P} \frac{q^{-[\bg \varpi_\al, \xi\bd ]\bg \la, \al^\vee\bd}}{1-q^{-\bg \la, \al^\vee\bd}}=(-1)^{|\Delta_P|}\vartheta_P^\xi(\la).$$
Ni le cardinal $|\Delta_P|$ ni le réseau $\ZZ(\Delta_P^\vee)$ ne dépendent du choix de $P\in \pc(M)$.
Il s'ensuit que pour tout $\la \in a_T^*$ tel que $\Re(\la_P)\in -(a_{P_0}^*)^{+}$ on a 
\begin{equation}
  \label{eq:egalitecarfct}
  \sum_{P\in \pc(M)} \vartheta_P^\xi(\la)= (-1)^{|\Delta_{P_0}|} \sum_{H\in \ZZ(\Delta_{P_0}^\vee)} \sigma_M(H-\xi) q^{-\bg\la,H\bd}.
\end{equation} 
Puisque  $\xi$ est en position générale, on a $\xi\notin  \ZZ(\Delta_{P_0}^\vee$. Par \eqref{eq:sigmaM}, le membre de droite de \eqref{eq:egalitecarfct} est nul, ce qu'il fallait démontrer. Par prolongement analytique, il est aussi nul pour tout $\la \in a_T^*$.

 \end{preuve}

\end{paragr}

\begin{paragr}[Fonctions $\Gamma_P^Q(\cdot,X,\xi)$.] --- Pour tout sous-groupe parabolique $P\in \fc(T)$, on introduit, en suivant Arthur (cf. \cite{trace_inv} p.13), la fonction de  variables $H$ et $X\in  a_{P,\RR}$

\begin{equation}
  \label{eq:Gamma}
  \Gamma_P^{Q}(H,X,\xi)=\sum_{R\in \fc^Q(P)} (-1)^{\dim(a_R^Q)} \tau_P^R(H-X)\hat{\tau}_R^Q(H-\xi).
\end{equation}
Cette fonction ne dépend que des projetés de $H$ et $X$ sur $a_{P,\RR}^G$. 
Le lemme suivant précise un lemme d'Arthur (cf.  \cite{trace_inv} lemme 2.1). 

\begin{lemme}\label{lem:maj-support}
  Il existe une constante $C>0$ (qui dépend de $\xi)$ tel que pour tous $H$ et $X$ dans $a_{P,\RR}$ tels que 
$$\Gamma_P^{G}(H,X,\xi)\not=0$$
on a, pour tout $\al\in \Delta_P^G$,
$$|\bg \al, H\bd| \leq C(1+\sum_{\beta\in \Delta_P^G} |\bg \beta, X\bd|).$$
En particulier, la fonction $\Gamma_P^{G}(\cdot,X,\xi)$ est à support compact sur $a_{P,\RR}^G$.
\end{lemme}

\begin{preuve}
Il s'agit d'une variante de la preuve du  lemme 2.1 de \cite{trace_inv}.   Pour tous sous-groupes paraboliques $Q\subset R$ contenant $P$, soit $\tilde{\tau}_{Q}^{R}$ la fonction caractéristique de l'ensemble 
$$\{H \in a_{P,\RR}\mid \bg \al, H\bd >0 \  \forall \al \in \Delta_P^R-\Delta_P^Q \}.$$
 On a $\tilde{\tau}_{P}^{R}=\tau_P^R$ et pour tout sous-groupe parabolique $S$ contenant $R$
$$   \tilde{\tau}_{Q}^{R}\tilde{\tau}_{R}^{S}=\tilde{\tau}_{Q}^{S} .$$

On introduit alors pour $P\subset R$ des sous-groupes paraboliques la fonction 
\begin{eqnarray*}
 \tilde{\Gamma}_P^R(H,X)&=&\sum_{Q\in \fc^R(P)} (-1)^{\dim(a_P^Q)}\tau_P^Q(H-X) \tilde{\tau}_{Q}^{R}(H) \\
&=& \prod_{\al \in \Delta_P^R }( \tau_\al(H)-\tau_\al(H-X))
\end{eqnarray*}
où $\tau_\al$ est la fonction caractéristique des $H\in a_{P,\RR}$ tels que $\bg \al,H\bd >0$. Cette fonction ne prend que les valeurs $0,1$ ou $-1$. De plus, on a $\tilde{\Gamma}_P^R(H,X)\not=0$ seulement si on a pour tout $\al\in \Delta_P^R$
$$|\bg \al,H\bd| \leq |\bg \al,X\bd|.$$
Par ailleurs, on a le calcul suivant :
\begin{eqnarray*}
  \sum_{P\subset R\subset S} (-1)^{\dim(a_P^R)} \tilde{\Gamma}_P^R(H,X)\tilde{\tau}_{R}^{S}(H)&=& \sum_{P\subset  Q \subset R\subset S} (-1)^{\dim(a_P^R)} (-1)^{\dim(a_P^Q)} \tau_P^Q(H-X) \tilde{\tau}_{Q}^{R}(H) \tilde{\tau}_{R}^{S}(H)\\
&=&  \sum_{P\subset  Q \subset R\subset S} (-1)^{\dim(a_Q^R)}  \tau_P^Q(H-X) \tilde{\tau}_{Q}^{R}(H) \tilde{\tau}_{R}^{S}(H)\\
&=&  \sum_{P\subset  Q \subset R\subset S} (-1)^{\dim(a_Q^R)}  \tau_P^Q(H-X) \tilde{\tau}_{Q}^{S}(H)\\
&=&  \sum_{P\subset  Q \subset S}  \tau_P^Q(H-X) \tilde{\tau}_{Q}^{S}(H) \sum_{  Q \subset R\subset S} (-1)^{\dim(a_Q^R)}\\
&=& \tau_P^S(H-X).
\end{eqnarray*}
La dernière égalité vient du fait que la somme alternée dans la ligne précédente est nulle sauf si $Q=S$. En utilisant l'expression précédente pour  $\tau_P^S(H-X)$, il vient
\begin{eqnarray*}
  \Gamma_P^G(H,X,\xi)&=&\sum_{P\subset S\subset G }   (-1)^{\dim(a_S^G)} \tau_P^S(H-X)\hat{\tau}_S^G(H-\xi)\\
&=&\sum_{P\subset R\subset S\subset G }(-1)^{\dim(a_P^R)}(-1)^{\dim(a_S^G)} \tilde{\Gamma}_P^R(H,X)\tilde{\tau}_{R}^{S}(H) \hat{\tau}_S^G(H-\xi)\\
&=& \sum_{P\subset R \subset G }(-1)^{\dim(a_P^R)} \tilde{\Gamma}_P^R(H,X) \sum_{R\subset S\subset G }\tilde{\tau}_{R}^{S}(H) \hat{\tau}_S^G(H-\xi)
\end{eqnarray*}
On va prouver que chaque terme indexé par $R$ dans la somme ci-dessus vérifie la conclusion du lemme. 

Lorsque $R=P$, on a $\tilde{\Gamma}_P^P(H,X)=1$ et la contribution correspondante ne dépend pas de $X$ : elle vaut $\Gamma_P^G(H,0,\xi)$. La conclusion est alors claire car on sait (cf. \cite{trace_inv} lemme 2.1) que si $H\in a_{P,\RR}^G$ et $\Gamma_P^G(H,0,\xi)\not=0$, alors $H$ appartient à un compact qui ne dépend que de $\xi$. 

Pour un terme $R\subsetneq P$, on raisonne par récurrence. Soit $H$ tel que
\begin{equation}
  \label{eq:nonnul}
  \tilde{\Gamma}_P^R(H,X) \sum_{R\subset S\subset G }\tilde{\tau}_{R}^{S}(H) \hat{\tau}_S^G(H-\xi)\not=0.
\end{equation}
En particulier, on a $\tilde{\Gamma}_P^R(H,X)\not=0$ et donc pour tout  $\al\in \Delta_P^R$ on a $|\bg \al,H\bd| \leq |\bg \al,X\bd|$. Il reste à majorer  $|\bg \al,H\bd|$ pour $\al\in  \Delta_P-\Delta_P^R$. Écrivons $H=H_R+H^R$ suivant la décomposition $a_P^G=a_R^G+a^R_P$. On a donc pour tout $\al\in \Delta_P^R$ 
\begin{equation}
  \label{eq:maj-HR}
  |\bg \al,H^R\bd|=|\bg \al,H\bd| \leq |\bg \al,X\bd|.
\end{equation}
Soit 
$$l(H^R)=\sum_{\al\in \Delta_P-\Delta_P^R} \bg \al, H^R\bd \pi_\al\in a_R^G$$
où $(\pi_\al)_{\al\in \Delta_P-\Delta_P^R}$ est la base de $a_R^G$ duale de la base obtenue par projection sur $a_R^{G,*}$ de $\Delta_P-\Delta_P^R$. On a donc pour tout $\al\in \Delta_P-\Delta_P^R$ l'égalité  
\begin{equation}
  \label{eq:intermediaire}
   \bg \al, l(H^R)\bd=  \bg \al, H^R\bd
 \end{equation}
 et donc l'égalité
$$\tilde{\tau}_{R}^{S}(H)=\tau_R^S(H_R+l(H^R)).$$
En tenant compte de cette dernière égalité, on voit que la condition \eqref{eq:nonnul} devient 
 $$\tilde{\Gamma}_P^R(H^R,X) \cdot \Gamma_R^G(H_R,-l(H^R),\xi)\not=0.$$
Par l'hypothèse de récurrence sur $\Gamma_R^G$, il existe $C>0$ tel que pour tout $H$ qui vérifie la condition précédente satisfait pour tout  $\al\in  \Delta_P-\Delta_P^R$ l'inégalité 
$$|\bg \al, H\bd |=|\bg \al, H_R\bd |\leq C(1+\sum_{\beta \in  \Delta_P-\Delta_P^R} |\bg \be, l(H^R) \bd |).$$
Pour conclure, on a par \eqref{eq:intermediaire}  $|\bg \be, l(H^R) \bd |= |\bg \be, H^R \bd |$ qu'on majore à une constante près par $\sum_{\gamma\in \Delta_P^R} | \bg \gamma, H^R\bd  |$ donc par $\sum_{\gamma\in \Delta_P^R} | \bg \gamma, X\bd  |$, cf. \eqref{eq:maj-HR}.

\end{preuve}
\end{paragr}

\begin{paragr}[Les triplets duaux.] --- \label{S:duaux} La définition suivante nous sera commode dans la suite.

\begin{definition}\label{def:nu} Soit $M\in \lc$. Un triplet dual pour $M$ est un triplet $(\mathfrak{N},b,b')$ formé d'une partie finie de  $\frac{2\pi i }{\log(q)}\ZZ(\Delta_P)$, où $P$ est un élément quelconque de $\pc(M)$,  et d'entiers $b$ et $b' \geq 1$ qui vérifient les deux conditions suivantes
  \begin{enumerate}
  \item pour tout $P\in \pc(M)$ et tout sous-groupe parabolique  $R$ contenant $P$, on a 
  \begin{equation}
    \label{eq:inclusion}
    \ZZ(\hat{\Delta}_P^{\vee})\subset \ZZ(\hat{\Delta}_P^{R,\vee})\oplus \frac{1}{b}\ZZ(\Delta_R^{\vee}) \subset \frac{1}{b'}\ZZ(\hat{\Delta}_P^{\vee})\  ;
  \end{equation}
\item Pour tout $P\in \pc(M)$, l'ensemble $\mathfrak{N}$ est un système de représentants du quotient 
$$\frac{2\pi i }{\log(q)}\ZZ(\Delta_P)/ \frac{2\pi i b'}{\log(q)}\ZZ(\Delta_P).$$
  \end{enumerate}
\end{definition}

\begin{remarque}
  \label{rq:nu}
De tels triplets existent. Dans \eqref{eq:inclusion}, la première inclusion résulte du fait que les  projections respectives  sur $a_P^R$ de $\hat{\Delta}_P^{\vee}-\hat{\Delta}_R^{\vee}$ et $\hat{\Delta}_R^{\vee}$ sont les ensembles $\hat{\Delta}_P^{R,\vee}$ et $\{0\}$. 

Tout triplet dual  $(\mathfrak{N},b,b')$ vérifie  la propriété suivante : pour tout $P\in \pc(M)$ et tout $H\in \frac{1}{b'}\ZZ(\hat{\Delta}_P^{\vee})$, on a
\begin{equation}
  \label{eq:orth-triplet}
  \frac{1}{|\mathfrak{N}|}\sum_{\nu\in \mathfrak{N}} q^{-\bg \nu, H \bd}= \left\lbrace \begin{array}{l} 1 \text{ si }   H \in \ZZ(\hat{\Delta}_P^{\vee})  \\ 0 \text{ sinon } \end{array}\right..
\end{equation}
\end{remarque}

\end{paragr}

\begin{paragr}[Une variante de la construction du §\ref{S:cphi}.] --- Soit $M\in \lc$ et  $(\mathfrak{N},b,b')$ un triplet dual pour $M$. Soit $P\in \pc(M)$ et  $\varphi$ une fonction sur $a_{P}^*$. Pour tout  $\la\in a_{P}^*$, on pose 
\begin{equation}
  \label{eq:cPQxi}
  c_{P,\mathfrak{N}}^{G,\xi}(\varphi,\la)=\frac{1}{|\mathfrak{N}|}\sum_{\nu\in \mathfrak{N}} \sum_{\{R\mid P\subset R\subset G \}} (-1)^{\dim(a_R^G)} \hat{\vartheta}_P^R(\la+\nu) \varphi^R(\la+\nu) \vartheta_R^{G,b \xi}((\la+\nu)/b).
\end{equation}

La construction dépend réellement du triplet $(\mathfrak{N},b,b')$ mais pour ne pas alourdir davantage la notation  on omet la dépendance en $b$ et $b'$.

  \begin{lemme}\label{lem:FourierGamma}
Soit $X \in \ZZ(\hat{\Delta}_P^\vee)$ et $X'=X-\sum_{\al \in \Delta_P} \varpi_\al^\vee$.
\begin{enumerate}
\item La série de Fourier 
$$\sum_{H\in \ZZ(\hat{\Delta}_P^\vee)} \Gamma_P^G(H,X',\xi) q^{-\bg \la , H\bd }$$
est un polynôme de Laurent en les $q^{-\bg \la, \varpi_{\al}^\vee\bd }$ pour $\al\in \Delta_P$. En particulier, elle est  holomorphe sur $a_P^*$.
\item Cette série est égale à $c_{P,\mathfrak{N}}^{G,\xi}(\varphi,\la)$ pour la fonction $\varphi(\la)=q^{-\bg \la , X\bd }$ et  pour tout triplet $(\mathfrak{N},b,b')$ dual pour $M$. 
\end{enumerate}

  \end{lemme}

  \begin{preuve}    L'assertion 1 est une conséquence immédiate de la compacité du support de $H\in a_P^G\mapsto  \Gamma_P^G(H,X',\xi)$ (cf. lemme \ref{lem:maj-support}).  L'assertion 2 est un calcul analogue à celui du  lemme 2.2 p.15 de \cite{trace_inv}. Pour la commodité du lecteur, nous donnons quelques détails. Soit  $(\mathfrak{N},b,b')$ un triplet dual pour $M$. La fonction $c_{P,\mathfrak{N}}^{G,\xi}(\varphi,\la)$ est au moins  méromorphe. Pour montrer qu'elle coïncide avec la série de Fourier de $\Gamma_P^G(H,X',\xi)$ sur $a_P^*$, il suffit donc de le vérifier sur l'ouvert défini par  $\Re(\la^G)\in (a_P^{G,*})^+$. Sur ce domaine, on peut calculer la série de Fourier terme à terme selon la définition \eqref{eq:Gamma}. Il s'agit donc de calculer pour  $\Re(\la^G)\in (a_P^{G,*})^+$ et $R\in \fc(P)$ la série de Fourier de l'expression

    \begin{equation}
      \label{eq:QQQ}
      (-1)^{\dim(a_R^G)} \tau_P^R(H-X')\hat{\tau}_R^G(H-\xi).
    \end{equation}
 Soit  $\mathfrak{N}$ un ensemble dual de $\ZZ(\Delta_{P}^\vee)$. Par l'inclusion \eqref{eq:inclusion} et l'égalité \eqref{eq:orth-triplet} de la remarque \ref{rq:nu}, la  série de Fourier de  \eqref{eq:QQQ} s'écrit
$$\frac{1}{|\mathfrak{N}|}\sum_{\nu\in \mathfrak{N}} \sum_{H\in \ZZ(\hat{\Delta}_P^{R,\vee})\oplus \frac{1}{b}\ZZ(\Delta_R^\vee)}  (-1)^{\dim(a_R^G)} \tau_P^R(H-X')\hat{\tau}_R^G(H-\xi) q^{-\bg \la+\nu,H\bd}.$$

La série intérieure est alors un produit de séries géométriques qui convergent pour $\la\in a_P^+$ et dont la somme se calcule aisément : on trouve ainsi
$$ \frac{1}{|\mathfrak{N}|}\sum_{\nu\in \mathfrak{N}} (-1)^{\dim(a_R^G)}\prod_{\al\in \Delta_P^R} \frac{q^{- \bg \al,X\bd  \bg \la+\nu,\varpi^\vee_\al\bd  }}{1-q^{-\bg \la+\nu,\varpi^\vee_\al\bd  }} \prod_{\al \in \Delta_R^G} \frac{q^{-[\bg  \varpi_\al,b\xi\bd]\bg \frac{\la+\nu}{b},\al^\vee\bd  }}{1-q^{-\bg \frac{\la+\nu}{b},\al^\vee\bd  }} .$$
En observant qu'on a 
$$\prod_{\al\in \Delta_P^R} q^{-\bg \al,X\bd \bg \la+\nu,\varpi^\vee_\al\bd}=q^{ -\bg (\la+\nu)^R,X\bd },$$
on obtient
$$\frac{1}{|\mathfrak{N}|}\sum_{\nu\in \mathfrak{N}}(-1)^{\dim(a_R^G)} \hat{\vartheta}_P^R(\la+\nu)  q^{ -\bg (\la+\nu)^R,X\bd }  \vartheta_R^{G,b\xi}((\la+\nu)/b)$$
ce qui conclut.

  \end{preuve}

  Le théorème suivant est une variante du théorème \ref{thm:lissite}.

    \begin{theoreme}
      \label{thm:lissitexi}
Soit $(\varphi_\al)_{\al \in \Delta_P}$ une famille de fonctions méromorphes sur le disque ouvert 
$$\{z\in \CC \mid |z|<q^{r} \}$$
avec $r>0$ et holomorphes en dehors de l'origine $z=0$ et
\begin{equation}
  \label{eq:modelephi}
  \varphi(\la)=\prod_{\al\in \Delta_P} \varphi_\al(q^{-\bg \la,\varpi_\al^\vee\bd}).
\end{equation}

Soit $P$ un sous-groupe parabolique de $G$. La fonction $c_{P,\mathfrak{N}}^{G,\xi}(\varphi,\la)$ ne dépend du choix du triplet  dual pour $M_P$ qui intervient dans sa définition. De plus, cette fonction est   holomorphe pour $\la$ tel que $\Re(\bg \la,\varpi_\al^\vee\bd)>-r$ pour tout $\al\in \Delta_P$. 
\end{theoreme}

\begin{remarque}
  \label{rq:notation}
Pour une fonction $\varphi$ comme dans le théorème ci-dessus on définit
\begin{equation}
  \label{eq:cPGxisansnu}
  c_{P}^{G,\xi}(\varphi,\la)=c_{P,\mathfrak{N}}^{G,\xi}(\varphi,\la)
\end{equation}
pour un triplet  $(\mathfrak{N},b,b')$  dual pour $M_P$ quelconque, puisque le membre de droite ne dépend pas de ce choix.
\end{remarque}

  \begin{preuve} Soit $P$ un sous-groupe parabolique de $G$ et $(\mathfrak{N},b,b')$ un triplet  dual pour $M_P$. Pour tout  $\eps>0$, soit  $\Omega_\eps$ le voisinage ouvert de  $ia_{P,\RR}^*$ défini par la condition  $|\Re(\bg \la,\varpi_\al^\vee\bd)|<\eps$ pour tout $\al\in \Delta_P$. 
 La fonction $\varphi(\la)$ est holomorphe sur l'ouvert défini par  $\Re(\bg \la,\varpi_\al^\vee\bd)>-r$.  La fonction $c_{P,\mathfrak{N}}^{G,\xi}(\varphi,\la)$ est donc méromorphe sur  cet ouvert. De plus, ses singularités éventuelles sont simples et portées par des hyperplans qui coupent tous $ia_{P,\RR}^*$.  À l'aide du lemme \ref{lem:Taylor}, on voit que pour obtenir l'holomorphie de $c_{P,\mathfrak{N}}^{G,\xi}(\varphi,\la)$, il suffit  de la  prouver sur un ouvert $\Omega_\eps$ pour $\eps>0$.

Introduisons le développement en série de Laurent de $\varphi_\al$
$$\varphi_\al(z)=\sum_{n\in \ZZ} \hat{\varphi}_\al(n) z^n$$
où $\hat{\varphi}_\al(n)=0$ dès que $n$ est assez négatif. On a donc aussi
$$\varphi(\la)=\sum_{X\in \ZZ(\hat{\Delta}_P^\vee)} \hat{\varphi}(X) q^{-\bg \la,X\bd}$$
où l'on pose pour toute famille $(k_\al)_{\al\in \Delta_P}\in \ZZ^{\Delta_P}$
$$\hat{\varphi}(\sum_{\al\in \Delta} k_\al \varpi_\al^\vee) =\prod_{\al \in \Delta_P}\hat{\varphi}_\al(k_\al).$$
Par conséquent, on a
\begin{eqnarray*}
  c_{P,\mathfrak{N}}^{G,\xi}(\varphi,\la)&=&\sum_{X\in \ZZ(\hat{\Delta}_P^\vee)} \hat{\varphi}(X) \cdot c_{P,\mathfrak{N}}^{G,\xi}( q^{-\bg \la,X\bd},\la)
\end{eqnarray*}
D'après le lemme \ref{lem:FourierGamma}, la fonction $c_{P,\mathfrak{N}}^{G,\xi}( q^{-\bg \la,X\bd},\la)$ ne dépend pas du choix du triplet dual  $(\mathfrak{N},b,b')$. Il en est donc de même pour  $c_{P,\mathfrak{N}}^{G,\xi}(\varphi,\la)$. En outre,  $c_{P,\mathfrak{N}}^{G,\xi}( q^{-\bg \la,X\bd},\la)$ est holomorphe sur $a_P^*$. Pour prouver l'holomorphie de $c_{P,\mathfrak{N}}^{G,\xi}(\varphi,\la)$ sur $\Omega_\eps$, il suffit de prouver que cette série converge uniformément sur cet ouvert. Toujours par le  lemme \ref{lem:FourierGamma}, on a, avec les notations de ce lemme,  pour tout $X\in \ZZ(\hat{\Delta}_P^\vee)$
\begin{eqnarray*}
c_{P,\mathfrak{N}}^{G,\xi}( q^{-\bg \la,X\bd},\la)&=&\sum_{H\in \ZZ(\hat{\Delta}_P^\vee)} \Gamma_P^G(H,X',\xi) q^{-\bg \la , H\bd }.
\end{eqnarray*}
À l'aide du lemme \ref{lem:maj-support}, on voit qu'il existe des constantes $c_1>0$ et $c_2>0$ telles que pour tous $X$ et  $H$ tels que  $\Gamma_P^G(H,X',\xi)\not=0$ et tout $\la\in \Omega_\eps$ on ait 
\begin{eqnarray*}
   |q^{-\bg \la , H\bd }|&\leq&  q^{\sum_{\al\in \Delta_P}  |\Re(\bg \la,\varpi_\al^\vee \bd)| \cdot|\bg \al ,H  \bd| }\\
&\leq &c_1 \prod_{\al\in \Delta_P} q^{ \eps c_2 |\Delta_P|\cdot |\bg \al ,X  \bd|   }.
\end{eqnarray*}
Le lemme \ref{lem:maj-support} implique aussi qu'il existe un polynôme en $|\Delta_P|$ variables à coefficients réels positifs tel que  le nombre de $H\in \ZZ(\hat{\Delta}_P^\vee)$ qui vérifient $\Gamma_P(H,X',\xi)\not=0$ est borné   par 
$$A_X=p((|\bg \al ,X  \bd|)_{\al\in \Delta_P}).$$ 
Ainsi on a la convergence uniforme cherchée si la série 
$$\sum_{X\in \ZZ(\hat{\Delta}_P^\vee)} |\hat{\varphi}(X)| A_X q^{ \eps c_2 |\Delta_P|\cdot |\bg \al ,X  \bd|   }$$
converge. Il suffit de traiter le cas où d'un polynôme $p$ qui est un monôme. On est alors ramené à montrer que pour tout $k\in \NN$, la série suivante converge
$$\sum_{n\in \ZZ} |\hat{\varphi}_\al(n)| |n|^k  (q^{ \eps c_2  |\Delta_P|})^{|n|}.$$
C'est le cas pourvu que $\eps c_2|\Delta_P|<r$.
  \end{preuve}
\end{paragr}

\begin{paragr}[Les $(G,M)$-cofamilles périodiques.] --- Introduisons la définition suivante qui est une variante  de la définition \ref{def:cofam}.

\begin{definition}\label{def:cofamper}
    Soit $M$ un sous-groupe de Levi de $G$ et  $\Omega$  un voisinage ouvert de  $ia_{M,\RR}^*$ invariant par le $\frac{2\pi i}{\log(q)} \ZZ(\Delta_P)$ où $P$ est un élément quelconque de $\pc(M)$.

Une $(G,M)$-cofamille $(\varphi_P)_{P\in \pc(M)}$  sur  $\Omega$  est  dite \emph{périodique} si chaque fonction $\varphi_P$  est invariante par le réseau $\frac{2\pi i}{\log(q)} \ZZ(\Delta_P)$.
\end{definition}

Voici le pendant du lemme \ref{lem:fctlisse}.

\begin{lemme}\label{lem:fctlisseper}
  Soit $(\varphi_P)_{P\in \pc(M)}$ une $(G,M)$-cofamille périodique sur $\Omega$. La fonction 
$$\varphi_M(\la)=\sum_{P\in \pc(M)} \hat{\vartheta}_P^G(\la) \varphi_P(\la)$$
est périodique et holomorphe sur $\Omega$.
\end{lemme}

\begin{preuve}
  La périodicité est évidente. Quant à l'holomorphie au voisinage de $\mu\in \Omega$, elle résulte du lemme \ref{lem:fctlisse} appliqué à la $(G,M)$-cofamille au voisinage de $\la=0$ définie par
  \begin{equation}
    \label{eq:etunecoGM}
    \big[\prod_{\al\in \Delta_P} \frac{\bg \la,\varpi_\al^\vee \bd }{ 1-q^{-\bg \la+\mu,\varpi_\al^\vee\bd}}\big] \varphi_P(\la+\mu).
  \end{equation}
  Vérifions qu'il s'agit bien d'une $(G,M)$-cofamille. Soit $P_1$ et $P_2$ deux éléments de $\pc(M)$ co-adjacents et soit $\varpi^\vee\in \Delta_{P_1}^\vee \cap (-\Delta_{P_2}^\vee)$ donné par l'assertion 3 du lemme \ref{lem:lesequi}. Il s'agit de voir que les expressions \eqref{eq:etunecoGM} pour $P=P_1$ et $P=P_2$ sont égales sur l'hyperplan $\bg \la,\varpi^\vee\bd=0$. Si $\bg \mu,\varpi^\vee\bd\not=0$, les deux sont nulles sur cet hyperplan. Si $\bg \mu,\varpi^\vee\bd=0$, on a  $\varphi_{P_1}(\la+\mu)= \varphi_{P_2}(\la+\mu)$ sur l'hyperplan $\bg \la,\varpi^\vee\bd=0$ puisque $(\varphi_P)_{P\in \pc(M)}$ est une $(G,M)$-cofamille. Par le lemme \ref{lem:lesequi}, le produit des  facteurs associés à un copoids distinct de $\pm\varpi^\vee$ dans le  crochet sont égaux sur cet hyperplan. Il reste à voir l'égalité des facteurs associés à $\pm\varpi^\vee$ : ceux-ci sont tous deux égaux à $(\log(q))^{-1}$
\end{preuve}
\end{paragr}

\begin{paragr}[Le principal théorème pour les $(G,M)$-cofamilles périodiques.] --- Le théorème suivant est l'analogue du théorème \ref{thm:ppal}.

\begin{theoreme}\label{thm:ppalper}
    Soit $(\varphi_P)_{P\in \pc(M)}$ une $(G,M)$-cofamille périodique telle que chaque $\varphi_P$ satisfait les hypothèses du théorème \ref{thm:lissitexi} (en particulier $\varphi_P$ est de la forme \eqref{eq:modelephi}). Pour tout $\xi$ en position générale, on a l'égalité suivante entre  fonctions holomorphes au voisinage de $ia_M^*$ 
$$\varphi_M(\la)=\sum_{P\in \pc(M)} c_{P}^{G,\xi}(\varphi_P,\la).$$
  \end{theoreme}

  \begin{preuve}
Elle est similaire à la démonstration du théorème \ref{thm:ppal}. Pour mener la preuve, on fixe un triplet  $(\mathfrak{N},b,b')$ dual pour $M$. D'après \eqref{eq:cPGxisansnu} et la définition \eqref{eq:cPQxi}, le membre de droite  s'écrit
$$    \sum_{P\in \pc(M)} \frac{1}{|\mathfrak{N}|}\sum_{\nu\in \mathfrak{N}} \sum_{\{R\mid P\subset R\subset G \}} (-1)^{\dim(a_R^G)} \hat{\vartheta}_P^R(\la+\nu) \varphi^R_P(\la+\nu) \vartheta_R^{G,b \xi}((\la+\nu)/b).$$
Par interversion des sommes, on obtient 
$$ \frac{1}{|\mathfrak{N}|}\sum_{\nu\in \mathfrak{N}} \sum_{R\in \fc^G(M)} (-1)^{\dim(a_R^G)}\vartheta_R^{G,b \xi}((\la+\nu)/b)  \sum_{P\in \pc^R(M)}  \hat{\vartheta}_P^R(\la+\nu) \varphi^R_P(\la+\nu).$$
Par une variante du  lemme \ref{lem:proc-indpQ}, la dernière somme ne dépend que de $M_R$ : on la note  
$$\varphi_M^{M_R}((\la+\nu)^{M_R}).$$
Avec cette observation et par une nouvelle interversion de sommes, l'expression précédente devient
$$\frac{1}{|\mathfrak{N}|}\sum_{\nu\in \mathfrak{N}} \sum_{L\in \lc^G(M)} (-1)^{\dim(a_L^G)}\varphi_{M}^{L}((\la+\nu)^L)\sum_{R\in \pc(L)} \vartheta_R^{G,b \xi}((\la+\nu)/b)$$
Le lemme \ref{lem:xigeneral} indique que pour $\xi$ en position générale, la somme intérieure est nulle sauf si $L=G$ auquel cas elle vaut $1$. La ligne précédente se simplifie en 
 $$\frac{1}{|\mathfrak{N}|}\sum_{\nu\in \mathfrak{N}}\varphi_{M}(\la+\nu).$$
Comme $\varphi_M$ est invariante par le réseau  $\frac{2\pi i}{\log(q)} \ZZ(\Delta_P)$ qui contient $\mathfrak{N}$, on trouve  $\varphi_{M}(\la)$ comme voulu.
  \end{preuve}
\end{paragr}

\section{Un calcul de séries}\label{sec:calculI}

\begin{paragr}[Adèles et mesures.] ---\label{S:calculI-adeles} Soit $C$ une courbe projective, lisse, géométriquement connexe, de corps des constantes un corps fini $\Fq$ de cardinal $q$.  Soit $g_C$ son genre, $F$ son  corps de fonctions,  $\AAA$ l'anneau des adèles de $F$ et $\oc\subset \AAA$ le sous-anneau de $\AAA$ ouvert, compact et maximal pour ces propriétés. Soit $D$ un diviseur sur $C$ qu'on écrit 
$$D=\sum_{c\in |F|} n_c [c] $$
comme somme formelle, à coefficients $n_c$ entiers et presque tous nuls, sur l'ensemble $|F|$ des places de $F$. Son degré est défini par
$$\deg(D)=\sum_{c\in |F|} n_c \deg(c)$$
où le degré de $c$ est définit comme le degré sur $\Fq$ du corps résiduel du complété $F_c$ de $F$ en $c$. Pour tout $c \in |F|$, soit  $\varpi_c$ une uniformisante du complété $F_c$. Soit
$$\varpi_D=\prod_{c\in |F|}\varpi_c^{-n_c} \in \AAA^\times.$$
On a un morphisme degré surjectif
$$\deg :\AAA^\times \to \ZZ$$
normalisé par $\deg(\oc^\times)=0$ et $\deg(\varpi_c)=\deg(c)$ pour tout  $c\in |F|$. On a $\deg(F^\times)=0$. On utilise sur le groupe $\AAA^\times$ des idèles de $F$ la valeur absolue $|\cdot|=q^{-\deg(\cdot)}$.

Soit $n\geq 1$ un entier et $G=GL(n)$ sur le corps fini $\Fq$. Soit $\ggo$ son algèbre de Lie ; plus généralement les algèbres de Lie des groupes algébriques qu'on considère sont notées par la lettre gothique correspondante. Dans cette section, on utilise librement les notations de la section \ref{sec:combi} pour le groupe $G$ et son  sous-tore maximal formé des matrices diagonales.  

Soit $\mathbf{1}_D$ la fonction sur $\ggo(\AAA)$ caractéristique de $\varpi_D\ggo(\oc)$. Le $\oc$-module $\varpi_D\ggo(\oc)$ ne dépend pas des choix des uniformisantes $\varpi_c$.

La composition du déterminant avec le morphisme degré fournit un morphisme degré surjectif
$$\deg : G(\AAA)\to \ZZ.$$
Pour tout $e\in \ZZ$, on note $G(\AAA)^e$ la fibre de ce morphisme en $e$.

 On fixe sur $G(\AAA)$ la mesure de Haar à gauche normalisée par $\vol(G(\oc))=1$. La mesure de Haar sur $\ggo(\AAA)$ est normalisée par $\vol(\ggo(F)\back \ggo(\AAA))=1$. Les mêmes normalisations sont utilisées pour les sous-groupes de $G$. Ces conventions de mesure donnent lieu au lemme suivant dont la démonstration est laissée au lecteur.

\begin{lemme}
    Soit $N$ le radical  unipotent d'un sous-groupe parabolique de $G$ et $\ngo$ l'algèbre de Lie . On a 
$$\vol(N(F)\back N(\AAA))= q^{\dim(N)(g_C-1)}$$
et
$$\int_{\ngo(\AAA)}\mathbf{1}_D(U) \, dU=q^{\dim(\ngo)(1-g_C+\deg(D))}.$$
\end{lemme}

Pour $s\in \CC$, soit $\zeta(s)$ la fonction zêta de la courbe $C$. C'est une fraction rationnelle en $q^{-s}$ qui coïncide sur l'ouvert $\Re(s)>1$ avec la série convergente en $q^{-s}$
$$\int_{\AAA^\times \cap \oc} |x|^s \ dx,$$
la mesure utilisée étant la mesure de Haar sur $\AAA^\times$ normalisée suivant nos conventions : $\vol(\oc^\times)=1$.  
\end{paragr}

\begin{paragr}[Choix d'un élément  nilpotent \og régulier par blocs\fg{}.] ---\label{S:calculI-nota} Soit $(e_i)_{1\leq i\leq n}$ la base canonique de $\Fq^n$.   Soit $d\geq 1$ un entier qui divise $n$ et $r=n/d$ le quotient.  Pour $1\leq j\leq r$, soit $V_j$ le sous-espace de dimension $d$ engendré par $\{e_{(j-1)d+1},\ldots, e_{jd}  \}$. On a donc une décomposition en somme directe 
  $$\Fq^n=V_1\oplus\ldots \oplus V_r.$$

Soit $P_0$ le stabilisateur du drapeau $V_1\subset V_1\oplus V_2\subset \ldots\subset V_1\oplus\ldots\oplus V_r$. Il s'agit d'un sous-groupe parabolique de $G$ dont le facteur de Levi $M_0=M_{P_0}$ (cf. §\ref{S:notations2}) s'identifie à $GL(V_1)\times \ldots \times GL(V_r)$. 
Soit 
$$X=X^G\in \ngo_{P_0}(\Fq)$$
 l'endomorphisme nilpotent de $\Fq^n$ défini par $X^G(e_{jd+i})=e_{(j-1)d+i}$ pour $2\leq j\leq r$ et $1\leq i \leq d$ et $X^G$ est nul sur $V_1$. L'orbite de $X$ est dite \og régulière par blocs\fg{}. Lorsque $d=n$ on obtient l'orbite nulle et lorsque $d=1$ on obtient l'orbite régulière au sens usuel. En général, l'orbite obtenue est celle d'un élément nilpotent dont la décomposition de Jordan possède $d$ blocs de taille $r$. On a donc dans la base canonique de $\Fq^n$
\begin{equation}
  \label{eq:lamat}
  X^G=\left(
\begin{array}{ccccc}
  0 & I & 0 &0 & 0\\
 &0 & I & 0 & 0\\
 & &  0 & \ddots & 0\\
 & & & 0 & I\\
 & &  &  & 0\\
\end{array}\right).
\end{equation}
Soit $P$ un sous-groupe parabolique de $G$ contenant $P_0$. Comme $X^G\in \pgo(\Fq)$, il se décompose en $X^G=X^P+X_P$ avec $X^P\in \mgo_P(\Fq)$ et $X_P\in \ngo_P(\Fq)$. Pour chaque bloc $GL$  de $M_P$, la matrice de $X_P$ est de la forme \eqref{eq:lamat}. On a $X^{P_0}=0$ pour $P=P_0$.

Soit $M_{X^P}$ le centralisateur de $X^P$ dans $M_P$. Le centralisateur de $X^G$ dans $G$ est le sous-groupe de $GL(n)$ formé des matrices du type suivant
\begin{equation}
  \label{eq:lamat2}
  \left(
\begin{array}{ccccc}
  A_1 & A_2 & A_3 &  & A_r\\
 &A_1  & A_2 & \ddots & \\
 & &  A_1 & \ddots & A_3\\
 & & & A_1& A_2\\ 
 & &  &  & A_1\\
\end{array}\right).
\end{equation}
où $A_i$ est une matrice carrée de taille $d$ et de surcroît inversible si $i=1$. Il y a une description similaire pour $M_{X^P}$ dans chaque bloc $GL$  de $M$. On observera qu'on a $G_{X^G}\subset P_0$ et $M_{X^P}N_P\subset P_0$. Le groupe $G_X$ a pour radical unipotent $N_0\cap G_X$ et  pour facteur de Levi $M_0\cap G_X$ (ce dernier groupe est isomorphe à $GL(d)$).
\end{paragr}

\begin{paragr}[Application $H_{0}$ et réseaux de copoids.] --- \label{S:reseaux} On continue avec les notations du paragraphe précédent. Pour tout sous-groupe parabolique semi-standard $P$ de $G$, l'espace $a_{P,\RR}$, qui est par définition $\Hom_{\RR}(a_{P,\RR}^*,\RR)$ cf. §\ref{S:notations2}, s'identifie aussi à $\Hom_\ZZ(X^*(P),\RR)$, où $X^*(P)$ est le groupe des caractères algébriques de $P$. Dans cet espace, on distingue le réseau 
$$\ago_P=\Hom_\ZZ(X^*(P),\ZZ).$$
On a un morphisme  surjectif, trivial sur $N_P(\AAA)$
$$H_P : P(\AAA) \to \ago_P$$
donné pour tous $p\in P(\AAA)$ et  $\chi\in X^*(P)$ par $\bg \chi, H_P(p)\bd=-\deg(\chi(p))$. La décomposition d'Iwasawa $G(\AAA)=P(\AAA)G(\oc)$ permet de prolonger $H_{P}$ en une application
$$H_P : G(\AAA) \to \ago_P$$
telle que pour tous  $p\in P(\AAA)$ et $k\in G(\oc)$ on ait $H_P(pk)=H_P(p)$.
On utilisera aussi l'application 
$$H_P^G : G(\AAA)\to a_{P,\RR}^G$$
obtenue par composition de $H_P$ avec la projection de $a_{P,\RR}$ sur  $a_{P,\RR}^G$ selon la décomposition  $a_{P,\RR}^G=a_{P,\RR}^G\oplus a_{G,\RR}$.

Ces notations s'appliquent en particulier au sous-groupe parabolique $P_0$ du paragraphe précédent. Par abus, on notera parfois $H_0$ et $H_0^G$ les applications $H_{P_0}$ et $H_{P_0}^G$. On notera que si $P_0'\in \pc(M_0)$ est distinct de $P_0$ alors les applications $H_{P_0}$ et $H_{P_0'}$ ne sont pas égales. Elles ont cependant même image, à savoir 
$$\ago_{M_0}=\Hom(X^*(M_0),\ZZ)=\ago_{P_0}.$$

 En prenant comme base de $X^*(M_0)$ les caractères donnés par le déterminant des blocs $GL$ de $M_0$,  on identifie $X^*(M_0)$ à $\ZZ^r$. De même, on identifie $X^*(T)$ à $\ZZ^n$. Dualement au morphisme de restriction $X^*(M_0)\to X^*(T)$, on a un morphisme $\ago_T\to \ago_{M_0}$ surjectif qui,  par extension des scalaires à $\RR$ redonne la projection $a_{T,\RR}\to  a_{M_0,\RR}$ (cf. §\ref{S:notations2}). Comme au  §\ref{S:notations2}, on identifie  $a_{M_0,\RR}$ à un sous-espace de $a_{T,\RR}$ et cette identification est concrètement donnée par l'application suivante
$$(x_1, \ldots, x_r) \mapsto \frac{1}{d}(x_1,\ldots,x_1,x_2,\ldots,x_2,\ldots, x_r,\ldots,x_r)$$
(chaque entrée est répétée $d$ fois). La projection  $a_{T,\RR}\to  a_{M_0,\RR}$ est la projection orthogonale sur le sous-espace $a_{M_0,\RR}$ où $a_{T,\RR}$ est identifié à $\RR^n$ muni du produit scalaire canonique.

 Pour $1\leq i \leq r-1$,  soit $\al_i\in \Delta_0=\Delta_{P_0}$ la racine obtenue par restriction de la racine $(t_1,\ldots,t_n)\mapsto t_{id}t_{id+1}^{-1}$ du tore $T$. On a alors  
$$\hat{\Delta}_0^\vee=\{\varpi_1^\vee,\ldots, \varpi_{r-1}^\vee \}$$
où $\bg \al_j , \varpi_i^\vee\bd=\delta_{i,j}$.

L'hyperplan  $a_{T,\RR}^G$ de  $a_{T,\RR}$ est formé des points dont la somme des coordonnées est nulle. La projection orthogonale sur $a_{T,\RR}^G$ du réseau $\ago_{M_0}$ (vu comme sous-groupe de  $a_{T,\RR}$) est le groupe $\frac{1}d \ZZ(\hat{\Delta}_0^\vee)$.

Soit $e\in \ZZ$ et $M_0(\AAA)^e=M_0(\AAA)\cap G(\AAA)^e$. Dans l'identification de $\ago_{M_0}$ avec $\ZZ^r$, l'ensemble  $H_0(M_0(\AAA)^e)$ s'identifie aux éléments de $\ZZ^r$ dont la somme des coordonnées vaut $-e$. L'ensemble $H_0^G(M_0(\AAA)^e)$, c'est-à-dire  la projection de  $H_0(M_0(\AAA)^e)$ sur  $a_{T,\RR}^G$, s'identifie au sous-ensemble de $\frac{1}d \ZZ(\hat{\Delta}_0^\vee)$ formé des éléments
$$H=\frac1d \sum_{j=1}^{r-1} n_j \varpi_j^\vee$$
qui vérifient
$$\sum_{j=1}^{r-1} j n_j \equiv -e \mod r .$$

Soit $W_0$ le groupe de Weyl de $M_0$. Il s'agit du quotient du normalisateur de $M_0$ dans le groupe de Weyl de $(G,T)$ par le groupe de Weyl de $(M_0,T)$. Ce groupe s'identifie au groupe des permutations de l'ensemble $\{1,\ldots,r\}$. Ce groupe agit sur $M_0$ et donc sur l'espace $a_{M_0}^*$ via l'action sur $\Fq^n$ qui permute les  sous-espaces $V_i$. De même, ce groupe agit simplement transitivement sur $\pc(M_0)$. Soit $(\varpi_1,\ldots,\varpi_{r-1})$ la base de $a_{M_0}^{G,\RR}$ duale de $\Delta_0^\vee$. On observe que pour tout $w\in W_0$ et tout $i\in \{1,\ldots,r-1\}$, on a  $w\cdot \varpi_{i}-\varpi_{i}\in d\ZZ(\Delta_0)$ et en  particulier, pour l'élément $\gamma$ défini ci-dessous
\begin{equation}
  \label{eq:invariance-W0}
 w\cdot \gamma-\gamma \in \frac{2\pi i }{\log q }\ZZ(\Delta_0).
\end{equation}

\begin{lemme}\label{lem:orthog-e}
Soit
\begin{equation}
  \label{eq:defgamma}
  \gamma=\frac{2\pi i }{r \log q }  \sum_{j=1}^{r-1}j \al_j= \frac{2\pi i }{d \log q }    \varpi_{r-1}
\end{equation}
Pour tout $H\in \frac{1}{d}\ZZ(\hat{\Delta}_{P_0}^\vee)$, on a
\begin{equation}
  \label{eq:orthog-e}
  \frac1r \sum_{k=0}^{r-1}   \exp(-2\pi i k e/r)  q^{-kd  \bg \gamma , H \bd }= \left\lbrace \begin{array}{l} 1 \text{ si }   H \in H_0^G(M_0(\AAA)^e)
 \\ 0 \text{ sinon } \end{array}\right..
\end{equation}

\end{lemme}

\end{paragr}

\begin{paragr}[Fonctions zêta $Z_{d,D}$ et $\tilde{Z}_{d,D}$.] --- Pour tout $s\in \CC$, on pose
$$Z_{d,D}(s)=\int_{GL(d,\AAA)} \mathbf{1}_D(X) |\det(X)|^s\, dX=\sum_{e\in \ZZ}\int_{GL(d,\AAA)^e} \mathbf{1}_D(X)\, dX \, q^{-es}  $$
et
$$  \tilde{Z}_{d,D}(s)= q^{sd\deg(D)}(1-q^{-s})q^{d^2(g_C-1)} \zeta(s+1)\zeta(s+2)\ldots \zeta(s+d).$$
La première expression est une série de Laurent en $q^{-s}$ alors que la seconde est une fraction rationnelle en $q^{-s}$.

\begin{lemme}\label{lem:zetas}
  \begin{enumerate}
  \item    La série qui définit $Z_{d,D}(s+d)$ converge absolument pour $\Re(s)>0$. Sur ce domaine, c'est une fraction rationnelle en $q^{-s}$ qui vaut
$$ Z_{d,D}(s+d)=\frac{q^{d^2(\deg(D)+1-g_C)}}{1-q^{-s}}\tilde{Z}_{d,D}(s).$$\\
\item La fraction rationnelle  $\tilde{Z}_{d,D}(s)$ en $q^{-s}$  est holomorphe pour $\Re(s)>-1$.
\item On a 
$$\tilde{Z}_{d,D}(0)=\vol(GL(d,F)\back GL(d,\AAA)^0).$$
\end{enumerate}
\end{lemme}

\begin{preuve}
 Par un changement de variables suivi de la décomposition d'Iwasawa on obtient
 \begin{eqnarray*}
   Z_{d,D}(s+d)&=&q^{d(s+d)\deg(D)}\int_{GL(d,\AAA)} \mathbf{1}_{\mathfrak{gl}(d,\oc)}(X) |\det(X)|^{s+d}\, dX\\
&=& q^{d(s+d)\deg(D)}\zeta(s+1)\zeta(s+2)\ldots \zeta(s+d).
 \end{eqnarray*}
Les assertions 1 et 2 se déduisent alors des propriétés classiques de la fonction $\zeta$. L'assertion 3 n'est autre que la formule classique pour le volume de $GL(d,F)\back GL(d,\AAA)^0$.
\end{preuve}

\end{paragr}

\begin{paragr}[Fonctions $K^P_{D,X}$.] --- \label{S:KPDX}On poursuit avec les mêmes notations qu'avant.  Soit $P\in \fc(P_0)$. Pour tout  $g\in G(\AAA)$, soit
\begin{equation}
  \label{eq:KPDX}
   K^{P}_{D,X}(g)= \int_{\ngo_P(\AAA)} \mathbf{1}_{D}(g^{-1}( X^P +U)g)\, dU.
 \end{equation}
  
 \begin{lemme}\label{lem:KPX}
   La fonction  $K^{P}_{D,X}$ est invariante à gauche par $M_{X^P}(F)N_P(\AAA)$ et à droite par $G(\oc)$. De plus, pour tout $m\in M(\AAA)$, on a 
$$K^{P}_{D,X}(m)=q^{(1-g_C+\deg(D))\dim(N)} q^{\bg 2\rho_P,H_P(m)\bd} \mathbf{1}_{D}(m^{-1}X^P  m).$$
 \end{lemme}

 \begin{preuve}
   Pour alléger les notations, on omet l'indice $P$ dans $M_P$ et $N_P$.  Pour $g=nmk$ selon la  décomposition d'Iwasawa $G(\AAA)=N(\AAA)M(\AAA)G(\oc)$, on a 
\begin{eqnarray*}
  K^P_{D,X}(g)&=&\int_{\ngo_P(\AAA)} \mathbf{1}_{D}(m^{-1}( X^P +U) m) \, dU \\
&=& \mathbf{1}_{D}(m^{-1}X^P  m) \int_{\ngo_P(\AAA)} \mathbf{1}_{D}(m^{-1}U m) \, dU\\
&=& \mathbf{1}_{D}(m^{-1}X^P  m) q^{\bg 2\rho_P,H_P(m)\bd}\int_{\ngo_P(\AAA)} \mathbf{1}_{D}(U ) \, dU \\
&=& q^{(1-g_C+\deg(D))\dim(N)} q^{\bg 2\rho_P,H_P(m)\bd} \mathbf{1}_{D}(m^{-1}X^P  m).
\end{eqnarray*}
Le lemme en résulte.
 \end{preuve}
\end{paragr}

\begin{paragr}[Intégrales $\hat{I}^P_{D,X}(H)$.] --- Soit  $H\in H_0(M_0(\AAA))$ et $P\in \fc(P_0)$. Soit $G(\AAA)^{H}$ l'ensemble des $g\in G(\AAA)$ tels que $H_{P_0}(g)=H$ et
  \begin{equation}
    \label{eq:IcPDX}
    \hat{I}^P_{D,X}(H)= \int_{ N_P(F) M_{X^P}(F)\back G(\AAA)^H} K_{D,X}^P(g)\,dg.
  \end{equation}

\begin{lemme}\label{lem:calculIc}
  Pour tout  $H\in H_0(M_0(\AAA))$, l'intégrale $\hat{I}^P_{D,X}(H)$ est convergente et est égale à 
$$\vol(GL(d,F)\back GL(d,\AAA)^0)^{|\Delta_P|+1}  q^{\deg(D)\dim(N_0)} \times$$
$$ \prod_{\al\in \Delta_0^P} \big[  q^{d^2(g_C-1-\deg(D))}  \int_{ GL(d,\AAA)^{\bg \al,dH\bd}}\mathbf{1}_D(X_\al)|\det(X_\al)|^d dX_\al\big]. $$
De plus, l'intégrale $\hat{I}^P_{D,X}(H)$ ne dépend que de la projection de $H$ sur $a_{P_0}^G$.
\end{lemme}

\begin{preuve}
Pour alléger les notations, on omet l'indice $P$ dans $M_P$ et $N_P$. Soit  $M(\AAA)^H=M(\AAA)\cap G(\AAA)^H$. On utilise la décomposition d'Iwasawa $G(\AAA)^H=N(\AAA)M(\AAA)^H G(\oc)$ pour dévisser l'intégrale. Cela fait apparaître au niveau des mesures la puissance $q^{-\bg 2\rho_P,H_P(m)\bd}$. En tenant compte des volumes $\vol(G(\oc))=1$ et $\vol(N_P(F)\back N_P(\AAA))=q^{(g_C-1)\dim(N_P)}$ ainsi que du \ref{lem:KPX}, on obtient 
$$\hat{I}^P_{D,X}(H)= q^{\dim(N)\deg(D)} \int_{ M_{X^P}(F)\back M(\AAA)^H}  \mathbf{1}_{D}(m^{-1}X^P  m)\, dm .$$
  On a  $M_{X^P}\subset M\cap P_0$ et si l'on pose $M_1=M_0\cap M_{X^P}$ et $N_1=M_{X^P}\cap N_{P_0}$, on a une  décomposition de Levi $M_{X^P}=N_1M_1$. On peut donc encore  appliquer la décomposition d'Iwasawa pour descendre l'intégrale ci-dessus au sous-groupe parabolique $M\cap P_0$ de $M$. On pose $N_0^P=M\cap N_{P_0}$ et  $M_0(\AAA)^H=M_0(\AAA)\cap G(\AAA)^H$. On note $2\rho_0^P$ la somme des racines du tore diagonal dans $\ngo_0^P$. On obtient 
$$\hat{I}^P_{D,X}(H)=q^{\dim(N)\deg(D)}\times q^{-\bg 2\rho_0^P,H\bd} \times $$
\begin{equation}
  \label{eq:privedeqN}
  \int_{    M_1(F)\back M_0(\AAA)^H  }      \int_{N_1(F)\back N_0^P(\AAA)}  \mathbf{1}_{D}(m^{-1} n^{-1}X^Pn^{-1}  m)\, dn\, dm .
\end{equation}

Calculons l'intégrale intérieure. Pour cela, on note  $(\ngo_0^P)'$ l'algèbre dérivée de $\ngo_0^P$ et $2(\rho_0^P)'$ la somme des racines du tore diagonal dans $(\ngo_0^P)'$.  Par le  lemme \ref{lem:mesure} ci-dessous, on obtient

\begin{eqnarray*}
& \displaystyle \int_{N_1(F)\back N_0^P(\AAA)}  \mathbf{1}_{D}(m^{-1} n^{-1}X^Pn^{-1}  m)\, dn = q^{ (g_C-1)\dim(N_0^P)}\displaystyle \int_{(\ngo_0^P)'(\AAA)  }\mathbf{1}_{D}(m^{-1} (X^P+U)  m)\, dU\\
 &  = q^{ (g_C-1)\dim(N_0^P)}\displaystyle q^{\bg2 (\rho_0^P)',H_0(m)\bd}  \displaystyle\mathbf{1}_{D}(m^{-1}X^P m) \int_{(\ngo_0^P)'(\AAA)  }\mathbf{1}_{D}(U)\, dU \\
&=  q^{ (g_C-1)\dim(N_0^P)} q^{\dim((\ngo_0^P)')(1-g_C+\deg(D))}  q^{\bg2 (\rho_0^P)',H_0(m)\bd}  \displaystyle\mathbf{1}_{D}(m^{-1}X^P m).\\
\end{eqnarray*}

 Soit $M_1(\AAA)^0=M_1(\AAA)\cap G(\AAA)^H$ pour $H=0$. On observera que  le quotient $M_1(F)\back M_1(\AAA)^0$ est de volume fini. On a alors

$$\hat{I}^P_{D,X}(H)  = q^{ (g_C-1)(\dim(N_0^P)-\dim((\ngo_0^P)'))+\deg(D)(\dim(N)+\dim(\ngo_0^P)')}\times\vol(M_1(F)\back M_1(\AAA)^0)\times$$
\begin{equation}
  \label{eq:lintegrale}
   q^{-\bg 2\rho_0^P-2(\rho_0^P)',H\bd}  \int_{    M_1(\AAA)^0\back M_0(\AAA)^H  }     \mathbf{1}_{D}(m^{-1} X^P  m)\, dm.
  \end{equation}

On a $2\rho_0^P-2(\rho_0^P)'=\sum_{\al\in \Delta_0^P}d^2\al$. Traitons ensuite l'intégrale dans \eqref{eq:lintegrale}. Soit $r_1d,\ldots ,r_sd$ la taille des blocs $GL$ de $M$ avec $r_1+\ldots +r_s=r$. Le groupe $M_0$ s'identifie au produit $GL(d)^r$. Le centralisateur $M_1$ de $X^P$ dans $M_0$ s'identifie au groupe $GL(d)^s$ par le plongement qui envoie diagonalement le $i$-ème facteur $GL(d)$ dans le facteur $GL(d,\AAA)^{r_i}$. On écrit suivant cette identification $m=(x_1,\ldots,x_r)$ un élément $m$ de  $M_0(\AAA)$. Alors l'application
$$(x_1,\ldots,x_r)\mapsto (x^{-1}_1x_2,\ldots, x_{r_1-1}^{-1}x_{r_1},x_{r_1+1}^{-1}x_{r_1+2},\ldots,x_{r_s-1}^{-1}x_{r_s})$$
induit une bijection
$$M_1(\AAA)^0\back M_0(\AAA)^{H} \simeq \prod_{\al\in \Delta_0^P} GL(d,\AAA)^{\bg \al,dH\bd},$$
où $ GL(d,\AAA)^{\bg \al,dH\bd}$ désigne l'ensemble des éléments de  $GL(d,\AAA)$ de degré $\bg \al,dH\bd\in \ZZ$. Cette bijection préserve les mesures de sorte qu'on a 
$$ q^{-\bg 2\rho_0^P-2(\rho_0^P)',H\bd} \int_{    M_1(\AAA)^0\back M_0(\AAA)^{H}  }     \mathbf{1}_{D}(m^{-1} X^P  m)\, dm= \prod_{\al\in \Delta_0^P}\int_{ GL(d,\AAA)^{\bg \al,dH\bd}}\mathbf{1}_D(x_\al)|\det(x_\al)|^d dx_\al.$$
Dans ce produit, les intégrales sont convergentes : ce sont des coefficients des séries $Z_{d,D}$ (cf. lemme \ref{lem:zetas}). On en déduit que l'intégrale $\hat{I}^P_{D,X}(H)$ est convergente.
L'égalité du lemme se déduit de ce qui précède et  des égalités suivantes
$$(g_C-1)(\dim(N_0^P)-\dim((\ngo_0^P)'))= (g_C-1)d^2 |\Delta_0^P|,$$
$$\deg(D)(\dim(N)+\dim(\ngo_0^P)')= \deg(D)(\dim(N_0)-d^2 |\Delta_0^P|) $$
et
$$\vol(M_1(F)\back M_1(\AAA)^0)=\vol(GL(d,F)\back GL(d,\AAA)^0)^{|\Delta_P|+1}.$$
Enfin la dernière assertion est  une conséquence immédiate du calcul précédent. 
\end{preuve}

Dans la preuve précédente, on a utilisé le lemme suivant.

  \begin{lemme}\label{lem:mesure}
    Soit $N_X$ le centralisateur de $X$ dans $N_0$ et $\ngo_0'$ la sous-algèbre dérivée de $\ngo_0$. On a l'égalité suivante pour toute fonction $\varphi\in \Cc(\ngo_0'(\AAA))$
$$\int_{N_X(F)\back N_0(\AAA)} \varphi(n^{-1}Xn-X)\, dn=q^{ (g_C-1)\dim(N_0)}  \int_{\ngo_0'(\AAA)}\varphi(U)\,dU$$
  \end{lemme}

  \begin{preuve}
On introduit l'isomorphisme de variétés algébriques sur $\Fq$
$$\Phi:
\begin{array}{ccc}
  N_X\back N_0 &\to& \ngo'_0 \\ n &\mapsto & n^{-1}Xn-X
\end{array}$$
On montre que le changement de variables $n\mapsto \Phi^{-1}(\Phi(n)+U_0)$ pour $U_0\in \ngo_0'(\AAA)$ préserve la mesure sur  $N_0(\AAA)$.  L'intégrale dans le membre de gauche définit donc une mesure de Haar sur $\ngo_0'(\AAA)$ ainsi qu'évidemment le membre de droite. L'égalité à démontrer est par conséquent vraie pour une constante $c_0$. Il reste à calculer $c_0$. 

Pour cela, on évalue les deux membres en la fonction $\varphi_0=\mathbf{1}_{\ngo_0'(\oc)}$. On a 
\begin{eqnarray*}
 \int_{N_X(F)\back N_0(\AAA)} \varphi_0(n^{-1}Xn-X)\, dn&= &\vol(N_X(F)\back N_X(\AAA))\int_{N_X(\AAA)\back N_0(\AAA)} \varphi_0(n^{-1}Xn-X)\, dn \\
&=& \vol(N_X(F)\back N_X(\AAA)) \times \vol(N_X(\oc)\back N_0(\oc))\\
&=& q^{(g_C-1)\dim(N_X)}
\end{eqnarray*}
Alors que pour le second membre on a 
$$\int_{\ngo_0'(\AAA)}\varphi_0(U)\,dU=\vol(\ngo_0'(\oc))=q^{(1-g_C)\dim(\ngo_0') }.$$
On obtient  $c_0=q^{ (g_C-1)(\dim(N_X)+\dim(\ngo_0'))}$ ce qui donne le résultat puisqu'on a aussi $\dim(N_X)+\dim(\ngo_0')=\dim(N_0)$.
\end{preuve}

\end{paragr}

\begin{paragr}[Les séries $I^{\xi}_{P,D}(X,\la)$ et $I^\xi_D(X,\la)$.] --- \label{S:calculI}  Soit  $\xi\in a_{P_0,\RR}$ et $P\in \fc(P_0)$.  Pour tout $\la\in a_{0}^*$, on introduit la série
  \begin{equation}
    \label{eq:laserieIxiPD}
    I^{\xi}_{P,D}(X,\la)= \sum_{H\in H_{0}^G(M_0(\AAA))} \hat{\tau}_P(H-\xi) \hat{I}^P_{D,X}(H)   q^{-\bg d\la, H\bd}.
  \end{equation}
  
On définit alors la série
\begin{eqnarray} \label{eq:laserieIxiD}
  I^{\xi}_D(X,\la)&= & \sum_{P\in \fc(P_0) } (-1)^{\dim(a_P^G)} I^{\xi}_{P,D}(X,\la).\\
  &=&\nonumber  \sum_{H\in H_{0}^G(M_0(\AAA))} \big[\sum_{P\in \fc(P_0) } (-1)^{\dim(a_P^G)}\hat{\tau}_P(H-\xi) \hat{I}^P_{D,X}(H) \big]   q^{-\bg d\la, H\bd}.
\end{eqnarray}

\end{paragr}

\begin{paragr}[Un calcul des séries $I^{\xi}_{P,D}(X,\la)$ et $I^\xi_D(X,\la)$.] --- \label{S:lapreuveduthm} Soit $\dim(\oc_X)$ la dimension de l'orbite de $X$ sous $G$. On introduit la fonction suivante de la variable $\la\in a_{P_0}^*$
$$\psi(\la)=\prod_{\al\in \Delta_0}\tilde{Z}_{d,D}(\bg \la,\varpi_\al^\vee \bd).$$

D'après le lemme \ref{lem:zetas}, c'est une fraction rationnelle en les $q^{\bg \la,\varpi_\al^\vee \bd}$ pour  $\al\in \Delta_0$. De plus, cette fonction est holomorphe sur l'ouvert défini par $\Re(\bg \la,\varpi_\al^\vee \bd)>-1$ pour tout $\al\in \Delta_0$.

L'énoncé ci-dessous fait intervenir les notations de la section \ref{sec:combi}.

  \begin{proposition}
    \label{prop:lecalculpourP} Soit $(\mathfrak{N},b,b')$ un triplet dual pour $M_0$ au sens de la définition \ref{def:nu}. Pour tout $P\in \fc(P_0)$, la série $I^{\xi}_{P,D}(X,\la)$ est absolument convergente sur l'ouvert  défini par $\Re(\bg \la, \varpi_\al^\vee\bd)>-1$ pour $\al\in \Delta_0^P$ et $\Re(\bg \la, \al^\vee\bd)>0$ pour $\al\in \Delta_P$.  De plus, sur cet ouvert,  elle est égale au produit de 
    
\begin{equation}
      \label{eq:factor1}
      \vol(GL(d,F)\back GL(d,\AAA)^0) q^{\deg(D) \dim(\oc_X)/2}
    \end{equation}
et
\begin{equation}
  \label{eq:factor2}
  \frac1{|\mathfrak{N}|}\sum_{\nu \in \mathfrak{N}}   \hat{\vartheta}_{P_0}^P(\la+\nu) \psi^P(\la+\nu)\vartheta_{P}^{G,bd\xi}((\la+\nu)/b).
\end{equation}
  \end{proposition}

  \begin{preuve}Pour alléger les notations, on omet l'indice $P$ dans $M_P$ et $N_P$. En tenant compte du lemme \ref{lem:calculIc} et des égalités
$$\dim(\oc_X)=\dim(N_0)/2$$
et 
$$\psi^P(\la)=\psi(\la^P)= \vol(GL(d,F)\back GL(d,\AAA)^0)^{|\Delta_P|} \prod_{\al\in \Delta_0^P} \tilde{Z}_d(\bg \la,\varpi_\al^\vee \bd),$$
(cf.  lemme \ref{lem:zetas} assertion 3), il suffit de prouver la convergence de la série
\begin{equation}
  \label{eq:laserie}
  \sum_{H\in H_{P_0}^G(M_0(\AAA))} \hat{\tau}_P(H-\xi) \prod_{\al\in \Delta_0^P} \big[  q^{d^2(g_C-1-\deg(D))}  \int_{ GL(d,\AAA)^{\bg \al,dH\bd}}\mathbf{1}_D(X_\al)|\det(X_\al)|^d dX_\al\big]   q^{-\bg d\la, H\bd}
\end{equation}
sur l'ouvert et de montrer qu'elle y est égale à l'expression \eqref{eq:factor2} où l'on remplace  $\psi^P$ par le produit  $\prod_{\al\in \Delta_0^P} \tilde{Z}_d(\bg \la,\varpi_\al^\vee \bd)$.

On a vu à la fin du § \ref{S:reseaux} que le réseau  $H_{P_0}^G(M_0(\AAA))$ n'est autre que $\frac{1}{d}\ZZ(\hat{\Delta}_0^\vee)$. Par un changement de variables évident, la série \eqref{eq:laserie} se réécrit

\begin{equation}
  \label{eq:laserie2}
  \sum_{H\in \ZZ(\hat{\Delta}_0^\vee)} \hat{\tau}_P(\frac{H}{d}-\xi) \prod_{\al\in \Delta_0^P} \big[  q^{d^2(g_C-1-\deg(D))}  \int_{ GL(d,\AAA)^{\bg \al,H\bd}}\mathbf{1}_D(X_\al) |\det(X_\al)|^d dX_\al\big]   q^{-\bg \la, H\bd}
\end{equation}

En utilisant les propriétés d'un triplet dual $(\mathfrak{N},b,b')$ pour $M_0$, on écrit \eqref{eq:laserie2} sous la forme

\begin{equation}
  \label{eq:laserie3}
  \frac1{|\mathfrak{N}|} \sum_{\nu \in \mathfrak{N}}   \sum_{H\in \ZZ(\hat{\Delta}_0^{P,\vee})\oplus \frac1b \ZZ(\Delta_P^\vee)} \hat{\tau}_P(\frac{H}{d}-\xi) \prod_{\al\in \Delta_0^P} \big[  q^{d^2(g_C-1-\deg(D))}  \int_{ GL(d,\AAA)^{\bg \al,H\bd}}\mathbf{1}_D(X_\al) |\det(X_\al)|^d dX_\al\big]   q^{-\bg (\la+\nu), H\bd}
\end{equation}

Dans \eqref{eq:laserie3}, on va calculer la somme intérieure sur $H$. Pour cela, on décompose $\la+\nu$ ainsi
$$\la= \sum_{\al \in \Delta_{0}^P} \bg \la+\nu, \varpi_\al^\vee \bd \al+\sum_{\al \in \Delta_P  } \bg \la+\nu, \al^\vee\bd \varpi_{\al^\vee}.$$
En observant que la famille  $ (\bg\al,H \bd)_{\al \in \Delta_{0}^P }$ parcourt $\ZZ^{\Delta_{0}^P}$ et que la famille  $ (\bg   \varpi_{\al^\vee},H \bd)_{\al \in \Delta_P }$  parcourt $(\frac1b \ZZ)^{\Delta_P }$, on s'aperçoit que la somme sur $H$ s'exprime comme un produit dont certains facteurs donnent une fonction zêta alors que les autres sont des sommes géométriques que nous avons déjà rencontrées. On obtient alors
\begin{eqnarray}
  \nonumber     \prod_{\al\in \Delta_0^P} \big[  q^{d^2(g_C-1-\deg(D))} Z_{d,D}(d+ \bg \la+\nu,\varpi_\al^\vee \bd)\big] \times q^{-\bg \la, [b d \xi]\bd }\prod_{\al\in \Delta_P} \frac{1}{1-q^{-\bg (\la+\nu)/b,\al^\vee \bd}}   \\
\nonumber  =    \hat{\vartheta}_{P_0}^P(\la+\nu)    \big[\prod_{\al\in \Delta_0^P} \tilde{Z}_{d,D}( \bg \la+\nu,\varpi_\al^\vee \bd)\big]   \vartheta_{P}^{G,b d\xi}((\la+\nu)/b).\\
\end{eqnarray}
Cela donne l'égalité cherchée. Au passage, on obtient aussi la convergence sur l'ouvert défini par $\Re(\bg \la, \varpi_\al^\vee\bd)>-1$ pour $\al\in \Delta_0^P$ et $\Re(\bg \la, \al^\vee\bd)>0$ pour $\al\in \Delta_P$.
\end{preuve}

Dans l'énoncé suivant, on utilise la définition \eqref{eq:cPGxisansnu} de la remarque \ref{rq:notation}.

  \begin{proposition}\label{thm:calculintegrale}
La série qui définit $I^{\xi}_D(X,\la)$ est absolument convergente sur l'ouvert défini par $\Re(\bg \la,\varpi_\al^\vee \bd)>-1$ pour tout $\al\in \Delta_0$. Sur cet ouvert, elle est égale à 
$$ \vol(GL(d,F)\back GL(d,\AAA)^0) \cdot q^{\deg(D) \dim(\oc_X)/2} \cdot c_{P_0}^{G,d \xi}(\psi,\la).$$
\end{proposition}

\begin{preuve}
D'un point de vue formel, l'égalité cherchée est une conséquence de la proposition \ref{prop:lecalculpourP} et des définitions \eqref{eq:cPQxi} et  \eqref{eq:cPGxisansnu}. D'après la proposition \ref{prop:lecalculpourP}, la série  $I^{\xi}_D(X,\la)$ est au moins absolument convergente sur l'ouvert défini par $\Re(\bg \la,\al^\vee \bd)>0$ pour tout $\al\in \Delta_0$. D'après le théorème \ref{thm:lissitexi}, la fonction $c_{P_0}^{G,d \xi}(\psi,\la)$ est holomorphe  sur l'ouvert défini par $\Re(\bg \la,\varpi_\al^\vee \bd)>-1$ pour tout $\al\in \Delta_0$. La proposition en résulte.
\end{preuve}
\end{paragr}

\begin{paragr}[La série $I^{e,\xi}_D(X,\la)$.] --- Pour tout $e\in \ZZ$, on introduit la variante suivante de la série $I^{\xi}_D(X,\la)$ du §\ref{S:calculI} 
  $$I^{e,\xi}_D(X,\la)= \sum_{H\in H_{0}^G(M_0(\AAA)^e)} \big[\sum_{P\in \fc(P_0) } (-1)^{\dim(a_P^G)}\hat{\tau}_P(H-\xi) \hat{I}^P_{D,X}(H) \big]   q^{-\bg d\la, H\bd}.$$
On a le corollaire suivant de la proposition \ref{thm:calculintegrale}.

\begin{corollaire}\label{cor:calculintegrale-e}
  La série qui définit $I^{e,\xi}_D(X,\la)$ a les mêmes propriétés de convergence que la série $I^{\xi}_D(X,\la)$. De plus, sur son ouvert de convergence, elle est égale à 
$$ \vol(GL(d,F)\back GL(d,\AAA)^0) \cdot q^{\deg(D) \dim(\oc_X)/2} \cdot \frac1r \sum_{k=0}^{r-1}   \exp(-2\pi i k e/r) c_{P_0}^{G,d \xi}(\psi,\la+k\gamma)$$
pour tout $\gamma\in \frac{2\pi i }{r \log q } \ZZ(\Delta_0)$ qui vérifie la relation \eqref{eq:orthog-e} du lemme \ref{lem:orthog-e}.
\end{corollaire}

\begin{preuve}
  Pour tout $\gamma$ comme ci-dessus,  la relation \eqref{eq:orthog-e} implique qu'on a 
$$I^{e,\xi}_D(X,\la)= \frac1r \sum_{k=0}^{r-1}   \exp(-2\pi i k e/r)I^{\xi}_D(X,\la+k\gamma).$$
Le corollaire est alors une conséquence immédiate de la proposition \ref{thm:calculintegrale}.
\end{preuve}
\end{paragr}

\begin{paragr}[Calcul de la moyenne des séries.] --- La construction de la série $I^{\xi}_D(X,\la)$ dépend du  choix de $P_0\in \pc(M_0)$ et de $X$.  Tout autre parabolique $P_0'\in \pc(M_0)$ est de la forme $w\cdot P_0$ pour un unique élément de $W_0$, le groupe de Weyl de $M_0$ (cf. § \ref{S:reseaux}). Soit $I^{\xi}_D(\la, P_0')$ la série précédente construite pour les données $(w\cdot P_0,w\cdot X)$. On observe que pour tout $w\in W_0$, on a 
$$I^{\xi}_D(\la, P_0')=I^{w\cdot \xi}_D(X,w\cdot \la).$$

On définit alors la moyenne de ces séries sur $\pc(M_0)$ :
$$S^{\xi}(\la)=\frac{1}{|\pc(M_0)|}\sum_{P\in \pc(M_0)} I^{\xi}_D(\la, P).$$
Pour tout $e\in \ZZ$, on définit de même  $I^{e,\xi}(\la, P)$ pour $P\in \pc(M_0)$ et $S^{e,\xi}(\la)$. On a d'ailleurs
\begin{equation}
  \label{eq:Sexi}
  S^{e,\xi}(\la)= \frac1r \sum_{k=0}^{r-1}   \exp(-2\pi i k e/r)S^{\xi}(\la+k\gamma)
\end{equation}
pour tout $\gamma\in \frac{2\pi i }{r \log q } \ZZ(\Delta_0)$ qui vérifie la relation \eqref{eq:orthog-e} du lemme \ref{lem:orthog-e}.

Pour formuler nos résultats, il est commode d'introduire les notations suivantes :
\begin{itemize}
\item la fraction rationnelle en $q^{-s}$
$$\varphi(s)=  q^{d^2(g_C-1)} q^{sd\deg(D)} \zeta(s+1)\zeta(s+2)\ldots \zeta(s+d) \ ;$$
\item la famille  $(\varphi_{P})_{P\in \pc(M_0)}$
$$\varphi_P(\la)=  \prod_{\al\in \Delta_P} \varphi(\bg \la, \varpi_\al^\vee\bd).$$
\end{itemize}

  \begin{proposition}\label{cor:sommezeta} Pour tout  $\xi\in a_{P_0,\RR}$, la  série  $S^{\xi}(\la)$ est  convergente et holomorphe sur un voisinage de $ia^*_{M_0,\RR}$. Si, de plus, $\xi$ est en position générale, la somme de cette série  ne dépend pas de $\xi\in a_{P_0,\RR}$ et, sur ce voisinage de convergence, elle est   égale à 
$$ \vol(GL(d,F)\back GL(d,\AAA)^0) \cdot q^{\deg(D)\cdot  \dim(\oc_X)/2}\cdot \frac{1}{|\pc(M_0)|} \sum_{P\in \pc(M_0)} \varphi_P(\la).$$
  \end{proposition}

  \begin{preuve}
La proposition \ref{thm:calculintegrale} donne la convergence  de $S^{\xi}(\la)$ ainsi que l'expression suivante pour sa somme privée de facteurs évidents,
$$\sum_{P\in \pc(M_0)}c_{P}^{G,d \xi}(\psi_P,\la)$$
où l'on définit la $(G,M_0)$-cofamille $(\psi_P)_{P\in \pc(M_0)}$ par
$$\psi_P(\la)=\prod_{\al\in \Delta_P}\tilde{Z}_{d,D}(\bg \varpi_\al^\vee,\la\bd).$$

Supposons, de plus,  $\xi$  en position générale. Le théorème \ref{thm:ppalper} implique que l'expression précédente est égale à
$$\sum_{P\in \pc(M_0)} \hat{\vartheta}_P^G(\la)\psi_P(\la)$$
Le résultat vient alors de l'égalité suivante $\hat{\vartheta}_P^G(\la) \psi_P(\la)=\varphi_P(\la)$.
  \end{preuve}

  \begin{theoreme}\label{cor:sommezeta2}  Prenons  $\xi=0$ et supposons que le degré $e\in \ZZ$ soit premier à l'entier $r$ (rappelons que $r$ est le rang du centre de $M_0$). Alors la valeur en $\la=0$ de la  série $I^{e,0}_D(X,\la)$ est égale à la valeur en $\la=0$ de  l'expression holomorphe   dans un voisinage de $ia^*_{M_0,\RR}$ 
$$ \vol(GL(d,F)\back GL(d,\AAA)^0) \cdot q^{\deg(D)\cdot  \dim(\oc_X)/2}\cdot \frac{1}{r |\pc(M_0)|}\sum_{k=0}^{r-1} \sum_{P\in \pc(M_0)} \exp(-2\pi i k e/r) \varphi_P(\la+k\gamma)$$
pour  tout $\gamma\in \frac{2\pi i }{r \log q } \ZZ(\Delta_0)$ qui vérifie la relation \eqref{eq:orthog-e} du lemme \ref{lem:orthog-e}.
De plus, cette valeur ne dépend pas de l'entier $e$ premier à $r$. 
  \end{theoreme}

  \begin{preuve}
L'hypothèse que le degré de $e$ soit premier à $r$ entraîne que pour tout $H\in H_0(M_0(\AAA)^e)$ et tout $\varpi\in \hat{\Delta}_{P_0}$, on a  $\bg \varpi,H\bd\not=0$. Dans le cas contraire, le poids $\varpi$ détermine un sous-groupe parabolique maximal $P\in \fc(P)$ par la condition $\hat{\Delta}=\{\varpi\}$. Il existe $1\leq r'<r$ tel que $P$ soit le stabilisateur du sous-espace $V_1\oplus\ldots\oplus V_{r'}$ (avec les notations du §\ref{S:calculI-nota}). La dimension de ce sous-espace est $r'd$. Pour tout $H\in H_0(M_0(\AAA)^e)$, il existe un entier $e'\in \ZZ$ tel qu'on ait 
$$\bg \varpi,H\bd= e'-\frac{r'd}n e = e'-\frac{r'}{r}e $$
où $n=rd$ est le rang de $G$. Donc l'égalité $\bg \varpi,H\bd=0$ entraîne $re'=r'e$ et donc $r$ divise $r'$ ce qui contredit $1\leq r'<r$.

Il s'ensuit que pour $P\in \fc(P_0)$, tout $H\in H_0(M_0(\AAA)^e)$, tout $w\in W_0$ et tout élément $\xi\in a_{M_0,\RR}$ assez proche de $0$  on a 
$$\hat{\tau}_P(H)=\hat{\tau}_P(H-w\cdot \xi).$$
Prenons alors $\xi$ proche de $0$ pour que vaille l'égalité précédente. 
On en déduit que pour un tel $\xi$ on a $I^{e,0}_D(X,\la)=I^{e,\xi}_D(X,\la)$. On a même pour tous $P\in \pc(M_0)$ et $w\in W_0$ tels que $P=w\cdot P_0$ 
$$I^{e,0}_D(\la,P)= I^{e,0}_D(X,w\cdot\la)= I^{e,\sigma(\xi)}_D(X,w\cdot\la)=I^{e,\xi}_D(\la,P).$$
Il s'ensuit que les valeurs en $\la=0$ de $I^{e,0}_D(X,\la)$ et $S^{e,\xi}(\la)$ sont égales. Comme on peut supposer de plus $\xi$ en position générale, le théorème se réduit à une simple combinaison de l'égalité \eqref{eq:Sexi} et de la proposition \ref{cor:sommezeta}. 
\end{preuve}

\begin{remarque}\label{rq:mobius}
Bien que nous ne l'ayons pas vérifié, nous nous attendons à ce que la valeur de la série $I^{e,0}_D(X,\la)$ en $\la=0$ soit indépendante de l'entier $e$ premier à $r$. Du moins, elle ne dépend que de la réduction de $e$ modulo $r$. Lorsque $r$ est premier, on peut esquisser une preuve de cette indépendance. Avec les notations ci-dessus, il s'agit de vérifier que la limite pour $\la\to 0$ de la somme suivante est indépendante de $e$ premier à $r$ :
$$\sum_{k=0}^{r-1} \sum_{w\in W_0}  \exp(-2\pi i k e/r) \prod_{\al\in \Delta_0} \varphi( \bg  w\cdot(\la+k\gamma),\varpi_\al^\vee\bd ).$$
Il est commode de prendre $\gamma$ défini par \eqref{eq:defgamma}. Avec les notations du § \ref{S:reseaux},  et en tenant compte de la relation \eqref{eq:invariance-W0}, on voit que l'expression précédente se réécrit
$$\sum_{k=0}^{r-1} \sum_{w\in W_0}  \exp(-2\pi i k e/r) \prod_{j=1}^{r-1} \varphi( \bg  w\cdot\la ,\varpi_j^\vee\bd    + \frac{2\pi i kj }{r\log q }).$$
La contribution de $k=0$ s'écrit
$$\sum_{w\in W_0}  \prod_{j=1}^{r-1} \varphi( \bg  w\cdot\la ,\varpi_j^\vee\bd ).$$
qui a une limite en $\la=0$, évidemment indépendante de $e$.  La contribution de $k\not=0$ a aussi une limite en $\la=0$ qui est à un coefficient $|W_0|$ près
\begin{eqnarray*}
  \sum_{k=1}^{r-1}   \exp(-2\pi i k e/r) \prod_{j=1}^{r-1} \varphi(  \frac{2\pi i kj }{r\log q })&=& \prod_{j=1}^{r-1} \varphi(  \frac{2\pi i j }{r\log q })  \sum_{k=1}^{r-1}   \exp(-2\pi i k e/r)\\
&=& - \prod_{j=1}^{r-1} \varphi(  \frac{2\pi i j }{r\log q }) 
\end{eqnarray*}
et qui est donc indépendante de $e$ premier au nombre premier $r$.
\end{remarque}

\end{paragr}

\section{Induction de Lusztig-Spaltenstein}\label{sec:induction}

\begin{paragr} Dans cette section, on rassemble quelques résultats utiles sur l'induction de  Lusztig-Spaltenstein (cf. \cite{lus-spal}). Le groupe $G$ est le groupe $GL(n)$ sur un corps quelconque. Les autres notations sont celles utilisées dans les sections \ref{sec:combi} et \ref{sec:calculI}.
\end{paragr}

\begin{paragr}[Induite de Lusztig-Spaltenstein.] --- Soit $P\subset G$ un sous-groupe parabolique et $P=MN_P$ sa décomposition de Levi standard. Soit $X\in \pgo$ et  $\oc^P$ la $P$-orbite de $X$. La variété $\oc^P\oplus\ngo_P$ est irréductible et est formée d'éléments nilpotents de $\ggo$. Comme il n'y a qu'un nombre fini de $G$-classes de conjugaison nilpotentes dans $\ggo$, il existe une unique orbite nilpotente $\oc$ telle que l'intersection 
$$\oc\cap (\oc^P\oplus \ngo_P)$$
soit un ouvert dense de $\oc^P\oplus \ngo_P$. Cette intersection est une unique classe de $P$-conjugaison. On pose alors
$$I_P^G(X)=\oc.$$
On dit que $\oc$ est l'induite de $X$ de $P$ à $G$. Elle ne dépend par construction que de la $P$-orbite de $X$. On voit aussi qu'elle ne dépend que de $P/M_P$-orbite de la projection de $X$ sur $\pgo/ \ngo_P$. En particulier, si $\oc^M$ est une $M$-orbite nilpotente dans $\mgo$, l'induite $I_P^G(X)$ ne dépend pas du choix du représentant  $X\in \oc^M$. Aussi on la note
$$I_P^G(\oc^M).$$
\end{paragr}

\begin{paragr}[Transitivité.] --- Soit $P_1\subset P_2\subset G$ des sous-groupes paraboliques et $X\in \pgo_1(F)$. On définit une $M_2$-orbite nilpotente dans $\mgo_2$ par
$$I_{P_1}^{P_2}(X)=I_{P_1\cap M_2}^{M_2}(X_2)$$
où $X=X_2+U$ avec $X_2\in \mgo_2\cap\pgo_1$ et $U\in \ngo_2$. 

\begin{lemme}\label{lem:transitivite}
  Pour tout $Y\in I_{P_1}^{P_2}(X)$, on a 
$$I_{P_2}^G(Y)=I_{P_1}^G(X).$$
\end{lemme}
  
\begin{preuve}
Écrivons $X=X_1+U+V$ avec $X_1\in \mgo_1$, $U\in \ngo_1\cap \mgo_2$ et $V\in \ngo_2$. Soit $\oc^{P_1}(X)$, $\oc^{M_1}(X_1)$ et $\oc^{P_1\cap M_2}(X_1+U)$ les orbites respectivement de $X$ sous $P_1$, $X_1$ sous $M_1$ et $X_1+U$ sous $P_1\cap M_2$. On a  $\oc^{P_1}(X)\oplus\ngo_1=\oc^{M_1}(X_1)\oplus \ngo_1$. Soit $\oc=I_{P_1}^G(X)$. On a donc $\oc\cap (\oc^{M_1}(X_1)\oplus \ngo_1)$ est un ouvert dense de $\oc^{M_1}(X_1)\oplus \ngo_1$. Soit $\oc'=I_{P_1}^{P_2}(X)$.   Comme on a $\oc^{P_1\cap M_2}(X_1+U)\oplus (\ngo_1\cap \mgo_2)=\oc^{M_1}(X_1)\oplus (\ngo_1\cap \mgo_2)$, on en déduit que $\oc' \cap [\oc^{M_1}(X_1)\oplus (\ngo_1\cap \mgo_2)]$ est un ouvert dense de $\oc^{M_1}(X_1)\oplus (\ngo_1\cap \mgo_2)$. Ainsi  
$$  [\oc'  \cap [\oc^{M_1}(X_1)\oplus (\ngo_1\cap \mgo_2)]] \oplus \ngo_2$$
est un ouvert dense de $\oc^{M_1}(X_1)\oplus \ngo_1$. Par conséquent, 
$$ \oc \cap  [[\oc'  \cap [\oc^{M_1}(X_1)\oplus (\ngo_1\cap \mgo_2)]] \oplus \ngo_2]$$
est aussi un ouvert dense de $\oc^{M_1}(X_1)\oplus \ngo_1$. Donc $\oc \cap (\oc'\oplus \ngo_2)$ est non vide et c'est un ouvert de $\oc'\oplus \ngo_2$ vu qu'on a  $\oc'\oplus \ngo_2\subset \oc^{M_1}(X_1)\oplus \ngo_1$. D'où le résultat.
\end{preuve}
\end{paragr}

\begin{paragr}[Croissance de l'induite selon les sous-groupes paraboliques.] ---  On note $\overline{\oc}$ l'adhérence d'une orbite nilpotente $\oc$.
  
\begin{lemme}\label{lem:adh}
Soit $P_1\subset P_2$ des sous-groupes paraboliques de $G$.  On a alors pour tout $X\in \pgo_1(F)$
$$I_{P_2}^G(X)\subset \overline{I_{P_1}^G(X)}.$$    
\end{lemme}
  
\begin{preuve}
Soit $\oc^{1}(X)$ l'orbite de $X$ sous $P_1$. Soit $Y\in I_{P_2}^G(X)$. Quitte à conjuguer $Y$, on peut supposer qu'on a $Y\in X+\ngo_2$. Mais on a 
$$X+\ngo_2 \subset X+\ngo_1 \subset \oc^1(X)\oplus \ngo_1 \subset  \overline{I_{P_1}^G(X)},$$
d'où le résultat.
\end{preuve}

\begin{lemme}
  Soit $X\in \pgo_1(F)$ et $P_1\subset P_2$ des sous-groupes paraboliques de $G$ tels que $I_{P_1}^G(X)=I_{P_2}^G(X)$.   Alors pour tout sous-groupe parabolique $Q$ tel que  $P_1\subset Q\subset  P_2$, on a 
$$I_Q^G(X)=I_{P_1}^G(X)=I_{P_2}^G(X).$$
\end{lemme}

\begin{preuve} D'après le lemme \ref{lem:adh}, on a 
  $$I_Q^G(X)\subset \overline{I_{P_1}^G(X)}.$$
Si on n'a pas l'égalité dans cette inclusion, $I_Q^G(X)$ se trouve dans un fermé propre de $\overline{I_{P_1}^G(X)}$ et on a 
$$\overline{I_Q^G(X)}\subsetneq \overline{I_{P_1}^G(X)}.$$
Mais alors par une nouvelle application du lemme \ref{lem:adh}, il vient 
$$I_{P_2}^G(X)\subset \overline{I_Q^G(X)}\subsetneq \overline{I_{P_1}^G(X)}$$
ce qui contredit l'égalité  $I_{P_1}^G(X)=I_{P_2}^G(X)$.
\end{preuve}
\end{paragr}

\section{Quelques résultats auxiliaires}\label{sec:aux}

\begin{paragr}[Transformée de Fourier et formule sommatoire de Poisson.] --- Les notations sont celles de la  section \ref{sec:calculI}. Soit $\psi_0 :\Fq\to \CC^{\times}$ un caractère additif non trivial. Soit $\om_C$ une forme différentielle non nulle sur la courbe $C$. Pour toute place $c$ de $C$, on associe à $\om_C$ et à $\psi_0$ un caractère additif $\psi$ du corps $F_c$ complété de $F$ en $c$ par 
$$\psi(x)=\psi_0(\trace_{k_c/\Fq}(\Res(x \cdot\om_C))).$$
Cela détermine un caractère additif, noté encore $\psi$, du groupe $\AAA$ des adèles de $F$ qui est trivial sur $F$.

Soit $\ngo$ et $\bar{\ngo}$ des espaces vectoriels sur $\Fq$ en dualité via un accouplement parfait $\bg \cdot,\cdot\bd$. Les groupes $\ngo(\AAA)$ et $\bar{\ngo}(\AAA)$ sont munis des mesures de Haar qui donnent le volume $1$ à  $\ngo(F)\back \ngo(\AAA)$ et $\bar{\ngo}(F)\back \bar{\ngo}(\AAA)$. Pour toute fonction  $f\in \Cc(\ngo(\AAA))$ (c'est-à-dire localement constante et à support compact), on définit une fonction $\hat{f}\in  \Cc(\bar{\ngo}(\AAA))$ par
$$\hat{f}(Y)= \int_{\ngo(\AAA)} f(X) \psi(\bg X,Y\bd)\, dX$$
La formule sommatoire de Poisson est alors l'égalité
$$\sum_{X\in \ngo(F)} f(X)=\sum_{Y\in \bar{\ngo}(F)} \hat{f}(Y).$$

Soit $\Omega_C$ le diviseur de $\omega_C$. On a $\deg(\Omega_C)=2g_C-2$.  Pour tout diviseur $D$ sur $C$, on dispose de la fonction $\mathbf{1}_D \in  \Cc(\ngo(\AAA))$ qui est la fonction caractéristique de $\varpi_D\ngo(\oc)$. On a la formule suivante
\begin{equation}
  \label{eq:calcul-TFexplicite}
  \hat{\mathbf{1}}_D=q^{\dim(\ngo)(1-g_C+\deg(D))}   \mathbf{1}_{\Omega_C-D}.
\end{equation}

\begin{lemme}\label{lem:maj}
Dans le cas $\ngo=\Fq$, il existe une constante $c\geq 1$ telle que pour $a\in \AAA^\times$ et tout diviseur $D$, on ait
$$\sum_{X\in F} \mathbf{1}_D(a^{-1} X)  \leq   c \sup(q^{\deg(D)-\deg(a)},1)$$   
\end{lemme}

\begin{preuve} On a $\deg(\varpi_D)=-\deg(D)$. En  remplaçant  $a$ par $\varpi_D a$, on voit qu'il suffit de traiter le cas $D=0$. Soit $\mathbf{1}$ la fonction caractéristique de $\oc$.  Si $\deg(a) >0$, on a $\mathbf{1}(a^{-1} X)=1$ si et seulement si pour $X=0$.  Supposons ensuite $\deg(a)<2-2g_C$. Par la formule sommatoire de Poisson, on a 
$$\sum_{X\in F} \mathbf{1}(a^{-1} X)=q^{1-g_C-\deg(a)}\sum_{X\in F} \mathbf{1}_{\Omega_C}(a X)$$
et seul $X=0$ intervient non trivialement dans le membre de droite. On a donc dans ce cas
$$\sum_{X\in F} \mathbf{1}(a^{-1} X)=q^{1-g_C-\deg(a)}.$$
Il reste à majorer notre expression pour $2-2g_C\leq \deg(a)\leq 0$  par une constante $c\geq 1$. L'existence d'un tel $c$ résulte de la compacité de l'ensemble 
$$\{a\in F\back \AAA  \mid 2-2g_C\leq \deg(a)\leq 0\}.$$
\end{preuve}

\end{paragr}

\begin{paragr}[Théorie de la réduction et fonctions d'Arthur.] --- \label{S:reduction}
On utilise les notations des sections \ref{sec:combi} et \ref{sec:calculI}. Le groupe $G$ est le groupe $GL(n)$ sur le corps fini $\Fq$, muni de $ B\subset G$ son sous-groupe de Borel standard  et $T_0\subset B$ son  sous-tore maximal standard. On réserve la lettre $T$ pour un paramètre de troncature ; c'est pour cela que le sous-tore maximal de $G$ se retrouve affublé d'un indice $0$. On omet souvent les indices $T_0$ et $B$ pour les constructions associées à $T_0$ et $B$. Par exemple, on note $\Delta=\Delta_B$ (c'est l'ensemble des racines simples de $T_0$ dans $B$).

  Soit $\eta\in \AAA^\times $ un  idèle de degré $1$. L'application 
\begin{equation}
  \label{eq:lappli}
  \begin{array}{ccc}
  X_*(T_0) &\to & T_0(\AAA)\\
\lambda &\mapsto& \lambda(\eta)^{-1}
\end{array}
\end{equation} 
induit une bijection du groupe  $X_*(T_0)$ des cocaractères de $T_0$ sur un sous-groupe discret noté $A$ de $T_0(\AAA)$. On observera que si $a=\lambda(\eta)^{-1}$, l'accouplement canonique entre la cocaractère $\la$ et un caractère $\al\in X^*(T_0)$ vaut  $\bg \al, \la\bd= \bg \al, H_B(a)\bd$.

Soit $T_1\in a_{\RR}$ tel  que pour tout $\al\in \Delta$, on a $\al(T_1)<-2g_C$.   Soit
$$A^G(T_1)$$
l'image dans $A$ de l'ensemble
$$\{\la\in X_*(T_0) \mid \forall \al\in \Delta \ \bg \al,\la-T_1\bd >0\}.$$
On a alors l'égalité (on renvoie le lecteur à \cite{Harder-reduc} pour une référence).

$$G(\AAA)= G(F) \cdot B(\AAA)^0\cdot A^{G}(T_1)\cdot G(\oc).$$
où
$$B(\AAA)^0=B(\AAA)\cap \Ker(H_B).$$
Comme $B(F)\back B(\AAA)^0$ est compact, il est loisible de fixer  $\Omega_B$ une partie compacte de $ B(\AAA)^0$ telle que 
$$G(\AAA)= G(F) \cdot\Omega_B\cdot A^G(T_1)\cdot G(\oc).$$

Soit $T\in \overline{a_B^+}$ et $A^G(T_1,T)$ l'image dans $A$ de l'ensemble
$$\{\la\in X_*(T_0) \mid \forall \al\in \Delta \ \bg \al,\la-T_1\bd >0 \text{  et   }  \forall \varpi\in \hat{\Delta} \ \bg \varpi,\la-T\bd \leq 0 \}.$$
Soit 
$$F^{G}(\cdot,T)$$
 la fonction caractéristique de l'ensemble $G(F) \cdot\Omega_B\cdot A^{G}(T_1,T)\cdot G(\oc)$. La fonction $F^G$ est invariante à gauche par $G(F)$. Pour tout $e\in \ZZ$, l'ensemble
$$\{g\in G(F)\back G(\AAA)^e \mid F(g,T)=1\}$$
est compact.  

La même construction vaut pour tout sous-groupe de Levi $M$ de $G$. La fonction $F^M(\cdot,T)$ sur $M(\AAA)$ est la fonction caractéristique de $M(F)\Omega_B^M A^M(T_1,T) M(\oc)$ pour $\Omega_B^M$ un certain compact de $B(\AAA)^0\cap M(\AAA)$. Il est loisible de supposer qu'on a $\Omega_B^M=\Omega_B\cap M(\AAA)$. Pour tout sous-groupe parabolique standard $P$, soit $F^P(\cdot,T)$ la fonction sur $N_P(\AAA)M_P(F)\back G(\AAA)/G(\oc)$ qui coïncide sur $M_P(F)\back M_P(\AAA)/M_P(\oc)$ avec $F^{M_P}(\cdot,T)$, l'existence et l'unicité de $F^P$ résulte de la décomposition d'Iwasawa. 

Pour tout $T\in \overline{a_B^+}$, on a (pour une preuve, cf.,  par exemple, la proposition 10.3.12 de \cite{LGII})
\begin{equation}
  \label{eq:HN}
  \sum_{P\in \fc(B)}  \sum_{\delta\in P(F)\back G(F)} F^P(\delta g,T) \tau_P(H_P(\delta g)-T)=1.
\end{equation}
On a noté  $H_P$ est la composée de $H_B$ avec la projection sur $a_B$. 

Pour tout $T\in \overline{a_B^+}$, on a 
$$F^G(g,T)=\sum_{P\in \fc(B)} (-1)^{\dim(a_P^G)} \sum_{\delta\in P(F)\back G(F)} \hat{\tau}_P(H_P(\delta g)-T).$$
 Cet énoncé est une version  du lemme 2.1 de  \cite{ar_unipvar} pour les corps de fonctions ; il se démontre comme \emph{loc. cit.}, le point clef étant jouée dans notre situation par l'égalité \eqref{eq:HN}.

Soit $P_1\subset P_2$ des sous-groupes paraboliques standard  de $G$. Pour tout $H\in a_{T_0,\RR}$ soit 

$$  \sigma_{P_1}^{P_2}(H)= \sum_{P_2\subset P\subset G}(-1)^{\dim(a_{P_2}^P)}\tau_{P_1}^{P}(H) \hat{\tau}_P(H).$$

Pour tout $g\in G(\AAA)$ et $T\in \overline{a_B^+}$ soit
$$\chi^{P_1,P_2}_T(g)=F^{P_1}(g,T)\sigma_{P_1}^{P_2}(H_{P_1}(g)-T).$$
Cette fonction est invariante à gauche par $N_{P_1}(\AAA)M_{P_1}(F)$ et à droite par $G(\oc)$. On va en donner une description alternative. Pour tout $T\in \overline{a_B^+}$, soit  $A^{P_1,P_2}(T_1,T)\subset A$ l'ensemble des $a=\la(\eta^{-1})\in A$ pour $\la\in X^*(T_0)$ qui vérifie les cinq conditions suivantes
\begin{enumerate}
\item $\forall \al\in \Delta^{M_1} \ \bg\al ,\la-T_1\bd >0$ ;
\item $ \forall \varpi\in \hat{\Delta}^{M_1} \ \bg \varpi,\la-T\bd \leq 0 $ ;
\item $\forall \al\in \Delta_{M_1}^{M_2} \ \bg\al ,\la-T\bd >0$ ;
\item $\forall \al\in \Delta_{M_1}-\Delta_{M_1}^{M_2} \ \bg\al ,\la-T\bd \leq 0$
\item  $ \forall \varpi\in \hat{\Delta}_{M_2} \ \bg \varpi,\la-T\bd > 0$.
\end{enumerate}

\begin{remarques} \label{rq:condi}
  \begin{itemize}
  \item En vertu des conditions $1$ et $2$, on a $A^{P_1,P_2}(T_1,T)\subset A^{M_1}(T_1,T)$ ;
  \item   Si $P_1=P_2=G$, les conditions 3 à 5 sont vides et on a $A^{G,G}(T_1,T)= A^{G}(T_1,T)$ ;
  \item Si $P_1=P_2\not=G$, la condition 3 est vide mais les conditions 4 et 5 sont incompatibles : on a alors $A^{P_1,P_2}(T_1,T)=\emptyset$.
  \item Si $T\in \overline{a_B^+}$, on a $A^{P_1,P_2}(T_1,T)\subset A^{M_2}(T_1)$ (cela résulte du lemme ci-dessous).
  \end{itemize}
\end{remarques}

 \begin{lemme}\label{lem:2-1}
  Soit $P$ un sous-groupe parabolique de $G$. Soit $T\in \overline{a_B^+}$ et $X\in a_{B,\RR}$ qui vérifient

  \begin{enumerate}
  \item $ \forall \varpi\in \hat{\Delta}^{P} \ \bg \varpi,X-T\bd \leq 0 $ ;
  \item $\forall \al\in \Delta_{P} \ \bg\al ,X-T\bd >0$ ;
  \end{enumerate}
  Alors pour tout $\al\in \Delta^G-\Delta^P$, on a $ \bg\al ,X-T\bd >0.$
  \end{lemme}

  \begin{preuve}
  Soit $p$ la projection de $a_B$ sur $a_P$. Soit $\hat{\Delta}^P=  \{ \varpi_\al \mid \al\in \Delta^P  \}$ la base  de $a_B^{P,*}$  duale de $(\al^\vee)_{ \al\in \Delta^P}$. On a  pour tout $X\in  a_{B,\RR}$
  $$X=\sum_{\al\in \Delta^P} \bg\varpi_\al ,X\bd \al^\vee +p(X).$$
  Pour tout $\be\in \Delta^G-\Delta^P$, la projection $\be_P$ de $\beta$ sur $a_P^*$ appartient à $\Delta_P$ et
  $$\bg \be,X\bd=\sum_{\al\in \Delta^P} \bg\varpi_\al ,X\bd  \bg \be,\al^\vee\bd + \bg \be_P,X\bd.$$
  Or $\bg \be,\al^\vee\bd\leq 0$. Donc si $X$ vérifie les conditions 1 et 2, on a 
  $$\bg \be,X\bd >\sum_{\al\in \Delta^P} \bg\varpi_\al ,T\bd  \bg \be,\al^\vee\bd + \bg \be_P,T\bd= \bg \be,T\bd$$
  d'où le résultat.
  \end{preuve}

  \begin{lemme}(Arthur)\label{lem:condi}
  La fonction $$m\in M_1(\AAA)\mapsto \chi^{P_1,P_2}_T(m)$$
est la fonction caractéristique de $M_1(F)\cdot\Omega_B^{M_1}\cdot A^{P_1,P_2}(T_1,T)\cdot M_1(\oc)$.
\end{lemme}

\begin{remarque}
  Si $P_1=P_2\not=G$ la fonction    $\chi^{P_1,P_1}_T$ est identiquement nulle (cf. remarque \ref{rq:condi}). 
\end{remarque}

\begin{preuve}
  Il suffit de combiner \cite{ar1} lemme 6.1 avec la description de la fonction $F^{P_1}$.

\end{preuve}

La  fonction  $\chi^{P_1,P_2}_T $ va  intervenir via  le lemme suivant.

\begin{lemme}(Arthur)\label{lem:Arthur-dvpt}
Pour tout sous-groupe parabolique standard $P$ de $G$, soit  $K^P$ une fonction sur $P(F)\back G(\AAA)$. On a l'égalité entre
$$  \sum_{B\subset P\subset G} (-1)^{\dim(a_P^G)} \sum_{\delta\in P(F)\back G(F)} \hat{\tau}_P(H_P(\delta g)-T) K^{P}(\delta g)$$
et
$$\sum_{B\subset P_1\subset P_2 \subset G} \sum_{\delta\in P_1(F)\back G(F)}\chi^{P_1,P_2}_T(\delta g) K^{P_1,P_2}(\delta g)$$
où 
\begin{itemize}
\item l'on pose 
  \begin{equation}
    \label{eq:KPP}
    K^{P_1,P_2}(g)=\sum_{P_1\subset P \subset P_2} (-1)^{\dim(a_P^G)}  K^P(g)
  \end{equation}
\item dans la somme ci-dessus les termes correspondant à $P_1=P_2$ sont nuls sauf si $P_1=P_2=G$ auquel cas on obtient
$$F^G(g,T) K^G(g).$$
\end{itemize}
\end{lemme}

\begin{preuve}
  Elle repose sur les arguments formels utilisés pp. 41-43 dans \cite{ar-intro}.
\end{preuve}

\end{paragr}

\begin{paragr}[Une majoration de sommes rationnelles.] --- \label{S:maj-utiles} Soit $\al\in \Delta$ et $Q$ le sous-groupe parabolique standard maximal de $G$  défini par la condition 
$$ \Delta -\Delta^Q=\{\al\}.$$
Soit $\Sigma_Q$ l'ensemble des racines de $T_0$ dans $\ngo_Q$. On a une décomposition en espace poids
$$\ngo_Q=\oplus_{\al\in \Sigma_Q} \ngo_\al.$$
Pour tout $\Phi\subset \Sigma_Q$, on pose 
$$\ngo_\Phi=\oplus_{\al\in \Phi} \ngo_\al.$$
Tout élément $U\in \ngo_\Phi$ s'écrit $U=\sum_{\al\in \Phi} U_\al$ selon cette décomposition. Pour tout polynôme $g\in \Fq[\ngo_\Phi]$ et tout $\al\in \Phi$, on note $\deg_\al(g)$ le degré de $g(U)$ en la variable $U_\al$. On note $\vc((g_i)_{i\in I})$ le fermé de $\ngo_\Phi$ défini par une famille de polynômes $(g_i)_{i\in I}$.

\begin{proposition}\label{prop:maj-var}
Soit $d$ un entier.  Pour tout diviseur $D$ sur $C$, il existe une constante $c >0$ telle que pour toute partie  non vide $\Phi\subset \Sigma_Q$,  tout $a\in A^G(T_1)$ et toute famille de polynômes  non tous nuls $(g_i)_{i\in I}$  de $\Fq[ \ngo_\Phi]$ qui vérifie 
$$\forall i\in I, \ \forall \be\in \Phi, \deg_\be(g_i)\leq d$$
on a  l'inégalité
$$\prod_{\be\in \Phi}   q^{-\bg \be,H_B(a)\bd }  \sum_{ U\in \ngo_\Phi(F)\cap\vc((g_i)_{i\in I})} \mathbf{1}_D(a^{-1}Ua) \leq c \cdot  q^{ -\bg \al,H_B(a)\bd  }.$$
\end{proposition}

\begin{preuve} On va raisonner par récurrence sur le cardinal de $\Phi$. Le cas d'un singleton $\Phi=\{\be\}$ amorce la récurrence. Dans ce cas, le nombre de zéros d'un polynôme $g_i$ non nul ne peut pas excéder $d$ et on a donc 

$$ q^{-\bg \be,H_B(a)\bd }\sum_{ U\in \ngo_\Phi(F)\cap\vc((g_i)_{i\in I})}\mathbf{1}_D(a^{-1}Ua) \leq d  q^{-\bg \be,H_B(a)\bd }.$$
Il existe $\Delta'\subset\Delta-\{\al\}$ tel que $\be=\al+\sum_{\gamma\in \Delta'} \gamma$. La majoration $-\bg \gamma, H_B(a)\bd < - \bg \gamma, T_1\bd$ pour $a\in A^G(T_1)$ et $\gamma\in \Delta'$ donne le résultat.

Supposons désormais  $\Phi=\Phi'\coprod\{\be\}$ et $|\Phi'|\geq 1$. Tout $U\in \ngo_\Phi$ s'écrit $U=U'+U_\be$ selon la décomposition $\ngo_\Phi=\ngo_{\Phi'}\oplus\ngo_\be$. Soit $I'=I\times \{0,1,\ldots,d\}$. On obtient une famille de polynômes $(g'_i)_{i\in I'}$ de $\Fq[\ngo_{\Phi'}]$ non tous nuls en écrivant pour tout $i\in I$
$$g_i(U)=\sum_{k=0}^d g_{i,k}'(U')U_\be^k.$$
Pour alléger les notations, on introduit $\vc= \vc((g_i)_{i\in I})$ et $ \vc'=\vc((g'_i)_{i\in I'})$.
Soit $\vc^0$ l'ouvert de $\vc$ formé des $U$ dont la projection $U'$ n'appartient pas au fermé  $\vc'$ de $\ngo_{\Phi'}$.  On a donc une réunion disjointe 

 $$\vc= (\vc'\oplus \ngo_\be) \cup \vc^0.$$
Pour conclure, il suffit de prouver que les  contributions suivantes satisfont toutes deux  la majoration cherchée 

\begin{equation}
  \label{eq:vc'}
\prod_{\gamma\in \Phi'}   q^{-\bg \gamma,H_B(a)\bd}  \sum_{ U'\in \ngo_{\Phi'}(F)\cap \vc'}  \mathbf{1}_D(a^{-1}U'a) \times  q^{-\bg \beta,H_B(a)\bd} \sum_{ U_\be\in \ngo_\be(F)}  \mathbf{1}_D(a^{-1}U_\be a) 
\end{equation}
et 
\begin{equation}
  \label{eq:vc0}
   \prod_{\gamma\in \Phi}  q^{-\bg \gamma,H_B(a)\bd} \sum_{ U    \in \ngo_\Phi(F)\cap \vc^0} \mathbf{1}_D(a^{-1}Ua). 
\end{equation}

Pour majorer (\ref{eq:vc'}), on utilise d'une part l'hypothèse de récurrence appliqué au facteur attaché à $\Phi'$ et d'autre part que l'expression  $q^{-\bg \beta,H_B(a)\bd} \sum_{ U_\be\in \ngo_\be(F)}  \mathbf{1}_D(a^{-1}U_\be a)$ est bornée pour $a\in A^G(T_1)$ (cela résulte immédiatement du lemme \ref{lem:maj}).

Pour majorer (\ref{eq:vc0}), on commence par réécrire cette expression sous la forme

$$\prod_{\gamma\in \Phi'}   q^{-\bg \gamma,H_B(a)\bd} \sum_{U'\in \ngo_{\Phi'}(F), U'\notin \vc'}  \mathbf{1}_D(a^{-1}U'a) \times q^{-\bg \beta,H_B(a)\bd}  \sum_{ U_\be\in \ngo_\be(F), U'+U_\beta \in \vc}  \mathbf{1}_D(a^{-1}U_\be a ) .$$
 En utilisant le cas du singleton $\{\be\}$ et de la famille de polynômes $(g_i(U'+\cdot))_{i\in I}$ (non tous nuls puisque $U'\notin \vc'$), on a la majoration 
$$q^{-\bg \beta,H_B(a)\bd}  \sum_{ U_\be\in \ngo_\be(F), U'+U_\beta \in \vc}  \mathbf{1}_D(a^{-1}U_\be a ) \leq d\cdot  q^{-\bg \al,H_B(a)\bd}.$$

Il reste donc à borner 
$$\prod_{\gamma\in \Phi'}   q^{-\bg \gamma,H_B(a)\bd} \sum_{U'\in \ngo_{\Phi'}(F), U'\notin \vc'}  \mathbf{1}_D(a^{-1}U'a).$$
On majore trivialement cette expression par
$$\prod_{\gamma\in \Phi'}   q^{-\bg \gamma,H_B(a)\bd} \sum_{U'\in \ngo_{\Phi'}(F)}  \mathbf{1}_D(a^{-1}U'a).$$
qui est encore le  produit sur $\gamma\in \Phi'$ de 
$$q^{-\bg \gamma,H_B(a)\bd}  \sum_{ U_\gamma\in \ngo_\gamma(F)}  \mathbf{1}_D(a^{-1}U_\gamma a).$$
Cette expression est bien bornée pour $a\in A^G(T_1)$ (cf. lemme \ref{lem:maj}).
\end{preuve}
\end{paragr}

\section{Asymptotique d'intégrales nilpotentes tronquées}\label{sec:asymp}

\begin{paragr}
  Les notations sont celles en vigueur à la section \ref{sec:aux}.
\end{paragr}

\begin{paragr}\label{S:asymp}
  Soit $\nc^G$ le cône nilpotent de $\ggo(F)$. Soit $(\nc^G)$ l'ensemble (fini) des orbites de $G(F)$ dans $\nc^G$.
Pour toute orbite nilpotente  $\oc\in (\nc^G)$ et tout sous-groupe parabolique standard $P$ de $G$, on introduit la fonction suivante de la variable $g\in M_P(F)N_P(\AAA)\back G(\AAA)/G(\oc)$
$$ K^{P}_{D,\oc}(g)=\sum_{X\in \mathcal{N}^{M_P}, I_P^G(X)=\oc} \int_{\ngo_P(\AAA)}  \mathbf{1}_D(g^{-1}(X+U)g)\,dU.$$
Il se peut que l'ensemble de sommation soit vide auquel cas on pose $ K^{P}_{D,\oc}(g)=0$. C'est le cas, par exemple, si $\oc$ est l'orbite nulle et $P$ est un sous-groupe parabolique propre.

Pour tout $g\in G(F)\back G(\AAA)$, on pose
\begin{equation}
  \label{eq:KDTo}
  K^G_{D,T,\oc}(g)= \sum_{B\subset P\subset G} (-1)^{\dim(a_P^G)} \sum_{\delta\in P(F)\back G(F)} \hat{\tau}_P(H_B(\delta g)-T) K^{P}_{D,\oc}(\delta g).
\end{equation}

On fixe une norme euclidienne $\|\cdot \|$ sur $a_{B,\RR}$ invariante par le groupe de Weyl de $(G,T_0)$. Voici le principal résultat de cette section.

\begin{theoreme}\label{thm:approx}
Soit $e\in \ZZ$. Soit $\oc \in (\nc)$ et $D$ un diviseur sur $C$. Il existe un point $T_D\in a_B^+$, des constantes $\eps>0$, $\eps'>0$ et $c>0$ telles que pour tout $T\in T_D+a_B^+$ qui vérifie pour tout $\al\in \Delta$
$$\bg \al,T \bd \geq \eps' \| T\|,$$
on a 
$$\int_{G(F)\back G(\AAA)^e}  |  F^G(g,T) \sum_{X\in \oc} \mathbf{1}_D(g^{-1}Xg) - K^G_{D,T,\oc}(g)|\, dg \leq c\cdot q^{-\eps \|T\|}.$$

\end{theoreme}

Ce théorème est une conséquence  immédiate du lemme \ref{lem:Arthur-dvpt} et de la proposition \ref{prop:reduc1} ci-dessous.

\begin{corollaire}
  \label{cor:approx}
Les hypothèses sont celles du théorème \ref{thm:approx}. Pour $T\in a_B$, l'intégrale 
  $$\int_{G(F)\back G(\AAA)^e}   K^G_{D,T,\oc}(g)\, dg$$
converge absolument.

De plus, la fonction 
$$T\in X_*(T_0) \mapsto \int_{G(F)\back G(\AAA)^e}   K^G_{D,T,\oc}(g)\, dg$$
est quasi-polynomiale au sens où il existe un ensemble  fini 
$$\mathfrak{f}\subset \frac{2\pi i}{\log(q)} X^*(T_0)\otimes_\ZZ \QQ$$ 
et pour tout $\nu\in \mathfrak{f}$ un polynôme $p_\nu$  sur $a_B$ tels que cette fonction coïncide sur $X_*(T_0)$ avec 
$$T\mapsto \sum_{\nu\in \mathfrak{f} } p_\nu(T) q^{ \bg \nu,T\bd }.$$
\end{corollaire}

\begin{preuve}
  Elle repose sur les méthodes d'Arthur. Elle est d'ailleurs très semblable aux preuves du corollaire 5.1.2 et du théorème 5.2.1 de \cite{scfh}.  D'après Arthur (cf. \cite{trace_inv} section 2), on a, pour tout $H$ et $T$ dans $ a_B$
  \begin{equation}
    \label{eq:def-Gamma}
    \hat{\tau}_P(H-T)=\sum_{P\subset Q\subset G} (-1)^{\dim(a_Q^G)} \hat{\tau}_P^Q(H) \Gamma'_Q(H,T)
  \end{equation}
  où $\Gamma_Q'(H,T)$ est la fonction $\Gamma_Q^G(H,0,T)$ définie en \eqref{eq:Gamma} ; pour $Q=G$, la fonction $\Gamma'_G$ est identiquement égale à $1$  . Il s'ensuit qu'on a
\begin{equation}
  \label{eq:dvp-kT}
K^G_{D,T,\oc}(g)=\sum_{B\subset Q\subset G}  \sum_{\delta\in Q(F)\back G(F)}\Gamma'_Q(H_Q(\delta g),T)   K^{Q}_{D,0,\oc}(\delta g)
\end{equation}
où l'on pose pour $g\in G(\AAA)$
\begin{equation}
  \label{eq:fct-Q}
   K^{Q}_{D,0,\oc}(g)= \sum_{B\subset P\subset Q  } (-1)^{\dim(a_P^Q)} \sum_{\eta\in P(F)\back Q(F)}\hat{\tau}_P^Q(H_B(\eta g))K^{P}_{D,\oc} (\eta g).
\end{equation}
Pour $Q=G$, on retrouve   $K^{G}_{D,T,\oc}$ pour $T=0$. Soit $P\subset Q$ des sous-groupes paraboliques standard. Soit $\Sc^Q\subset (\nc^{M_Q})$ l'ensemble (fini) des orbites $\oc^Q\in (\nc^{M_Q})$ telles que $I_Q^G(\oc^Q)=\oc$.  On a successivement les égalités pour tout $m\in M_Q(\AAA)$ et $n\in N_Q(\AAA)$ (on note $\ngo_P^Q=\mgo_Q\cap \ngo_P$) :
\begin{eqnarray*}
   K^{P}_{D,\oc}(n m)&=&\sum_{X\in \mathcal{N}^{M_P}, I_P^G(X)=\oc} \int_{\ngo_P(\AAA)}  \mathbf{1}_D(m^{-1}n^{-1}(X+U)nm)\,dU\\
&=& \sum_{X\in \mathcal{N}^{M_P}, I_P^G(X)=\oc} \int_{\ngo_P(\AAA)}  \mathbf{1}_D(m^{-1}(X+U)m)\,dU\\
&=&\sum_{X\in \mathcal{N}^{M_P}, I_P^G(X)=\oc} \int_{\ngo_P^Q(\AAA)}\int_{\ngo_Q(\AAA)}  \mathbf{1}_D(m^{-1}(X+U+U')m)\,dU\, dU'\\
&=& q^{\bg2\rho_Q,H_Q(m) \bd} q^{\dim(N_Q)(1-g_C+\deg(D))}\times \\
& & \sum_{\oc^Q\in \Sc^Q} \sum_{\{X\in \mathcal{N}^{M_P} \mid I_{P\cap M_Q}^{M_Q}(X)=\oc^Q\}}\int_{\ngo_P^Q(\AAA)}  \mathbf{1}_D(m^{-1}(X+U)m)\,dU\\
&=& q^{\bg2\rho_Q,H_Q(m) \bd} q^{\dim(N_Q)(1-g_C+\deg(D))} \sum_{\oc^Q\in \Sc^Q} K^{P\cap M_Q}_{D,\oc^Q}(m).
\end{eqnarray*}

On en déduit qu'on a
\begin{equation}
  \label{eq:descente-par}
    K^{Q}_{D,0,\oc}(nm)=q^{\bg2\rho_Q,H_Q(m) \bd} q^{\dim(N_Q)(1-g_C+\deg(D))} \sum_{\oc^Q\in \Sc^Q} K^{M_Q}_{D,0,\oc^Q}(m)
\end{equation}

Nous allons faire une série de manipulations formelles qui seront \emph{a posteriori} justifiées. Nous avons d'après \eqref{eq:dvp-kT} et par l'utilisation de la décomposition d'Iwasawa
\begin{eqnarray*}
  \int_{G(F)\back G(\AAA)^e} \big(K^G_{D,T,\oc}(g)-K^G_{D,0,\oc}(g)\big)\,dg = \sum_{B\subset Q\subsetneq G}\int_{Q(F)\back G(\AAA)^e}\Gamma'_Q(H_Q(g),T)   K^{Q}_{D,0,\oc}(g) \, dg\\
= \sum_{B\subset Q\subsetneq G} q^{\dim(N_Q)\deg(D)} \sum_{\oc^Q\in \Sc^Q} \int_{M_Q(F)\back M_Q(\AAA)\cap G(\AAA)^e} \Gamma'_Q(H_Q(m),T)   K^{M_Q}_{D,0,\oc^Q}(m) \, dm \\
=   \sum_{B\subset Q\subsetneq G} q^{\dim(N_Q)\deg(D)} \sum_{\oc^Q\in \Sc^Q} \sum_{H\in H_Q(M_Q(\AAA)\cap G(\AAA)^e)} \Gamma'_Q(H,T)  \int_{M_Q(F)\back M_Q(\AAA)^H}  K^{M_Q}_{D,0,\oc^Q}(m) \, dm 
\end{eqnarray*}
où, dans la dernière ligne, $M_Q(\AAA)^H$ est l'ensemble des $m\in M_Q(\AAA)$ tels que $H_Q(m)=H$. Pour $Q\subsetneq G$, le sous-groupe de Levi $M_Q$ est un produit de facteurs linéaires, chacun de rang strictement plus petit que celui de $G$. Les fonctions $K^{M_Q}_{D,0,\oc^Q}$ et leurs intégrales sont aussi des produits indexés par ces facteurs. En raisonnant par récurrence sur le rang, on peut donc supposer que les intégrales  $ \int_{M_Q(F)\back M_Q(\AAA)^H} K^{M_Q}_{D,0,\oc^Q}(m)\, dm$ sont absolument convergentes. Comme la fonction $H\mapsto  \Gamma'_Q(H,T)$ est à support fini sur $H_Q(M_Q(\AAA)\cap G(\AAA)^e)$ (cf. lemme \ref{lem:maj-support}), cela justifie les manipulations précédentes et montre que  l'intégrale $\int_{G(F)\back G(\AAA)^e} K^G_{D,T,\oc}(g)\, dg $ est absolument convergente pour un $T \in a_B$ si et seulement elle l'est pour tout $T\in a_B$. Puisque l'intégrale 
$$\int_{G(F)\back G(\AAA)^e}    F^G(g,T) \sum_{X\in \oc} \mathbf{1}_D(g^{-1}Xg)\, dg$$
est absolument absolument convergente (l'intégrande est à support compact), le théorème \ref{thm:approx} montre qu'il existe au moins un tel $T$.

La fonction  $m\in M_Q(\AAA)\mapsto K^{M_Q}_{D,0,\oc^Q}(m)$ est invariante sous l'action du centre $Z_{Q}(\AAA)$ de $M_Q(\AAA)$. Il s'ensuit que l'application 
$$H\in H_Q(M_Q(\AAA)\cap G(\AAA)^e)\mapsto  \int_{M_Q(F)\back M_Q(\AAA)^H}  K^{M_Q}_{D,0,\oc^Q}(m) \, dm$$
 est invariante sous le sous-groupe $H_Q(Z_Q(\AAA)\cap G(\AAA)^e)$. Finalement pour tout $H\in  H_Q(M_Q(\AAA)\cap G(\AAA)^e)$ la fonction
$$T\mapsto \sum_{H'\in H_Q(Z_Q(\AAA)\cap G(\AAA)^e) } \Gamma_Q'(H+H',T)$$
est quasi-polynomiale (cf. proposition 4.5.5 de \cite{scfh}). Cela conclut.

\end{preuve}

\begin{proposition} \label{prop:reduc1}Sous les hypothèses du théorème \ref{thm:approx}, on a pour  tous sous-groupes paraboliques standard $P_1\subsetneq  P_2$ de $G$
$$\int_{G(F)\back G(\AAA)^e}    \sum_{\delta\in P_1(F)\back G(F)}\chi^{P_1,P_2}_T(\delta g)| K^{P_1,P_2}_{D,\oc}(\delta g)|\, dg \leq c\cdot q^{-\eps \|T\|}.$$
où l'on pose
\begin{equation}
  \label{eq:KP1P2}
  K^{P_1,P_2}_{D,\oc}(g)=\sum_{P_1\subset P \subset P_2} (-1)^{\dim(a_P^G)}  K^P_{D,\oc}(g).
\end{equation}
\end{proposition}

\begin{preuve} Elle se trouve au §\ref{S:preuvereduc1} après de longs préparatifs.  
\end{preuve}

\end{paragr}

\begin{paragr}[Une première majoration dans la preuve de  la proposition \ref{prop:reduc1}.] --- \label{S:reduc1} Dans toute la suite, on fixe  une orbite nilpotente $\oc \in (\nc)$ et des sous-groupes paraboliques standard $P_1\subsetneq  P_2$ de $G$. Pour tout $X\in \pgo_1(F)$, on pose
\begin{equation}
  \label{eq:xiP}
   \xi_{\oc}^{P_1,P_2}(X)= \sum_{P_1\subset P \subset P_2, I_P^G(X)=\oc} (-1)^{\dim(a_P^G)}.
\end{equation} 
Cette expression est invariante par conjugaison par $M_1(F)$ et ne dépend que la projection de $X$ sur $\mgo_1(F)$.

Pour alléger les notations, on note $P_1=M_1N_1$ (au lieu de $M_{P_1}N_{P_1}$) la décomposition de Levi standard du sous-groupe parabolique standard $P_1$. Idem pour $P_2$.

 \begin{proposition}
  \label{lem:step5} Pour tout diviseur $D$, il existe un point $T_D\in a_B^+$ tel que pour tous $T\in T_D+a_B^+$,  $n \in  N_1(\AAA)$ et $m\in  M_1(\AAA)$ tels que 
$$ \chi^{P_1,P_2}_T(m)\not=0$$
on a 
$$K^{P_1,P_2}_{D,\oc}(nm)= q^{\dim(\ngo_2)(1-g_C+\deg(D)) +\bg2\rho_2,H_2(m) \bd }  \sum_{X\in \mathcal{N}^{M_2}\cap \pgo_1(F)} \mathbf{1}_D(m^{-1}Xm) \xi_{\oc}^{P_1,P_2}(X).$$
\end{proposition}

\begin{preuve}
  C'est une conséquence immédiate du lemme \ref{lem:step4} ci-dessous.
\end{preuve}

\begin{corollaire}\label{cor:step5} Pour tout diviseur $D$, il existe un point $T_D\in a_B^+$ tel que pour tout $T\in T_D+a_B^+$ on ait l'inégalité 
  
 $$\int_{G(F)\back G(\AAA)^e}    \sum_{\delta\in P_1(F)\back G(F)}\chi^{P_1,P_2}_T(\delta g)| K^{P_1,P_2}_{D,\oc}(\delta g)|\, dg \leq  q^{\dim(N_1)(g_C-1)+\dim(\ngo_2)(1-g_C+\deg(D))} \times$$
$$\int_{M_1(F)\back M_1(\AAA)\cap G(\AAA)^e} q^{-\bg2\rho_1^2,H_1(m) \bd }  \chi^{P_1,P_2}_T(m) \sum_{X\in \mathcal{N}^{M_2}\cap \pgo_1(F)} \mathbf{1}_D(m^{-1}Xm)  |\xi_{\oc}^{P_1,P_2}(X)|\, dm$$
\end{corollaire}

\begin{preuve}
  On a 

$$\int_{G(F)\back G(\AAA)^e}    \sum_{\delta\in P_1(F)\back G(F)}\chi^{P_1,P_2}_T(\delta g)| K^{P_1,P_2}_{D,\oc}(\delta g)|\, dg=\int_{P_1(F)\back G(\AAA)^e}   \chi^{P_1,P_2}_T(g)| K^{P_1,P_2}_{D,\oc}(g)|\, dg.$$ 
 Par la décomposition d'Iwasawa, on obtient
$$\int_{N_1(F)\back N_1(\AAA)} \int_{M_1(F)\back M_1(\AAA)\cap G(\AAA)^e} \int_{G(\oc)} q^{-\bg 2\rho_{P_1}, H_{P_1}(m) \bd}  \chi^{P_1,P_2}_T(nmk) |K^{P_1,P_2}_{D,\oc}(nmk)|\,dk \,dm\,dn. $$
La fonction   $\chi^{P_1,P_2}_T$ est invariante à gauche par $N_1(\AAA)$ et à droite par $G(\oc)$. L'expression $|K^{P_1,P_2}_{D,\oc}(\cdot)|$ est invariante à droite par $G(\oc)$.  En tenant compte du volume $\vol(G(\oc))=1$, on obtient 
$$ \int_{N_1(F)\back N_1(\AAA)} \int_{M_1(F)\back M_1(\AAA)\cap G(\AAA)^e} q^{-\bg 2\rho_{P_1}, H_{P_1}(m) \bd}  \chi^{P_1,P_2}_T(m)|K^{P_1,P_2}_{D,\oc}(nm)| \,dm\,dn. $$
Pour terminer on utilise la formule $\vol(N_1(F)\back N_1(\AAA)) = q^{\dim(N_1)(g_C-1)}$ et l'expression donnée par la proposition \ref{lem:step5} qu'on majore de manière évidente.
\end{preuve}

\begin{lemme}
  \label{lem:step2}
Soit $P$ un sous-groupe parabolique de $G$ tel que $P_1\subset P\subset P_2$.
La fonction $g\in G(\AAA)\mapsto K^{P}_{D,\oc}(g)$ est invariante à gauche par $N_P(\AAA)$. De plus, pour tout $u\in  M_P(\AAA)\cap N_1(\AAA)$ et tout $m\in M_1(\AAA)$, on a 
$$K^{P}_{D,\oc}(um)= q^{\dim(\ngo_P)(\deg(D)+1-g_C) +\bg2\rho_P,H_P(m) \bd }  \sum_{X\in \mathcal{N}^{M_P}, I_P^G(X)=\oc} \mathbf{1}_D((um)^{-1}Xum).$$
\end{lemme}

\begin{preuve}
Soit $n\in N_P(\AAA)$ et $g\in G(\AAA)$. On a 
$$K^{P}_{D,\oc}(ng)=\sum_{X\in \mathcal{N}^{M_P}, I_P^G(X)=\oc} \int_{\ngo_P(\AAA)}  \mathbf{1}_D((ng)^{-1}(X+U)ng)\,dU.$$
On effectue le changement de variable $U'=n^{-1}(X+U)n-X$ ce qui donne
$$K^{P}_{D,\oc}(ng)=\sum_{X\in \mathcal{N}^{M_P}, I_P^G(X)=\oc} \int_{\ngo_P(\AAA)}  \mathbf{1}_D(g^{-1}(X+U)g)\,dU$$
et qui prouve l'invariance cherchée. 

Soit $u\in M_P(\AAA)\cap N_1(\AAA)$ et $m\in M_1(\AAA)$. À l'aide du changement de variable $U'=(um)^{-1}Uum$, on obtient l'égalité
$$K^{P}_{D,\oc}(um)=q^{\bg2\rho_P,H_P(m) \bd } \sum_{X\in \mathcal{N}^{M_P}, I_P^G(X)=\oc} \int_{\ngo_P(\AAA)}  \mathbf{1}_D((um)^{-1}Xum+U)\,dU.$$
Comme on a $(um)^{-1}Xum\in \mgo_P(\AAA)$, l'expression $\mathbf{1}_D((um)^{-1}Xum+U)$ est le produit de  $\mathbf{1}_D((um)^{-1}Xum)$ par $ \mathbf{1}_D(U)$. 
En utilisant l'égalité
$$\int_{\ngo_P(\AAA)}  \mathbf{1}_D(U)\,dU= q^{\dim(\ngo_P)(\deg(D)+1-g_C)},$$
on obtient l'égalité cherchée.
\end{preuve}

\begin{lemme}\label{lem:step4} Soit $P$ un sous-groupe parabolique de $G$ tel que $P_1\subset P\subset P_2$. Pour tout diviseur $D$, il existe un point $T_D\in a_B^+$ tel que pour tous $T\in T_D+a_B^+$,  $n \in  N_1(\AAA)$ et $m\in  M_1(\AAA)$ tels que 
$$ \chi^{P_1,P_2}_T(m)\not=0$$
on ait  
$$K^{P}_{D,\oc}(nm)=q^{\dim(\ngo_2)(\deg(D)+1-g_C) +\bg2\rho_2,H_{P_2}(m) \bd }  \sum_{X\in \mathcal{N}^{M_2}\cap \pgo_1(F), I_P^G(X)=\oc} \mathbf{1}_D(m^{-1}Xm).$$
\end{lemme}

\begin{preuve}
Par le lemme \ref{lem:step2}, on voit qu'il s'agit d'obtenir l'égalité suivante 
$$q^{\dim(\ngo_P)(\deg(D)+1-g_C) +\bg2\rho_P,H_P(m) \bd }  \sum_{X\in \mathcal{N}^{M_P}, I_P^G(X)=\oc} \mathbf{1}_D((nm)^{-1}Xnm)=$$
$$q^{\dim(\ngo_2)(\deg(D)+1-g_C) +\bg2\rho_2,H_{P_2}(m) \bd }  \sum_{X\in \mathcal{N}^{M_2}\cap \pgo_1(F), I_P^G(X)=\oc} \mathbf{1}_D(m^{-1}Xm)$$
pour $n \in  M_P(\AAA)\cap N_1(\AAA)$ et  $m\in  M_1(\AAA)$ tels que $ \chi^{P_1,P_2}_T(m)\not=0$.

On va d'abord montrer que pour tout  $n \in  M_P(\AAA)\cap N_1(\AAA)$ et  $m\in  M_1(\AAA)$ tel que $ \chi^{P_1,P_2}_T(m)\not=0$, on a
\begin{equation}
  \label{eq:firstly}
  \sum_{X\in \mathcal{N}^{M_P}, I_P^G(X)=\oc} \mathbf{1}_D((nm)^{-1}Xnm)= \sum_{X\in \mathcal{N}^{M_P}, I_P^G(X)=\oc} \mathbf{1}_D(m^{-1}Xm).
\end{equation}
Le résultat est tautologique si $P=P_1$. On suppose donc $P_1\subsetneq P$.
Puisque  $ \chi^{P_1,P_2}_T(m)\not=0$, on a $m\in M_1(F)\cdot\Omega_B^{M_1}\cdot A^{P_1,P_2}(T_1,T)\cdot M_1(\oc)$. Quitte à remplacer $n$ par un conjugué sous $M_1(F)$, on peut supposer $m\in \Omega_B^{M_1}\cdot A^{P_1,P_2}(T_1,T)\cdot M_1(\oc)$. Les deux membres de \eqref{eq:firstly} sont invariants par translations à droite de $m$ par $M_1(\oc)$ et par translations à gauche de $n$ par $M_P(F)\cap N_1(F)$. On peut donc  supposer $m\in \Omega_B^{M_1}\cdot A^{P_1,P_2}(T_1,T)$ et $n$ dans un compact fixé $C$ de  $M_P(\AAA)\cap N_1(\AAA)$. Pour tout $a\in A^{P_1,P_2}(T_1,T)$ et toute racine $\al$ de $T_0$ dans $M_P\cap N_1$, il existe $\Phi^1 \subset \Delta^{P_1}$ et $\Phi^P_1\not=\emptyset \subset \Delta^{P}-\Delta^{P_1}$
$$\bg \al, H_B(a)\bd > \sum_{\beta\in \Phi^1} \bg \beta, T_1\bd +  \sum_{\beta\in \Phi^P_1} \bg \beta, T\bd   $$
(ici on utilise le lemme \ref{lem:2-1}). En prenant $T_D$ convenablement, il s'ensuit que si  $T\in T_D+a_B^+$, alors pour tout $a\in A^{P_1,P_2}(T_1,T)$ on a 
$$a^{-1} (\Omega_B^{M_1})^{-1}C \Omega_B^{M_1}a\subset M_P(\oc)\cap N_1(\oc)$$
et donc aussi $m^{-1} C m \subset M_P(\oc)\cap N_1(\oc)$. L'égalité \eqref{eq:firstly} s'ensuit.

Montrons ensuite que pour tout   $m\in  M_1(\AAA)$ tel que $ \chi^{P_1,P_2}_T(m)\not=0$, on a
\begin{equation}\label{eq:secondly}
    \sum_{X\in \mathcal{N}^{M_P}, I_P^G(X)=\oc} \mathbf{1}_D(m^{-1}Xm)=\sum_{X\in \mathcal{N}^{M_P}\cap \pgo_1(F), I_P^G(X)=\oc} \mathbf{1}_D(m^{-1}Xm).
\end{equation}
On peut supposer $m\in \Omega_B^{M_1}\cdot A^{P_1,P_2}(T_1,T)$. Si l'égalité \eqref{eq:secondly} est en défaut, il existe $X\in  \mgo_P(F)-\pgo_1(F)\cap \mgo(F)$ tel que $ \mathbf{1}_D(m^{-1}Xm)\not=0$. Il existe alors une racine $\al$ de $T_0$ dans $M_P\cap \overline{N}_1$ (où $ \overline{N}_1$ est le radical unipotent \og opposé \fg à $N_1$) telle que la projection sur $\ngo_\al$ de $m^{-1}Xm$ soit non nulle. Si, de plus, $\al$ est minimale (pour l'ordre défini par $B$) on a nécessairement 
$$-\deg(\bg \al, H_B(a)\bd) \leq \deg(D).$$  
Là encore, on peut trouver $T_D$ pour que ce ne soit pas possible lorsque $T\in T_D+a_B^+$. Cela prouve  \eqref{eq:secondly}.

Enfin, pour $m\in  M_1(\AAA)$, l'expression
$$\sum_{X\in \mathcal{N}^{M_2}\cap \pgo_1(F), I_P^G(X)=\oc} \mathbf{1}_D(m^{-1}Xm)$$
est le produit de 
$$\sum_{X\in \mathcal{N}^{M_P}\cap \pgo_1(F), I_P^G(X)=\oc} \mathbf{1}_D(m^{-1}Xm)$$
par
\begin{equation}
  \label{eq:third1}
\sum_{Y\in  (\ngo_P\cap \mgo_2)(F)} \mathbf{1}_D(m^{-1}Xm).
\end{equation}

Pour conclure, il suffit de montrer que pour $m\in  M_1(\AAA)$ tel que $ \chi^{P_1,P_2}_T(m)\not=0$, l'expression \eqref{eq:third1} est égale à 
$$q^{\dim(\ngo_P\cap \mgo_2) (1-g_C+\deg(D))+\bg 2\rho_P^{2},H_B(m)\bd }.$$
Cette égalité se démontre à l'aide de la formule de Poisson et par des arguments semblables à ceux développés plus haut.
\end{preuve}

\end{paragr}

\begin{paragr}[Une deuxième majoration dans la preuve de  la proposition \ref{prop:reduc1}.] --- \label{S:reduc2}

Pour tout sous-groupe parabolique $Q\in \fc^{P_1}(B)$, soit
$$A^{1,2}_{Q,e}(T)$$
l'ensemble des $a\in A^{P_1,P_2}(T_1,T)\cap A^{M_Q}(T_1,0)\cap G(\AAA)^e$ qui vérifient $\tau_Q^{P_1}(H_Q(a))=1$.

  \begin{proposition}
    \label{prop:2ndemaj}
Pour tout diviseur $D$, il existe
\begin{itemize}
\item un diviseur $D'$ sur $C$
\item un point $T_D\in a_B^+$
\item pour tout $Q\in \fc^{P_1}(B)$ une constante $c_Q>0$
\end{itemize}
tels que pour tout $T\in T_D+a_B^+$, l'intégrale 

    \begin{equation}
      \label{eq:integ-amaj}
      \int_{M_1(F)\back M_1(\AAA)\cap G(\AAA)^e} q^{-\bg2\rho_1^2,H_1(m) \bd }  \chi^{P_1,P_2}_T(m) \sum_{X\in \mathcal{N}^{M_2}\cap \pgo_1(F)} \mathbf{1}_D(m^{-1}Xm)  |\xi_{\oc}^{P_1,P_2}(X)|\, dm
    \end{equation}

se majore par la somme sur $Q\in \fc^{P_1}(B)$ et $a\in A^{1,2}_{Q,e}(T)$ de 
$$c_Q\cdot q^{-\bg2\rho_Q^2,H_Q(a) \bd } \sum_{X\in \mathcal{N}^{M_2}\cap \pgo_1(F)} \mathbf{1}_{D'}(a^{-1}Xa)  |\xi_{\oc}^{P_1,P_2}(X)|.$$
  \end{proposition}

  \begin{preuve}
On rappelle que la somme 
$$\sum_{X\in \mathcal{N}^{M_2}\cap \pgo_1(F)} \mathbf{1}_D(m^{-1}Xm)  |\xi_{\oc}^{P_1,P_2}(X)|$$
est, comme fonction de $m\in M_1(\AAA)$ invariante à gauche par $M_1(F)$. On insère dans l'intégrale \eqref{eq:integ-amaj}l'identité \eqref{eq:HN} pour le groupe $G=M_1$ et le paramètre $T=0$. Par conséquent, l'intégrale \eqref{eq:integ-amaj} est égale à la somme sur $Q\in \fc^{P_1}(B)$ de l'intégrale sur $m\in (Q\cap M_1)(F)\back M_1(\AAA)\cap G(\AAA)^e$ de

$$ q^{-\bg2\rho_1^2,H_1(m) \bd } F^Q(m,0)\tau_Q^{P_1}(H_Q(m)) \chi^{P_1,P_2}_T(m) \sum_{X\in \mathcal{N}^{M_2}\cap \pgo_1(F)} \mathbf{1}_D(m^{-1}Xm)  |\xi_{\oc}^{P_1,P_2}(X)|.$$
Fixons un tel $Q$. Par la décomposition d'Iwasawa pour $M_1$ et le choix de nos mesures, on obtient l'intégrale sur $m\in  M_Q(F)\back M_Q(\AAA)\cap G(\AAA)^e$ et sur $n\in (N_Q\cap M_1)(F)\back (N_Q\cap M_1)(\AAA)$ de 
$$ q^{-\bg2\rho_Q^2,H_Q(m) \bd } F^{M_Q}(m,0)\tau_Q^{P_1}(H_Q(m)) \chi^{P_1,P_2}_T(n m) \sum_{X\in \mathcal{N}^{M_2}\cap \pgo_1(F)} \mathbf{1}_D((nm)^{-1}Xnm)  |\xi_{\oc}^{P_1,P_2}(X)|.$$
On va majorer cette dernière expression. On peut supposer $F^{M_Q}(m,0)\not=0$. Il s'ensuit qu'on peut supposer $m\in \Omega^{M_Q}_B A^{M_Q}(T_1,0)$ avec les notations de § \ref{S:reduction}. Comme $(N_Q\cap M_1)(F)\back (N_Q\cap M_1)(\AAA)$ est compact, on va supposer que $n$ appartient à un compact fixé, noté  $\Omega_Q^1$, de  $(N_Q\cap M_1)(\AAA)$. Pour $a\in  A^{M_Q}(T_1,0)$ tel que  $\tau_Q^{P_1}(H_Q(a))=1$, le lemme \ref{lem:2-1} montre qu'on a $\bg \al, H_B(a)\bd >0$ pour $\al\in \Delta^{P_1}-\Delta^Q$. Par ailleurs, pour $\al\in \Delta^{Q}$, on a $\bg \al, H_B(a)\bd >\bg \al, T_1\bd$. Il s'ensuit que pour $a\in  A^{M_Q}(T_1,0)$ tel que $\tau_Q^{P_1}(H_Q(a))=1$, l'ensemble $a^{-1} \Omega^1_Q \Omega_B^{M_Q}a$  reste dans un compact indépendant de $a$. Il existe donc un diviseur $D'$ tel que pour $n\in \Omega_Q^1$, $m\in \Omega^{M_Q}_B$ et $a\in  A^{M_Q}(T_1,0)$ on a 
$$\mathbf{1}_D((nm)^{-1}Xnm)\leq \mathbf{1}_{D'}(a^{-1}Xa).$$
On a aussi $H_Q(ma)=H_Q(a)$ et  $\chi^{P_1,P_2}_T(n m a)= \sigma_{P_1}^{P_2}(H_{P_1}(a)-T) F^{P_1}(nma,T)$. La condition  $   \chi^{P_1,P_2}_T(n m a)\not=0$ implique qu'on a $a\in A^{P_1,P_2}(T_1,T)$. 

La proposition s'ensuit en prenant pour la constante $c_Q$ le produit des volumes de $(N_Q\cap M_1)(F)\back (N_Q\cap M_1)(\AAA)$ et de $M_Q(F)\back M_Q(\AAA)^0$ (l'exposant $0$ signifie qu'on se limite au sous-ensemble des $m\in  M_Q(\AAA)$ qui vérifient $H_Q(m)=0$).
  \end{preuve}
\end{paragr}

\begin{paragr}[Une troisième majoration.] --- \label{S:reduc3}Dans ce paragraphe, on fixe $P\subsetneq G$ un sous-groupe parabolique standard \emph{propre}. Soit $P=MN$ sa décomposition de Levi standard.

\begin{proposition} \label{prop:reduc4}
  Soit $\oc\in (\nc^G)$ et $D$ un diviseur sur $C$.  Pour tout $\al\in \Delta-\Delta^P$, il existe une constante $c>0$ telle que pour tout $a\in A^G(T_1)$, on a 

$$  q^{-\bg2\rho_P,H_P(a) \bd } \sum_{X\in \nc^G \cap \pgo(F)} \mathbf{1}_D(a^{-1}Xa )  |\xi_{\oc}^{P,G}(X)|  \leq c\cdot  q^{-\bg\al,H_B(a)\bd}\cdot \sum_{X\in \mathcal{N}^{M}} \mathbf{1}_D(a^{-1}Xa).$$
\end{proposition}

Avant de donner la preuve de cette proposition, mentionnons le corollaire suivant. Les notations sont celles des paragraphes précédents.

\begin{corollaire} \label{cor:reduc4}
Pour tout sous-groupe parabolique $Q\in \fc^{P_1}(B)$, toute orbite nilpotente $\oc\in (\nc^G)$ et tout diviseur $D$ sur $C$, il existe des constantes $c>0$ et $r>0$ telles que pour tout $a\in A^{M_2}(T_1)$ on ait la majoration suivante
\begin{equation}
  \label{eq:ilfautmaj}
 q^{-\bg2\rho_Q^2,H_Q(a) \bd } \sum_{X\in \mathcal{N}^{M_2}\cap \pgo_1(F)} \mathbf{1}_{D}(a^{-1}Xa)  |\xi_{\oc}^{P_1,P_2}(X)|\leq
\end{equation}
\begin{equation}
  \label{eq:lamaj}
  c\cdot  q^{-\bg2\rho_Q^1,H_Q(a) \bd } \prod_{\al\in \Delta^2-\Delta^1} q^{-r \bg\al,H_B(a)\bd} 
  \cdot \sum_{X\in \mathcal{N}^{M_1}} \mathbf{1}_D(a^{-1}Xa).
\end{equation}
\end{corollaire}

\begin{preuve}
Observons qu'il suffit de majorer pour toute racine $\al\in \Delta^2-\Delta^1$ l'expression \eqref{eq:ilfautmaj} par
\begin{equation}
  \label{eq:lemajorant}
  c\cdot  q^{-\bg2\rho_Q^1,H_Q(a) \bd }  q^{- \bg\al,H_B(a)\bd}   \cdot \sum_{X\in \mathcal{N}^{M_1}} \mathbf{1}_D(a^{-1}Xa).
\end{equation}
On obtient alors immédiatement la majoration \eqref{eq:lamaj} pour $r^{-1}=|\Delta^2-\Delta^1|$. On obtient le majorant \eqref{eq:lemajorant} comme conséquence d'une part du lemme \ref{lem:desc} ci-dessous  et d'autre part de la proposition \ref{prop:reduc4} appliquée au groupe $G=M_2$.
\end{preuve}

\begin{lemme} \label{lem:desc}  Soit $\Sc \subset (\nc^{M_2})$ l'ensemble (éventuellement vide) formé des $M_2(F)$-orbites des éléments $Y\in \mgo_2(F)$ tels que $I_{P_2}^G(Y)=\oc$. Pour tout $X\in \mgo_2(F)\cap \pgo_1(F)$ nilpotent, on a 
$$\xi_{\oc}^{P_1,P_2}(X)=(-1)^{\dim(a_{P_2}^G)}\sum_{\oc'\in \Sc} \xi_{\oc'}^{P_1\cap M_2,M_2}(X).$$
 \end{lemme}

  \begin{preuve}
L'application $P\mapsto P\cap M_2$ induit une bijection entre l'ensemble des sous-groupes paraboliques $P$ de $G$ tels que $P_1\subset P \subset P_2$ et l'ensemble des sous-groupes paraboliques de $M_2$ compris entre $P_1\cap M_2$ et $M_2$. Soit $P$ un tel sous-groupe parabolique de $G$ et $Y\in I_{P\cap M_2}^{P_2}(X)$. D'après le lemme \ref{lem:transitivite}, on a 
$$I_P^G(X)=I_{P_2}^G(Y)$$
Il s'ensuit que $I_P^G(X)=\oc$ si et seulement si $I_{P_2}^G(Y)=\oc$.
Donc on a $I_P^G(X)=\oc$ si et seulement si  $I_{P\cap M_2}^{P_2}(X)\in \Sc$.
L'égalité cherchée s'en déduit.

  \end{preuve}

\begin{preuve}(de la proposition \ref{prop:reduc4}) Soit  $\al\in \Delta-\Delta^P$. Nous allons démontrer la proposition \ref{prop:reduc4} par récurrence sur la dimension de l'orbite $\oc$. Le cas de l'orbite nulle $(0)$ amorce la récurrence car dans ce cas $\xi_{(0)}^{P,G}(X)=0$ sauf si $X=0$ de sorte que l'expression à majorer se réduit à $  q^{-\bg2\rho_P,H_P(a) \bd }$ pour $a\in A^G(T_1)$. En se rappelant que $2\rho_P$ est la somme des racines de $T_0$ dans $N_P$, on voit que la majoration est obtenue pour la constante $c=q^{-\bg 2\rho_P -\al,T_1\bd}$.

 Soit $\overline{\oc}$ l'adhérence de l'orbite $\oc$ et 
 \begin{eqnarray*}
   \xi_{\overline{\oc}}^{P,G}(X)&=& \sum_{\oc'\in \nc^G, \oc'\subset\overline{\oc}}\xi_{\oc'}^{P,G}(X)\\ &=& \sum_{R} (-1)^{\dim(a_R^G)}
 \end{eqnarray*}
où la somme est prise sur les sous-groupes paraboliques $R$ de $G$ qui vérifient
\begin{itemize}
\item $P\subset R \subset  G$
\item $I_P^G(X)\subset \overline{\oc}$.
\end{itemize}
En utilisant l'hypothèse de récurrence, on voit qu'il suffit  de majorer l'expression
\begin{equation}
  \label{eq:ob}
    q^{-\bg2\rho_P,H_P(a) \bd } \sum_{X\in \nc^G \cap \pgo(F)} \mathbf{1}_D(a^{-1}Xa)  |\xi_{\overline{\oc}}^{P,G}(X)|.
  \end{equation}

On utilise ensuite le lemme d'annulation suivant.

\begin{lemme}
  \label{lem:annul} Soit $X\in \nc^G\cap \pgo(F)$.  
Supposons qu'il existe un sous-groupe parabolique maximal de $G$ qui contienne  tous les éléments minimaux  pour l'inclusion de l'ensemble
\begin{equation}
  \label{eq:lens}
  \{R\in \fc^G(P) \mid I_P^G(X)\subset \overline{\oc}  \}.
\end{equation}
Alors on a 
$$\xi_{\overline{\oc}}^{P,G}(X)=0.$$
\end{lemme}

\begin{preuve}
On suppose l'ensemble (\ref{eq:lens}) non vide sinon l'annulation cherchée est évidente. D'après le lemme \ref{lem:adh}, si $R$ appartient à l'ensemble (\ref{eq:lens}), il en est de même de tout sous-groupe parabolique qui le contient. Donc si $\{R_1,\ldots,R_n\}$ désigne l'ensemble des éléments minimaux de \eqref{eq:lens}, celui-ci se décrit comme l'ensemble des sous-groupes paraboliques de $G$ qui contiennent un des $R_i$. Par conséquent, on a 
$$\xi_{\overline{\oc}}^{P,G}(X)=\sum_{\{R \mid \exists i R_i\subset R\}} (-1)^{\dim(a_R^G)}.$$

Soit $Q\subset G$ un sous-groupe parabolique maximal de $G$ qui contient tous les $R_i$. Prouvons par récurrence sur $n$ que pour toute famille $\{R_1,\ldots,R_n\}$ de sous-groupes paraboliques de $G$ inclus dans $Q$, on a  
$$\sum_{\{R \mid \exists i R_i\subset R\}} (-1)^{\dim(a_R^G)}=0.$$
Le résultat est bien connu si $n=1$ et résulte de la bijection entre les sous-groupes paraboliques $R$ contenant $R_1$ et l'ensemble des parties de $\Delta_{R_1}$. Si $n>1$, on peut écrire

$$\sum_{\{R \mid \exists i R_i\subset R\}} (-1)^{\dim(a_R^G)}=\sum_{R_n\subset R\  } (-1)^{\dim(a_R^G)} +  \sum_{\{R \mid \exists i<n \, R_i\subset R\}} (-1)^{\dim(a_R^G)}  - \sum_{\{R \mid \exists i<1 \, (R_i,R_n)\subset R\}} (-1)^{\dim(a_R^G)} ,$$
où l'on note $(R_i,R_n)$ le plus petit sous-groupe parabolique contenant $R_i$ et $R_n$. On a  $(R_i,R_n)\subset Q$. L'hypothèse de récurrence implique que chaque somme dans le membre de droite est nulle. D'où le résultat.
\end{preuve}

 Soit $Q$ le sous-groupe parabolique maximal de $G$ qui contient $P$ défini par la condition $\Delta -\Delta^Q=\{\al\}$. Soit $X\in \nc^G\cap \pgo(F)$ tel que $\xi_{\overline{\oc}}^{P,G}(X)\not=0$. D'après le lemme \ref{lem:annul}, il existe un sous-groupe parabolique $R$ tel que
$$ P\subset R \not\subset Q$$
et qui est un élément minimal de l'ensemble (\ref{eq:lens}).  Soit $S=R\cap Q$.

\begin{lemme}\label{lem:non-induit}
 On a  $S \subsetneq R$. De plus, si l'on décompose 
$$X=X_S+U+V$$
avec $X_S\in \mgo_S(F)$, $U\in (\mgo_R\cap\ngo_S)(F)$ et $V\in \ngo_R(F)$, on a 
$$X_S+U\notin I_{S}^{R}(X_S).$$
  \end{lemme}

  \begin{preuve} Si $S=R$ alors $R\subset Q$ ce qui n'est pas. Supposons en raisonnant par l'absurde qu'on a $X_S+U\in I_{S}^{R}(X_S)$. La transitivité de l'induction (cf. lemme \ref{lem:transitivite}) implique qu'on a 
$$I_S^G(X)=I_S^G(X_S)= I_R^G(I_S^R(X_S))= I_R^G(X_S+U)=I_R^G(X).$$
Or on a $I_R^G(X)\subset \overline{\oc}$ puisque $R$ appartient à l'ensemble (\ref{eq:lens}). On a ainsi prouvé qu'on a $I_S^G(X)\subset \overline{\oc}$ donc que $S$ appartient à l'ensemble (\ref{eq:lens}). Mais  ceci contredit la minimalité de $R$.
  \end{preuve}

On observe que la fonction $|\xi_{\overline{\oc}}^{P,G}(X)|$  est toujours majorée par le nombre de sous-groupes paraboliques compris entre $P$ et $G$.  Il résulte alors des lemmes \ref{lem:annul} et \ref{lem:non-induit} que pour prouver la proposition \ref{prop:reduc4}, il suffit de majorer pour  tout sous-groupe parabolique $R\not\subset Q$  qui contient $P$  et toute orbite nilpotente $\oc^S \in (\nc^{M_S})$ (avec  $S=R\cap Q$) l'expression  suivante 

\begin{equation}
  \label{eq:XUV1}
   q^{-\bg2\rho_P^S,H_P(a)) \bd } \sum_{X\in \oc^S \cap \pgo(F)} \mathbf{1}_D(a^{-1}X a) \times
 \end{equation}
\begin{equation}
  \label{eq:XUV2}
   q^{-\bg2\rho_S^R,H_S(a)) \bd } \sum_{U} \mathbf{1}_D(a^{-1}U a),
 \end{equation}
 où la somme est prise sur les $U\in (\mgo_R\cap\ngo_S)(F)$ tels  que  $ X+U\notin I_S^R(X)$, multipliée par  
\begin{equation}
  \label{eq:XUV3}
   q^{-\bg2\rho_R,H_R(a)) \bd } \sum_{V\in \ngo_R(F)} \mathbf{1}_D(a^{-1}V a).
 \end{equation}

Le gain recherché du facteur  $q^{-\bg \al, H(a)\bd }$ va être apporté par l'expression \eqref{eq:XUV2}. Pour cela, on fixe $X\in \oc^S$ et on va majorer \eqref{eq:XUV2}. L'induite $I_S^R(\oc^S)$ est l'unique orbite  dans $(\nc^{M_R})$ telle que l'intersection
$$I_S^R(\oc^S)\cap \big(\oc^S \oplus (\ngo_S\cap\mgo_R)\big)$$
soit un ouvert dense. Il existe un nombre fini de polynôme disons $f_1,\ldots,f_r$ dans $F[\sgo\cap \mgo_R]$ tels que 
$$\big(\oc^S\oplus (\ngo_S\cap\mgo_R)\big) \cap \vc,$$
pour $\vc=\vc(f_1,\ldots,f_r)$, en soit le fermé complémentaire.

\begin{lemme} \label{lem:nonnuls} Pour tout $X\in \oc^S$, il existe $1\leq i\leq r$ tel que le polynôme $f_i(X+\cdot)\in F[ \ngo_S\cap\mgo_R ]$ soit non nul.
\end{lemme}

\begin{preuve}
  Supposons l'assertion en défaut. Il existe $X\in \oc^S$ tel que pour tout $1\leq i\leq r$  le polynôme  $f_i(X+\cdot)\in F[ \ngo_S\cap\mgo_R ]$ est nul. Prenons un point $Y\in I_S^R(\oc^S)$. Quitte à conjuguer $X$ par un élément de $M_S$, on peut supposer que $Y=X+U$ avec $U\in (\ngo_S\cap\mgo_R)(F)$. Par conséquent, pour tout $1\leq i\leq r$ on a  $f_i(Y)=f_i(X+U)=0$. Donc $Y\in \vc(f_1,\ldots,f_r)$ ce qui est la contradiction cherchée.
\end{preuve}

Pour $X\in \oc^S$ fixé,  soit $\vc_X=\vc(f_1(X+\cdot),\ldots,f_r(X+\cdot))$. Les polynômes $f_i(X+\cdot)$ sont non tous nuls d'après le lemme \ref{lem:nonnuls} ci-dessus et leurs degrés partiels sont bornés par ceux des $f_i$, donc indépendamment de $X$. L'expression \eqref{eq:XUV2} se réécrit donc
 $$q^{-\bg2\rho_S^R,H_S(a)) \bd } \sum_{U\in (\mgo_R\cap\ngo_S)(F)\cap \vc_X} \mathbf{1}_D(a^{-1}U a).$$
Par la proposition \ref{prop:maj-var}, cette expression est bornée, à une constante près qui ne dépend pas de $X$, par $q^{-\bg \al, H(a)\bd }$ lorsque $a\in A^G(T_1)$.

Pour conclure, il suffit de prouver que, pour $a\in A^G(T_1)$, les expressions \eqref{eq:XUV1} et  \eqref{eq:XUV3} sont majorées par une constante qui ne dépend que de $D$ et $T_1$. Par exemple, il existe $c\geq 1$ tel que l'expression \eqref{eq:XUV1} se majore à l'aide du lemme \ref{lem:maj} par
\begin{eqnarray*}
   q^{-\bg2\rho_P^S,H_P(a)) \bd } \prod_{\beta \in \Sigma^{N_P\cap M_S}} \sum_{X\in F}\mathbf{1}_D(\beta(a)^{-1}X)&\leq &c \cdot \prod_{\beta \in \Sigma^{N_P\cap M_S}} \sup(q^{\deg(D)}, q^{-\bg \beta, H(a)\bd })\\
   &\leq & c \cdot \prod_{\beta \in \Sigma^{N_P\cap M_S}} \sup(q^{\deg(D)}, q^{-\bg \beta, T_1\bd }).
\end{eqnarray*}
L'expression \eqref{eq:XUV3} se traite de la même façon.
\end{preuve}

\end{paragr}

\begin{paragr}[Démonstration de la proposition \ref{prop:reduc1}.] ---\label{S:preuvereduc1}
  Si l'on fait la synthèse du corollaire \ref{cor:step5}, de la proposition \ref{prop:2ndemaj} et du corollaire \ref{cor:reduc4}, on voit qu'il suffit de majorer pour  $T\in a_B^+$ l'expression suivante
$$\sum_{a\in A^{1,2}_{Q,e}(T)} q^{-\bg2\rho_Q^1,H_Q(a) \bd } \prod_{\al\in \Delta^2-\Delta^1} q^{-r \bg\al,H_B(a)\bd} 
  \cdot \sum_{X\in \mathcal{N}^{M_1}} \mathbf{1}_D(a^{-1}Xa),$$
où $D$ est un diviseur sur $C$, $Q\in \fc^{P_1}(B)$ et $r>0$, les autres notations sont celles des §§\ref{S:reduc2} et \ref{S:reduc3}.
 On majore tout d'abord la somme intérieure ainsi :
 \begin{eqnarray*}
   \sum_{X\in \mathcal{N}^{M_1}} \mathbf{1}_D(a^{-1}Xa) &\leq& \sum_{X\in \mgo_1(F)} \mathbf{1}_D(a^{-1}Xa)\\
&=& \sum_{X\in \mgo_1(F)\cap \qgo(F)} \mathbf{1}_D(a^{-1}Xa) \times  \sum_{X\in \mgo_1(F)\cap \ngo_{\bar{Q}}(F)} \mathbf{1}_D(a^{-1}Xa),
 \end{eqnarray*}
où $\bar{Q}$ est le sous-groupe parabolique opposé à $Q$. Pour tout $a\in A^{M_1}(T_1)$, on a déjà vu (par exemple à la toute fin de la preuve de la proposition \ref{prop:reduc4}) que l'expression 
$$ q^{-\bg2\rho_Q^1,H_Q(a) \bd } \sum_{X\in \mgo_1(F)\cap \qgo(F)} \mathbf{1}_D(a^{-1}Xa)$$
est bornée. Par le lemme \ref{lem:maj}, il existe une constante $c>0$ telle que
$$\sum_{X\in \mgo_1(F)\cap \ngo_{\bar{Q}}(F)} \mathbf{1}_D(a^{-1}Xa) \leq c\cdot \prod_{\al \in \Sigma^{M_1\cap N_Q}} \sup(q^{\deg(D)-\bg \al,H_B(a)\bd},1).$$
Pour $a\in A^{M_1}(T_1)$, l'expression du membre de droite ci-dessus est également majorée. On est donc ramené à majorer l'expression suivante 

$$\sum_{a\in A^{1,2}_{Q,e}(T)} \prod_{\al\in \Delta^2-\Delta^1} q^{-r \bg\al,H_B(a)\bd} $$
 
Soit $H=H_B(a)$ pour $a\in A^{1,2}_{Q,e}(T)$. On écrit 
$$H=H^1+H_1^2+H_2^G+H_G$$
suivant la décomposition $a_B=a_B^{P_1}\oplus a_{P_1}^{P_2}\oplus a_{P_2}^G \oplus a_G$. De même, on écrit $T=T^1+T_1^2+T_2^G+T_G$. La composante $H_G$ ne dépend pas de $H$ (seulement de $e$). La composante sur $a_B^{P_1}$ s'écrit
$$H^1=\sum_{\beta \in \Delta^{P_1}} \bg \varpi_\beta,H \bd \beta^\vee.$$
Comme $a$ appartient en particulier à $A^{M_1}(T_1)$, on a $\bg \beta,H \bd >\bg \beta,T_1 \bd$ pour tout $\beta\in \Delta^{P_1}$. \emph{A fortiori}, on a  $\bg \varpi_\beta,H \bd >\bg \varpi_\beta,T_1 \bd$ pour tout $\beta\in \Delta^{P_1}$. Par ailleurs, comme $a$ appartient à $A^{P_1,P_2}(T_1,T)$, on a aussi  $\bg \varpi_\beta,H \bd \leq \bg \varpi_\beta,T \bd$ pour tout $\beta\in \Delta^{P_1}$. Il s'ensuit que la composante $H^1$ vit dans un compact : il n'y en a donc qu'un nombre fini possible. De plus, pour tout $\al\in \Delta^2-\Delta^1$ et  tout $\beta\in \Delta^{P_1}$, on  a   $\bg\al,\beta^\vee\bd\leq 0$ et donc 
\begin{eqnarray*}
  - \bg\al,H^1\bd &=& -\sum_{\beta \in \Delta^{P_1}} \bg \varpi_\beta,H \bd  \bg\al,\beta^\vee\bd\\
&\leq & -\sum_{\beta \in \Delta^{P_1}} \bg \varpi_\beta,T \bd  \bg\al,\beta^\vee\bd=   - \bg\al,T^1\bd.
\end{eqnarray*}

Il existe donc une constante $c>0$ telle que pour tout $T\in \overline{a_B^+}$, on ait
\begin{equation}
  \label{eq:1comp}
  \sum_{H^1} \prod_{\al\in \Delta^2-\Delta^1} q^{-r \bg\al,H^1\bd} \leq c\cdot \prod_{\beta \in \Delta^{P_1} } (1+| \bg \varpi_\beta,T \bd  |) \cdot \prod_{\al\in \Delta^2-\Delta^1}  q^{-r \bg\al,T^1\bd},
\end{equation}
où la somme est prise sur l'ensemble (fini) des composantes $H^1$ possibles.

Fixons $H_1^2$ et regardons les contraintes qui pèsent sur $H_2^G$. Ce sont celles qui traduisent les conditions 4 et 5 (cf. les conditions qui précèdent les remarques \ref{rq:condi}). On a donc pour tout $\al\in \Delta_1-\Delta_1^2$
$$\bg \al, H_2^G \bd \leq \bg \al, T-H_1^2\bd,$$
qui traduit l'inégalité  $\bg \al, H \bd \leq \bg \al, T\bd$.
Par ailleurs, pour tout $\varpi\in \hat{\Delta}_2$, on a 
$$\bg \varpi,H_2^G\bd = \bg \varpi,H \bd >  \bg \varpi, T\bd.$$
La composante $H_2^G$ est astreinte à rester dans un polyèdre compact qui dépend de $H_1^2$ et de $T$. En particulier, le nombre de tels $H_2^G$ est borné par $|P(T,H_1^2)|$ pour  un certain polynôme $P$.

Retenons des conditions sur $H_1^2$ la suivante (c'est la condition 3  qui précède les remarques \ref{rq:condi}) : pour tout $\al\in \Delta^2-\Delta^1$, on a $\bg \al, H_1^2 \bd > \bg \al, T_1^2 \bd $. On observera qu'on a $\bg \al, T_1^2 \bd\geq 0$ pour $T\in \overline{a_B^+}$. En rassemblant les majorations précédentes, on voit qu'il existe $c>0$ tel que
$$ \sum_{a\in A^{1,2}_{Q,e}(T)} \prod_{\al\in \Delta^2-\Delta^1} q^{-r \bg\al,H_B(a)\bd} \leq$$

 $$ c\cdot \prod_{\beta \in \Delta^{P_1} } (1+| \bg \varpi_\beta,T \bd  |) \cdot \prod_{\al\in \Delta^2-\Delta^1}  q^{-r \bg\al,T^1\bd} \sum_{H_1^2} |P(T,H_1^2)|  \prod_{\al\in \Delta^2-\Delta^1}  q^{-r \bg\al,H_1^2\bd},$$
où l'on somme sur l'ensemble des composantes $H_1^2$ possibles. La somme sur $H_1^2$ porte sur un réseau intersecté avec le cône  défini  $\bg \al, H_1^2 \bd > \bg \al, T_1^2 \bd $  pour tout $\al\in \Delta^2-\Delta^1$. On la majore alors par  $|Q(T)|\cdot \prod_{\al\in \Delta^2-\Delta^1}  q^{-r \bg\al,T_1^2\bd}$ pour un certain polynôme $Q$. Finalement, on a prouvé qu'il existe un polynôme $P$ tel que

$$ \sum_{a\in A^{1,2}_{Q,e}(T)} \prod_{\al\in \Delta^2-\Delta^1} q^{-r \bg\al,H_B(a)\bd} \leq  |P(T)| \prod_{\al\in \Delta^2-\Delta^1}  q^{-r \bg\al,T\bd}  . $$
Le résultat est alors évident.
\end{paragr}

\section{Le calcul d'une intégrale nilpotente \og régulière par blocs\fg{}}\label{sec:reg-blocs}

\begin{paragr}[Énoncé du résultat principal.] ---  On continue avec les notations des sections précédentes. Soit $n\geq 1$ un entier et $G=GL(n)$ sur le corps fini $\Fq$. Soit   $B_0\subset G$ et  $T_0\subset G$ respectivement le sous-groupe de Borel des matrices triangulaires supérieures et le sous-tore maximal diagonal. On identifie le groupe $X^*(T_0)$ des caractères de $T_0$ à $\ZZ^n$ de la manière suivante : le caractère donné par la $i$-ème entrée correspond au $i$-ème vecteur de la base canonique de $\ZZ^n$. Soit $\Gc=GL(n,\CC)$. Soit  $\Bc_0\subset \Gc$ et  $\Tc_0\subset \Gc$ respectivement le sous-groupe de Borel des matrices triangulaires supérieures et le sous-tore maximal diagonal . Comme précédemment, on identifie le groupe $X^*(\Tc_0)$ des caractères de $\Tc_0$ à $\ZZ^n$. Dualement, on a donc une identification des groupes des caractères $X_*(T_0)$ et $X_*(\Tc_0)$ avec $\ZZ^n$. On en déduit des isomorphismes
$$X^*(T_0)\simeq X_*(\Tc_0) \text{  et  } X_*(T_0)\simeq X^*(\Tc_0)$$
qui envoient l'ensemble $\Delta$ des racines simples (resp. $\Delta^\vee$ des coracines simples) de $T_0$ dans $B_0$ dans celui des coracines simples (resp. des racines) de $\Tc_0$ dans $\Bc_0$.

Soit $G'=SL(n)$ le groupe dérivé de $G$ et $T_0'=T_0\cap G'$. On a une suite exacte courte entre groupes de caractères  dont les flèches intermédiaires sont données par les restrictions
$$0\longrightarrow X^*(G) \longrightarrow X^*(T_0) \longrightarrow  X^*(T_0')\longrightarrow 0.$$
Les groupes ci-dessus sont des $\ZZ$-modules libres et, en appliquant le foncteur $\Hom_\ZZ(\cdot,\ZZ)$, on obtient une nouvelle suite exacte 
$$0\longrightarrow X_*(T_0')  \longrightarrow X_*(T_0) \longrightarrow  \ago_G\longrightarrow 0,$$
où $\ago_G=\Hom_\ZZ(X^*(G),\ZZ)$ et $X_*(\cdot)$ désigne le groupe des cocaractères. Le groupe $X_*(T_0')$ n'est autre que le sous-groupe $\ZZ(\Delta^\vee)$  de $ X_*(T_0) $ engendré par l'ensemble $\Delta^\vee$ des coracines simples. En utilisant l'isomorphisme $X_*(T_0)\simeq X^*(\Tc_0)$ qu'on a fixé, on obtient une identification 
$$\ago_G\simeq X^*(\Tc_0)/ \ZZ(\Delta^\vee)$$
où l'on interprète maintenant $\Delta^\vee$ comme l'ensemble des racines simples de $\Tc_0$ dans $\Gc$. Soit $Z_{\Gc}\subset \Tc_0$ le centre du groupe $\Gc$. On a donc un accouplement naturel
$$Z_{\Gc} \times \ago_G \to \CC^\times$$
Le groupe $X^*(G)$ a un générateur canonique à savoir le déterminant ce qui permet d'identifier $\ago_G$ à $\ZZ$. On a donc aussi un accouplement
$$Z_{\Gc} \times \ZZ \to \CC^\times$$
qu'on note
$$(z,e)\mapsto z^e.$$
  
Soit $d$ et $r$ des entiers $\geq 1$ tels que $n=rd$. Soit $\Pc\subset \Gc$ le sous-groupe parabolique dont le facteur de Levi standard $\Mc$ est   isomorphe à $GL(d)^r$. Soit $\Mc'=\Mc\cap \Gc'$ et $Z_{\Mc'}^0$ la composante neutre du centre de $\Mc'$. C'est un tore défini sur $\CC$ de rang $r-1$.

Soit $Z_C$ la fonction zêta de la courbe $C$ définie par la série formelle
$$Z_C(t)=\exp(\sum_{k=1}^\infty |C(\mathbb{F}_{q^k})|t^k/k),$$
où $\mathbb{F}_{q^k}$ désigne un corps fini à $q^k$ éléments. Il est bien connu que $Z_C$ est en fait une fraction rationnelle en $t$ de dénominateur  $(1-t)(1-qt)$. On a d'ailleurs 
$$\zeta(s)=Z_C(q^{-s})$$
avec les notations du §\ref{S:calculI-adeles}.
On introduit aussi la fonction
$$Z^*_C(t)=(1-qt)Z_C(t)$$
qui n'a pas de pôle en $t=q^{-1}$.

Soit $D$ un diviseur sur $C$ et
$$Z_{C,D}^d(t)=t^{-d\deg(D)}Z_C(q^{-1}t)Z_C(q^{-2}t)\ldots Z_C(q^{-d}t).$$
Soit la fraction rationnelle sur $\Tc_0'$ définie pour $t\in \Tc_0'$ par   
$$\Phi_{C,D}^d(t)=\prod_{\varpi \in \hat{\Delta}_{\Pc}} Z_{C,D}^d(t^{\varpi})$$
où  $\hat{\Delta}_{\Pc}$ désigne l'ensemble des poids  fondamentaux $\Pc$-dominants  de $\Tc_0'$ : ce sont des caractères de $\Tc_0'$. Lorsque $\Pc=\Gc$ (i.e. lorsque $d=n$), l'ensemble $\hat{\Delta}_{\Gc}$ est vide et par convention, $\Phi_{C,D}^n(t)=1$. Hormis le cas $\Pc=\Gc$, la fraction rationnelle $\Phi_{C,D}^d(t)$ a un pôle en $t=1$.

Soit $W_{\Mc}$ le groupe de Weyl de $\Mc$ : c'est le quotient du sous-groupe de $\Gc$ qui normalise à la fois $\Tc_0$ et $\Mc$ par le normalisateur de $\Tc_0$ dans $\Mc$. Le groupe $W_{\Mc}$ a pour ordre $r!$ et il agit sur $Z_{\Mc}$. On note $w\cdot t$ l'action de $w\in W_{\Mc}$ sur un élément $t\in Z_{\Mc}$. On observera que cette action est triviale sur le sous-groupe $Z_{\Gc}$. Pour tout $t\in Z_{\Mc}$ et tout $z\in Z_{\Gc}$, on a donc
$$w\cdot(tz)=(w\cdot t) z$$
ce qu'on note simplement $ w\cdot t z$.
Pour tout  $e\in \ZZ$ et $t\in Z_{\Mc'}^0$ soit
$$\Psi_{C,D}^{d,e}(t)=\frac{1}{n\cdot |W_{\Mc}|} \sum_{z\in Z_{\Gc'}} \sum_{w\in W_{\Mc}}z^{-ed}  \Phi_{C,D}^d( w\cdot t z).$$
Cela définit une fraction rationnelle sur $Z_{\Mc'}^0$. Le groupe $Z_{\Gc'}$ est cyclique d'ordre $n$. En particulier, pour tout $z\in Z_{\Gc'}$, le nombre complexe $z^{-ed}$ est une racine $r$-ième de l'unité qui ne dépend que de la classe de $e$ modulo $r$. Lorsque $d=n$, on a $\Psi_{C,D}^{n,e}(t)=1$ pour tout $t$.

\begin{theoreme}\label{cor:calcul}
  \begin{enumerate}
  \item La fraction rationnelle $\Psi_{C,D}^{d,e}(t)$ n'a pas de pôle en $t=1$.
  \item Soit $\oc\subset \ggo(F)$ l'orbite nilpotente des éléments dont la décomposition de Jordan est formée de $d$ blocs de taille $r$. Soit $ K^G_{D,0,\oc}$  la fonction sur $G(\AAA)$ définie en \eqref{eq:KDTo} pour le paramètre $T=0$.

Lorsque $e$ est premier à $r$, l'intégrale ci-dessous, qui est absolument convergente par le corollaire \ref{cor:approx},
$$\int_{G(F)\back G(\AAA)^e}   K^G_{D,0,\oc}(g)\, dg$$
est égale à 
$$q^{n(n-d)\deg(D)/2} q^{nd(g_C-1)} Z_C^*(q^{-1})Z_C(q^{-2})\ldots Z_C(q^{-d}) \Psi_{C,D}^{d,e}(1).$$
  \end{enumerate}
\end{theoreme}

\begin{remarque}\label{rq:calcul}
  Dans l'expression ci-dessous, on a $n(n-d)=\dim(\oc)$. L'intégrale et l'expression ci-dessus ne dépendent que de $e$ modulo $r$. Lorsque $d=n$, (cas de l'orbite nulle), l'intégrande est identiquement $1$ et l'intégrale n'est autre que le volume du quotient $G(F)\back G(\AAA)^e$. L'expression ci-dessus se réduit à 
$$ q^{n^2(g_C-1)} Z_C^*(q^{-1})Z_C(q^{-2})\ldots Z_C(q^{-n})$$
qui est la formule de Siegel pour le volume. La démonstration du théorème \ref{cor:calcul} utilise d'ailleurs cette formule sans la redémontrer.
\end{remarque}

\begin{preuve}
  Elle se trouve au paragraphe \ref{par:preuveduthmcorcalcul}.
\end{preuve}
\end{paragr}

\begin{paragr}[L'élément $X\in \oc$.] --- On continue avec les entiers $n=dr$ et les notations sont celles de la section \ref{sec:calculI}. Rappelons qu'on note $\oc\subset \ggo(F)$ l'orbite  des éléments nilpotents qui possèdent $d$ blocs de Jordan de taille $r$. C'est encore l'orbite de l'élément $X$ défini au §\ref{S:calculI-nota}. Pour tout sous-groupe parabolique standard $P$ de $G$, il existe une orbite $\oc^P\in (\nc^{M_P})$ tel que $I_P^G(\oc^{P})=\oc$ si et seulement si $P_0\subset P$. Dans ce cas, l'orbite $\oc^P$ est uniquement définie. On a d'ailleurs $\oc^P=I_{P_0}^P(0)$. Avec les notations de §\ref{S:calculI-nota}, c'est aussi l'orbite de $X^P$.
\end{paragr}

\begin{paragr}[Fonctions $\tilde{K}^P_{D,X}$.] --- \label{S:Ktilde} Pour tout sous-groupe parabolique standard $P$ de décomposition de Levi $P=MN$, on pose $\tilde{K}^P_{D,X}(g)=0$ sauf si $P_0\subset P$ auquel cas on définit
  \begin{equation}
    \label{eq:KtildePDX}
    \tilde{K}^P_{D,X}(g)=\sum_{\delta\in M_{X^P}(F)N(F) \back P_0(F)} \int_{\ngo(\AAA)} \mathbf{1}_D((\delta g)^{-1}(X^P+U)\delta g)\, dU.
  \end{equation}
  On notera que, par construction, $\tilde{K}^P_{D,X}$ est une fonction sur le quotient $P_0(F)\back G(\AAA)$. Rappelons que le groupe $M_{X^P}$ est le centralisateur dans $M$  de l'élément $X^P\in \mgo(\Fq)$ défini au \ref{S:calculI-nota}. On a $M_{X^P}\subset M\cap P_0$.
Avec la notation \eqref{eq:KPDX} du §\ref{S:KPDX}, on a 
\begin{equation}
  \label{eq:KKtilde0}
  \tilde{K}^P_{D,X}(g)=\sum_{\delta\in M_{X^P}(F)N(F) \back P_0(F)} K^P_{D,X}(\delta g).
\end{equation}
Soit $\oc$ l'orbite nilpotente $I_{P_0}^G(0)$.  Pour faire le lien avec les objets du §\ref{S:asymp}, on observera la relation suivante :
\begin{equation}
  \label{eq:KKtilde}
   K^P_{D,\oc}(g)=\sum_{\delta\in P_0(F)\back P(F)}\tilde{K}^P_{D,X}(\delta g)
 \end{equation}
d'où l'on déduit pour $T\in a_B$ et $g\in G(\AAA)$
 \begin{equation}
  \label{eq:KKtildeII}
K^G_{D,T,\oc}(g)=\sum_{\delta\in P_0(F)\back P(F)} \big[\sum_{B\subset P\subset G} (-1)^{\dim(a_P^G)} \hat{\tau}_P(H_P(\delta g)-T) \tilde{K}^P_{D,X}(\delta g) \big].
\end{equation}

Le théorème suivant sera utilisé dans la démonstration du théorème \ref{cor:calcul}.

\begin{theoreme}\label{thm:cvabs}
Soit $e\in \ZZ$. L'intégrale suivante est finie pour tout $T\in a_B$
$$\int_{P_0(F)\back G(\AAA)^e} |\sum_{B\subset P\subset G} (-1)^{\dim(a_P^G)} \hat{\tau}_P(H_P(g)-T) \tilde{K}^P_{D,X}(g) |\,dg<\infty$$
     \end{theoreme}
  
     \begin{preuve} Elle occupe entièrement le  paragraphe \ref{S:preuve-cvabs} ci-dessous.
     \end{preuve}
\end{paragr}

\begin{paragr}[Démonstration du théorème \ref{cor:calcul}.] --- \label{par:preuveduthmcorcalcul}Par le  théorème \ref{thm:cvabs} ci-dessous et la relation \eqref{eq:KKtildeII} ci-dessus, on  a l'égalité
$$\int_{G(F)\back G(\AAA)^e}   K^G_{D,0,\oc}(g)\, dg=$$
\begin{equation}
  \label{eq:integ-alternee}
  \int_{P_0(F)\back G(\AAA)^e} \sum_{B\subset P\subset G} (-1)^{\dim(a_P^G)} \hat{\tau}_P(H_P(g)) \tilde{K}^P_{D,X}(g) \,dg.
\end{equation}
Bien sûr, dans l'intégrande ci-dessus, on peut et on va imposer dans la somme sur $P$ de ne prendre que les sous-groupes paraboliques contenant $P_0$. On peut montrer qu'il existe une constante $c$ telle que si $g\in G(\AAA)$ est tel que 
$$\sum_{P_0\subset P\subset G} (-1)^{\dim(a_P^G)} \hat{\tau}_P(H_P(g)) \tilde{K}^P_{D,X}(g)\not=0$$
alors  pour tout $\varpi\in \hat{\Delta}_{P_0}$, on  a $\bg \varpi,H_{P_0}(g)\bd \geq c$ (en fait cette propriété est déjà vraie pour chaque terme de la somme ; pour le voir il suffit d'expliciter la fonction $\tilde{K}^P_{D,X}(g)$, cf. lemme \ref{lem:KPX} et la démonstration du lemme \ref{lem:calculIc}).

À l'aide du théorème \ref{thm:cvabs}, on en déduit que pour tout $\la$ dans l'ouvert $\Omega$ de $a_{P_0}^*$ défini par  $\Re(\bg \la,\al^\vee\bd)>0$ pour tout $\al\in \Delta_{P_0}$, l'intégrale
\begin{equation}
  \label{eq:integrale-lambda}
  \int_{P_0(F)\back G(\AAA)^e} q^{-\bg \la,H_{P_0}(g) \bd}\sum_{P_0\subset P\subset G} (-1)^{\dim(a_P^G)} \hat{\tau}_P(H_P(g)) \tilde{K}^P_{D,X}(g) \,dg
\end{equation}
est absolument convergente et que sa limite quand $\la$ tend vers $0$ est égale à l'intégrale \eqref{eq:integ-alternee}. Par ailleurs, pour $\la$ dans l'ouvert $\Omega$, l'intégrale \eqref{eq:integrale-lambda} se calcule terme à terme ainsi
$$\sum_{P_0\subset P\subset G} (-1)^{\dim(a_P^G)} \int_{P_0(F)\back G(\AAA)^e} q^{-\bg \la,H_{P_0}(g) \bd} \hat{\tau}_P(H_P(g)) \tilde{K}^P_{D,X}(g) \,dg.$$
L'interversion se justifie car, comme on va le voir, chaque intégrale ci-dessus est absolument convergente. Pour le voir, on écrit à l'aide de \eqref{eq:KKtilde0}
$$\int_{P_0(F)\back G(\AAA)^e} q^{-\bg \la,H_{P_0}(g) \bd} \hat{\tau}_P(H_P(g)) \tilde{K}^P_{D,X}(g) \,dg$$
$$=\int_{ M_{X^P}(F)N_P(F)(F)\back G(\AAA)^e} q^{-\bg \la,H_{P_0}(g) \bd} \hat{\tau}_P(H_P(g)) K^P_{D,X}(g) \,dg$$
$$=\sum_{H\in H_{P_0}(G(\AAA)^e)} \hat{\tau}_P(H)\cdot  \hat{I}^{P}_{D,X}(H)\cdot   q^{-\bg \la,H\bd }.$$
Le coefficient  $\hat{I}^{P}_{D,X}(H)$ est celui défini en \eqref{eq:IcPDX}. En particulier, on reconnaît dans la somme ci-dessus la série  $I^{0}_{P,D}(X,\la/d)$ définie en \eqref{eq:laserieIxiPD}, dont on connaît la convergence absolue sur $\Omega$ par la proposition \ref{prop:lecalculpourP}. On en déduit que l'intégrale \eqref{eq:integrale-lambda} est égale à la série $I^0_D(X,\la/d)$ définie en \eqref{eq:laserieIxiD}.
Le théorème est alors une simple traduction du théorème \ref{cor:sommezeta2} qui tient compte de la formule de Siegel rappelée dans la remarque \ref{rq:calcul}. Il suffit d'identifier par l'exponentielle le tore $Z_{\Mc'}^0$ à un quotient de $a_M^{G,*}$. La somme sur $\pc(M)$ dans  le théorème \ref{cor:sommezeta2} est remplacée par la somme sur $w\in W_{\Mc}$. Il existe un générateur $z$ de $Z_{\Gc'}$ tel qu'on ait 
$$z^{-d}=\exp(-2\pi i /r)$$
et les accouplements de $z$ et $\gamma$ avec les  poids $\Pc$-fondamentaux soient égaux. Ces accouplements sont des racines $r$-ièmes de l'unités. Ce générateur est bien défini modulo le sous-groupe d'ordre $d$ de  $Z_{\Gc'}$. On peut alors remplacer la somme sur $k$ par la somme des $r$ premières puissances de $z$, puis, quitte à diviser par $n$ au lieu de $r$, par la somme sur tous les éléments de $Z_{\Gc'}$. Le résultat s'ensuit.
\end{paragr}

\begin{paragr}[Preuve du théorème \ref{thm:cvabs}.] ---\label{S:preuve-cvabs} Comme dans la preuve du  corollaire \ref{cor:approx}, on a l'égalité entre

$$\sum_{B\subset P\subset G} (-1)^{\dim(a_P^G)} \hat{\tau}_P(H_P(g)-T) \tilde{K}^P_{D,X}(g)$$
et
$$\sum_{B\subset Q\subset G} \Gamma'_Q(H_Q(g),T) \sum_{B\subset P\subset Q} (-1)^{\dim(a_P^Q)} \hat{\tau}_P^Q(H_P(g)) \tilde{K}^P_{D,X}(g) \ ;$$
la fonction $\Gamma'_Q$ ci-dessus est celle qui apparaît dans \eqref{eq:def-Gamma}. On se restreint dans la somme ci-dessus aux sous-groupes paraboliques $Q$ qui contiennent $P_0$ sans quoi leur contribution est nulle. On a la formule suivante de descente  pour $P_0\subset P\subset Q$, $m\in M_Q(\AAA)$ et $n\in N_Q(\AAA)$
\begin{eqnarray*}
   \tilde{K}^{P}_{D,X}(n m)&=& q^{\bg2\rho_Q,H_Q(m) \bd} q^{\dim(N_Q)(1-g_C+\deg(D))} \tilde{K}^{P\cap M_Q}_{D,X}(m).
 \end{eqnarray*}
Ici $ \tilde{K}^{M_Q\cap P}_{D,X}$ est l'analogue évident de \eqref{eq:KtildePDX} pour le groupe $M_Q$. Il s'ensuit qu'on a pour $P_0\subset Q$ par décomposition d'Iwasawa
$$\int_{P_0(F)\back G(\AAA)^e}  \Gamma'_Q(H_Q(g),T) |\sum_{B\subset P\subset Q} (-1)^{\dim(a_P^Q)} \hat{\tau}_P^Q(H_P(g)) \tilde{K}^P_{D,X}(g)| \, dg$$
 $$= q^{\dim(N_Q)\deg(D)} \int_{(M_Q\cap P_0)(F)\back M_Q(\AAA)\cap G(\AAA)^e} |\sum_{B \subset P\subset Q} (-1)^{\dim(a_P^Q)} \hat{\tau}_P^Q(H_P(m)) \tilde{K}^{M_Q\cap P}_{D,X}(m)| \, dm.$$
Si $Q\subsetneq G$, alors $M_Q$ est un produit de groupes linéaires de rang strictement plus petit que celui de $G$. L'intégrale ci-dessus est aussi un produit. En raisonnant par récurrence comme dans la preuve du corollaire  \ref{cor:approx}, on peut donc supposer que l'intégrale 
$$\int_{P_0(F)\back G(\AAA)^e}  \Gamma'_Q(H_Q(g),T) |\sum_{B\subset P\subset Q} (-1)^{\dim(a_P^Q)} \hat{\tau}_P^Q(H_P(g)) \tilde{K}^P_{D,X}(g)| \, dg$$
est absolument convergente. Il suffit alors de prouver le théorème \ref{thm:cvabs}  pour un seul élément $T$, par exemple un élément  assez profond dans la chambre de Weyl positive ce que l'on suppose dans la suite. 

  \begin{lemme}
\label{lem:Ireecr}
    Pour tout $T\in \overline{a_B^+}$, on a l'égalité entre
$$\sum_{B\subset P\subset G} (-1)^{\dim(a_P^G)} \hat{\tau}_P(H_P(g)-T) \tilde{K}^P_{D,X}(g)$$
et
$$\sum_{B\subset P_1\subset P_2\subset G} \sum_{\delta\in P_1(F)\back G(F)} \chi_T^{P_1,P_2}(\delta g) \sum_{\{P\in \fc^{P_2}(P_1) \mid \delta\in P(F)\}}  (-1)^{\dim(a_P^G)} \tilde{K}^P_{D,X}(g).$$
  \end{lemme}

  \begin{preuve}
On utilise \eqref{eq:HN} pour avoir
$$\sum_{B\subset P\subset G} (-1)^{\dim(a_P^G)} \hat{\tau}_P(H_P(g)-T) \tilde{K}^P_{D,X}(g)=$$
$$=\sum_{B\subset P\subset G} (-1)^{\dim(a_P^G)} \big[ \sum_{B\subset P_1\subset P}  \sum_{\delta\in P_1(F)\back P(F)}F^{P_1}(\delta  g,T) \tau_{P_1}^{P}(H_{P_1}(\delta  g) -T)    \big] \hat{\tau}_P(H_P( g)-T)\tilde{K}^P_{D,X}(g)$$
On a $ \hat{\tau}_P(H_P( g)-T)= \hat{\tau}_P(H_B(\delta g)-T)$ pour $\delta\in P(F)$. Puis on utilise la formule d'inversion d'Arthur $\tau_{P_1}^{P} \cdot  \hat{\tau}_P=\sum_{ P \subset P_2\subset G}  \sigma_{P_1}^{P_2}$ (cf. \cite{ar-intro} formule (8.2)). On obtient
$$ =\sum_{B\subset P_1\subset P\subset P_2 \subset G} (-1)^{\dim(a_P^G)} \sum_{\delta\in P_1(F)\back P(F)} F^{P_1}(\delta g,T)\sigma_{P_1}^{P_2}(H_{P_1}(\delta g)-T) \tilde{K}^P_{D,X}(g).$$
Pour obtenir le lemme, il suffit d'inverser les sommes sur $P$ et sur $\delta$.
  \end{preuve}

En tenant compte du lemme \ref{lem:Ireecr} ci-dessus, on voit que le théorème \ref{thm:cvabs} résulte du lemme suivant.

\begin{lemme} 
  \label{lem:cvabs}
Soit $e\in \ZZ$ et $P_1\subset P_2$ des sous-groupes paraboliques standard de $G$. Il existe un point $T\in \overline{a_B^+}$ ,tel que l'intégrale ci-dessous soit finie
\begin{equation}
  \label{eq:abs-cv}
 \int_{P_0(F)\back G(\AAA)^e} \sum_{\delta\in P_1(F)\back G(F)} \chi_T^{P_1,P_2}(\delta g) \, |\!\!\sum_{\{P\in \fc^{P_2}(P_1)\mid \delta\in P(F)\}}  (-1)^{\dim(a_P^G)} \tilde{K}^P_{D,X}(g)| \, dg<\infty.
\end{equation}
\end{lemme}

\begin{preuve}Traitons d'abord le cas $P_1=P_2$. Si $P_1\not=G$, la fonction $\chi_T^{P_1,P_2}$ est nulle et le résultat évident. Si $P_1=P_2=G$, l'intégrale se réécrit
$$\int_{P_0(F)\back G(\AAA)^e} F^G(g,T)  \sum_{ \delta\in G_X(F)\back P_0(F)}  \mathbf{1}_D(g^{-1}Xg) \, dg$$
ou encore 
$$\int_{G(F)\back G(\AAA)^e} F^G(g,T)  \sum_{ Y\in \oc}  \mathbf{1}_D(g^{-1}Yg) \, dg.$$
Cette dernière est prise sur un ensemble compact ; la convergence est donc évidente. On suppose désormais  $P_1\subsetneq P_2$. Effectuons quelques manipulations sur l'intégrale \eqref{eq:abs-cv}. Elle s'écrit encore
$$\int_{G(F)\back G(\AAA)^e} \sum_{\delta_0\in P_0(F)\back G(F)} \sum_{\delta\in P_1(F)\back G(F)} \chi_T^{P_1,P_2}(\delta \delta_0g) \, |\!\!\sum_{\{P\in \fc^{P_2}(P_1)\mid \delta\in P(F)\}}   (-1)^{\dim(a_P^G)} \tilde{K}^P_{D,X}(\delta_0g)| \, dg.$$
Par le changement de variables $\delta\mapsto \delta\delta_0$, elle est égale à 
$$\int_{G(F)\back G(\AAA)^e}  \sum_{\delta\in P_1(F)\back G(F)} \chi_T^{P_1,P_2}(\delta g) \, \sum_{\delta_0\in P_0(F)\back G(F)}|\!\!\sum_{\{P\in \fc^{P_2}(P_1)\mid \delta\delta_0^{-1}\in P(F)\}}   (-1)^{\dim(a_P^G)} \tilde{K}^P_{D,X}(\delta_0g)| \, dg.$$
Rappelons   $\tilde{K}^P_{D,X}$ est invariante à gauche par $P_0$ et  même  nulle  si $P_0\not\subset P$. Dans ce cas,  la condition  $\delta\delta_0^{-1}\in P(F)$ ne dépend que de la classe modulo $P_0(F)$ de $\delta_0$.

Par un nouveau changement de variables  $\delta_0\mapsto \delta_0\delta^{-1}$, on obtient
$$\int_{G(F)\back G(\AAA)^e}  \sum_{\delta\in P_1(F)\back G(F)} \chi_T^{P_1,P_2}(\delta g) \, \sum_{\delta_0\in P_0(F)\back G(F)}|\!\!\sum_{\{P\in \fc^{P_2}(P_1)\mid \delta_0^{-1}\in P(F)\}}  (-1)^{\dim(a_P^G)} \tilde{K}^P_{D,X}(\delta_0\delta g)| \, dg$$
$$=\int_{P_1(F)\back G(\AAA)^e} \chi_T^{P_1,P_2}(g) \, \sum_{\delta_0\in P_0(F)\back G(F)}|\!\!\sum_{\{P\in \fc^{P_2}(P_1)\mid \delta_0\in P(F)\}}  (-1)^{\dim(a_P^G)} \tilde{K}^P_{D,X}(\delta_0g)| \, dg.$$
À ce stade,  on utilise la décomposition d'Iwasawa pour voir que \eqref{eq:abs-cv} est égale à
$$=\int_{M_1(F)\back M_1(\AAA)\cap G(\AAA)^e} q^{-\bg2\rho_1,H_{P_1}(m)\bd} \chi_T^{P_1,P_2}(m) \times$$
$$\int_{N_1(F)\back N_1(\AAA)} \, \sum_{\delta\in P_0(F)\back G(F)}|\!\!\sum_{\{P\in \fc^{P_2}(P_1)\mid \delta\in P(F)\}}  (-1)^{\dim(a_P^G)} \tilde{K}^P_{D,X}(\delta nm)| \, dndm.$$ 

On utilise ensuite les deux lemmes suivants.

\begin{lemme}\label{lem:stepX1} Pour tous $\delta\in P_0(F)\back P(F)$, $n\in N_P(\AAA)$ et $m\in M_P(\AAA)$, on a 

  $$\tilde{K}^P_{D,X}(\delta nm)=q^{\dim(\ngo_P)(\deg(D)+1-g_C) +\bg2\rho_P,H_P(m) \bd }  \sum_{\{Y\in \mgo_P(F)\cap \delta^{-1}\ngo_0(F)\delta\mid I_P^G(X)=\oc\}} \mathbf{1}_D(m^{-1}Y m).$$
\end{lemme}

\begin{preuve}
On observera que les deux membres de l'égalité sont nuls si $P_0\not\subset P$. On suppose dans la suite $P_0\subset P$.
Pour tout $g\in G(\AAA)$,  on a 
$$\tilde{K}^P_{D,X}(g)=\sum_{\delta\in M_{X^P}(F)\back (P_0\cap M)(F)} \int_{\ngo(\AAA)} \mathbf{1}_D(g^{-1}(\delta^{-1}X^P\delta +U) g)\, dU.$$
  Pour $\delta\in M_{X^P}(F)\back (P_0\cap M)(F)$, les éléments $\delta^{-1}X^P\delta$ sont précisément les éléments $Y\in  \ngo_0(F)\cap \mgo(F)$ tels que $I_P^G(Y)=\oc$ (cf. section \ref{sec:induction}). Par conséquent, pour tout $\delta\in P(F)$, on a 
$$\tilde{K}^P_{D,X}(\delta g)= \sum_{\{Y\in \mathcal{N}^{M_P}\cap \delta^{-1}\ngo_0(F)\delta\mid I_P^G(X)=\oc\}} \int_{\ngo(\AAA)} \mathbf{1}_D(g^{-1}(Y+U)g)\, dU.$$
Le lemme s'en déduit aisément (cf. la preuve du lemme \ref{lem:step2}).
\end{preuve}

\begin{lemme}
  \label{lem:stepX2} Il existe un point $T_D$ tel que pour tous $T\in T_D+a_B^+$,  $\delta\in P(F)$,  $n \in  N_1(\AAA)$ et $m\in  M_1(\AAA)$ tels que $ \chi^{P_1,P_2}_T(m)\not=0$, on a
\begin{equation}
  \label{eq:equality}
  \tilde{K}^{P}_{D,X}(\delta nm)=q^{\dim(\ngo_2)(\deg(D)+1-g_C) +\bg2\rho_2,H_{P_2}(m) \bd }  \sum_{\{Y \in \pgo_1(F)\cap \mgo_2(F)\cap \delta^{-1}\ngo_0(F)\delta \mid  I_P^G(Y)=\oc\}} \mathbf{1}_D(m^{-1}Ym).
\end{equation}
\end{lemme}

\begin{preuve} Elle est semblable à celle du lemme \ref{lem:step4}, le rôle du lemme \ref{lem:step2} étant ici joué par le lemme \ref{lem:stepX1}.
\end{preuve}

Posons
$$\xi^{1,2}_{\oc,\delta}(Y)= \sum_{\{P\in \fc^{P_2}(P_1)\mid \delta\in P(F), I_P^G(Y)=\oc\}}(-1)^{\dim(a_P^G)}$$
et 
$$\Xi^{1,2}_{\oc}(Y)=\sum_{\{\delta \in P_0(F)\back G(F) \mid Y\in  \delta^{-1}\ngo_0(F)\delta \}} |\xi^{1,2}_{\oc,\delta}(Y)|.$$
Des lemmes \ref{lem:stepX1} et \ref{lem:stepX2}, on déduit la majoration suivante pour l'intégrale \eqref{eq:abs-cv} (à une puissance de $q$ près qu'il est inutile ici de préciser)
$$\int_{M_1(F)\back M_1(\AAA)\cap G(\AAA)^e} q^{-\bg2\rho_1^2,H_{P_1}(m)\bd} \chi_T^{P_1,P_2}(m) \sum_{Y \in \nc^{M_2}\cap \pgo_1(F)} \mathbf{1}_D(m^{-1}Ym)\cdot \Xi^{1,2}_{\oc}(Y)\,dm.$$

On utilise ensuite le lemme suivant dont la preuve (omise) est mot pour mot celle de la proposition \ref{prop:2ndemaj} (les notations sont celles du §\ref{S:reduc2}).

\begin{lemme}
  Pour tout diviseur $D$, il existe
\begin{itemize}
\item un diviseur $D'$ sur $C$
\item un point $T_D\in a_B^+$
\item pour tout $Q\in \fc^{P_1}(B)$ une constante $c_Q>0$
\end{itemize}
tels que pour tout $T\in T_D+a_B^+$, l'intégrale 

    \begin{equation}
      \label{eq:integ-amajII}
      \int_{M_1(F)\back M_1(\AAA)\cap G(\AAA)^e} q^{-\bg2\rho_1^2,H_1(m) \bd }  \chi^{P_1,P_2}_T(m) \sum_{Y\in \mathcal{N}^{M_2}\cap \pgo_1(F)} \mathbf{1}_D(m^{-1}Y m) \cdot  \Xi_{\oc}^{1,2}(Y)\, dm
    \end{equation}

se majore par la somme sur $Q\in \fc^{P_1}(B)$ et $a\in A^{1,2}_{Q,e}(T)$ de 
$$c_Q\cdot q^{-\bg2\rho_Q^2,H_Q(a) \bd } \sum_{Y \in \mathcal{N}^{M_2}\cap \pgo_1(F)} \mathbf{1}_{D'}(a^{-1}Y a)\cdot   \Xi_{\oc}^{P_1,P_2}(Y).$$
\end{lemme}

Comme dans le preuve de la proposition \ref{prop:reduc1} (cf. §\ref{S:preuvereduc1}), ce qui précède montre que le lemme \ref{lem:cvabs} est  une conséquence du lemme \ref{lem:dernierlemme} suivant. Cela conclut la preuve du lemme \ref{lem:cvabs}.
\end{preuve}

\begin{lemme}\label{lem:dernierlemme}
Soit $P_1$ un sous-groupe parabolique standard de $G$.   Pour tout  diviseur $D$ sur $C$ et tout $\al\in \Delta-\Delta^{P_1}$, il existe une constante $c>0$ telle que pour tout $a\in A^G(T_1)$, on a 

$$  q^{-\bg2\rho_1,H_1(a) \bd } \sum_{Y\in \nc^G \cap \pgo_1(F)} \mathbf{1}_D(a^{-1}Ya )\cdot  \Xi_{\oc}^{P_1,G}(Y)  \leq c\cdot  q^{-\bg\al,H_B(a)\bd}\cdot \sum_{Y\in \mathcal{N}^{M_{1}}} \mathbf{1}_D(a^{-1}Ya).$$
\end{lemme}

\begin{preuve}1
  Soit $Y\in \nc^G \cap \pgo_1(F)$. Rappelons qu'on note $P_0$ le sous-groupe parabolique standard qui vérifie $I_{P_0}^G(0)=\oc$. Son facteur de Levi est isomorphe à $GL(d)^r$.  On a
  \begin{equation}
    \label{eq:val-absolue}
    \Xi_{\oc}^{P_1,G}(Y)=\sum_{P_0'}  | \sum_{P} (-1)^{\dim(a_P^G)}  |
  \end{equation}
  où la première somme porte sur les sous-groupes paraboliques $P_0'$ tels que $P_0'$ est conjugué à $P_0$ et $ Y\in \ngo_0'(F)$. La second somme porte sur les sous-groupes paraboliques $P$ qui contiennent $P_1$ et $P_0'$ et tels que $I_P^G(Y)=\oc$.

Soit $P_0'$ un  sous-groupe parabolique  conjugué à $P_0$ qui contient $ Y$ dans son radical unipotent. Soit $M_0'$ un facteur de Levi de $P_0'$ et $A_0'$ son centre. Soit $\Delta_0'$ les \og racines simples\fg    de $A_0'$ dans $\ngo_0'$. Pour tout $\al$ soit $\ngo_\al$ l'espace de poids $\al$. On identifie en tant qu'espace vectoriel $\ngo_\al$ à $\mathfrak{gl}(d)$.  La somme $\oplus_{\al\in \Delta_0'} \ngo_\al$ est un supplémentaire dans $\ngo_0'$ de l'algèbre dérivée $[\ngo_0',\ngo_0']$. Soit $Q$ le  sous-groupe parabolique contenant $P_0'$ tel que  l'ensemble des racines simples de $A_0'$ dans $N_Q$ soit exactement l'ensemble des racines  $\al\in  \Delta_0'$ telles que la projection de $Y$ sur $\ngo_\al$ soit de rang $<d$. L'ensemble des sous-groupes paraboliques $P$ tels que $P_0'\subset P$ et $I_P^G(Y)=\oc$ admet alors une description simple : c'est exactement l'ensemble des sous-groupes paraboliques compris entre $P_0'$ et $Q$. 

Dans  \eqref{eq:val-absolue}, le terme associé à $P_0'$ est toujours nul sauf si $\tilde{P}_1$, le plus petit sous-groupe parabolique qui contient $P_0'$ et $P_1$, est égal à $Q$, auquel cas il vaut $1$. Par conséquent, on a

\begin{equation}
  \label{eq:derniere-maj}
  q^{-\bg2\rho_1,H_1(a) \bd } \sum_{Y\in \nc^G \cap \pgo_1(F)} \mathbf{1}_D(a^{-1}Ya )  \Xi_{\oc}^{P_1,G}(Y)\leq \sum_{Q} \sum_{P_0' }   q^{-\bg2\rho_1,H_1(a) \bd }  \sum_Y \mathbf{1}_D(a^{-1}Ya )
\end{equation}
où les sommes portent sur 
\begin{itemize}
\item les sous-groupes paraboliques $Q$ contenant $P_1$,
\item les sous-groupe paraboliques $P_0'\subset Q$ tels que $P_0'$ est conjugué à $P_0$ sous $G$ et le plus petit sous-groupe parabolique contenant $P_1$ et $P_0'$ est $Q$,
\item les éléments $Y\in \pgo_1(F)\cap\ngo_0'(F)$ tels que la projection de $Y$ sur $\mgo_Q$ appartient à $I_{P_0'}^Q(0)$ et celle de $Y$ sur $\ngo_Q$ appartient à l'ensemble des éléments de $\ngo_Q(F)$ dont l'image dans $\ngo_0'/ [\ngo_0',\ngo_0']$ (identifié à $\mathfrak{gl}(d)^r$) est formé de \og blocs\fg{} de rang $<d$. 
\end{itemize}

Il reste à majorer le terme de droite de \eqref{eq:derniere-maj}. Pour cela, on fixe un sous-groupe parabolique $Q$ comme dans ce terme (il n'y en a qu'un nombre fini). 
Soit $\al\in \Delta-\Delta^{P_1}$. Cette racine détermine un sous-groupe parabolique maximal $R$ qui contient $P_1$.

Supposons tout d'abord $Q\subset R$ et donc $\ngo_R\subset \ngo_Q$. Pour tout $P_0'\subset Q$ comme ci-dessus, soit $\vc_{P_0'}\subset \ngo_R$ le fermé formé des $Z\in \ngo_R$ dont l'image  dans $\ngo_0'/ [\ngo_0',\ngo_0']$ (identifié à $\mathfrak{gl}(d)^r$) est formé de \og blocs\fg{} de rang $<d$. Il s'agit d'un fermé défini par des polynômes dont le degré est borné. Dans ce cas, on a la majoration suivante ($P_0'$ et $Y$ sont comme ci-dessus)

\begin{equation}
  \label{eq:lastbut}
   \sum_{P_0' }   q^{-\bg2\rho_1,H_1(a) \bd }  \sum_Y \mathbf{1}_D(a^{-1}Ya )\leq \sum_{P_0' }   q^{-\bg2\rho_1,H_1(a) \bd }  \sum_{Y^Q, Y_Q^R, Y_R}  \mathbf{1}_D(a^{-1}(Y^Q+Y_Q^R+Y_R)a )
\end{equation}

où, dans le membre de droite, $P_0'$ est comme ci-dessus, $Y^Q$ parcourt les éléments de 
$$\mgo_Q(F)\cap \pgo_1(F)\cap \ngo_0'(F)\cap I_{P_0'}^Q(0),$$
$Y_Q^R$ parcourt $\ngo_Q(F)\cap \mgo_R(F)$, $Y_R$ parcourt $\vc_{P_0'}$. Il résulte de la proposition \ref{prop:maj-var} qu'il existe une constante $c>0$ telle que pour $a\in A^G(T_1)$, 
$$  q^{-\bg2\rho_R,H_R(a) \bd }\sum_{Y_R\in \vc_{P_0'}}  \mathbf{1}_D(a^{-1}(Y_R)a )\leq c \cdot   q^{-\bg\al,H_B(a)\bd}. $$

Observons ensuite que l'orbite $I_{P_0'}^Q(0)$ ne dépend pas du choix de $P_0'\subset Q$. C'est d'ailleurs l'orbite qu'on a noté $\oc^Q$ précédemment. De plus, si $Y\in \mgo_Q(F)\cap  I_{P_0'}^Q(0)$, il existe un unique  $P_0'\subset Q$ tel que $Y\in \ngo_0'$. Il s'ensuit qu'on peut obtient le majorant suivant pour  \eqref{eq:lastbut} 
$$c\cdot  q^{-\bg\al,H_B(a)\bd}\cdot \big[ q^{-\bg2\rho_1^R,H_1(a) \bd } \sum_{Y\in \pgo_1(F)\cap \mgo_R(F)}  \mathbf{1}_D(a^{-1} Ya ) \big].$$
Le terme entre crochets est un majorant grossier mais qui est borné pour $a\in A^G(T_1)$ (cf. lemme \ref{lem:maj}). Cela conclut lorsque $Q\subset R$.

Supposons ensuite qu'on a $Q\not\subset R$. Alors $R\cap M_Q\subset M_Q$ est un sous-groupe parabolique maximal de $M_Q$. Toujours par le lemme  \ref{lem:maj}, pour $a\in  A^G(T_1)$, la somme
$$q^{-\bg2\rho_Q,H_Q(a) \bd } \sum_{Y\in \ngo_Q(F)}  \mathbf{1}_D(a^{-1} Ya )$$
est bornée. Quitte à remplacer $G$ par $M_Q$ et $P_1$ par $M_Q\cap P_1$, on voit qu'il suffit de traiter le cas $Q=G$. Il s'agit donc de majorer la quantité suivante
\begin{equation}
  \label{eq:least}
  \sum_{P_0'} q^{-\bg2\rho_1,H_1(a) \bd }\sum_{Y\in \oc\cap \pgo_1(F)\cap \ngo_0'(F)} \mathbf{1}_D(a^{-1} Ya )
\end{equation}
où la première somme est prise sur les sous-groupes paraboliques de $G$ conjugués à $P_0$ tel que le plus petit sous-groupe parabolique qui contienne à la fois $P_0'$ et $P_1$ est $G$ lui-même. Pour tout $Y\in \pgo_1(F)$, on a  $Y\in \overline{I_R^G(Y)}$. Si, de plus,  $Y \in \oc\cap \pgo_1(F)\cap \ngo_0'(F)$ où $P_0'$ est comme ci-dessus alors $Y$ appartient au fermé de $\overline{I_R^G(Y)}$ complémentaire de $I_R^G(Y)$. En effet, supposons le contraire : si $Y\in I_R^G(Y)$, on a  $\overline{I_R^G(Y)}=\oc$ et donc la projection $Y^R$ de $Y$ sur $\mgo_R$ (relativement à la décomposition $\rgo=\mgo_R\oplus\ngo_R$) appartient à l'orbite $\oc^R=I_{P_0''}^R(0)$ pour un certain sous-groupe parabolique $P_0''\subset R$ qui est $G$-conjugué à $P_0$. On peut et on va supposer que $Y\in \ngo_0''(F)$. Un tel sous-groupe parabolique $P_0''$ est nécessairement le stabilisateur des drapeaux des images itérées de $Y$. On a donc $P_0'=P_0''$. En particulier, on a $P_0'\subset R$. Maintenant $R$ contient $P_0'$ et $P_1$ donc par hypothèse sur $P_0'$ on a $R=G$ ce qui contredit l'hypothèse que $R$ est maximal. 

Comme on l'a déjà utilisé, pour tout $Y\in \oc$, il existe un unique sous-groupe parabolique $P_0'$ conjugué à $P_0$ tel que $Y\in \ngo_0'$. On peut donc majorer \eqref{eq:least}  par

\begin{equation}
  \label{eq:least2}
q^{-\bg2\rho_1,H_1(a) \bd }\sum_{\{Y\in \oc\cap \pgo_1(F) \mid Y\notin I_R^G(Y) } \mathbf{1}_D(a^{-1} Ya ).
\end{equation}
On obtient alors le majorant voulu par le même raisonnement que dans la preuve de la proposition \ref{prop:reduc4} (plus précisément dans la majoration de \eqref{eq:XUV2}).
\end{preuve}

\end{paragr}

\bibliographystyle{plain}
\bibliography{bibliog}

\begin{flushleft}
Pierre-Henri Chaudouard \\
Université Paris Diderot (Paris 7)\\ 
 Institut de Mathématiques de Jussieu-Paris Rive Gauche \\ 
 UMR 7586 \\
 Bâtiment Sophie Germain \\
 Case 7012 \\
 F-75205 PARIS Cedex 13 \\
 France 
\medskip

Gérard Laumon \\
CNRS et Universit\'{e} Paris-Sud \\
 UMR 8628 \\
 Math\'{e}matique, B\^{a}timent 425 \\
F-91405 Orsay Cedex \\
France \\
\bigskip

Adresses électroniques :\\
Chaudouard@math.jussieu.fr \\
Gerard.Laumon@math.u-psud.fr\\

\end{flushleft}

\end{document}